\documentclass{article}
\usepackage{amsmath}[1996/11/01]
\usepackage{amssymb,amsthm,amsxtra}

\def\[#1\]{\begin{equation}#1\end{equation}}
\makeatletter
\def\beq{%
   \relax\ifmmode
      \@badmath
   \else
      \ifvmode
         \nointerlineskip
         \makebox[.6\linewidth]%
      \fi
      $$
   \fi
}
\def\eeq{%
   \relax\ifmmode
      \ifinner
         \@badmath
      \else
         $$
      \fi
   \else
      \@badmath
   \fi
   \ignorespaces
}

\def\enddisplaymath{\eeq\global\@ignoretrue}
\makeatother

\newtheorem{thm}{Theorem}
\newtheorem{cor}[thm]{Corollary}
\newtheorem{lem}[thm]{Lemma}
\newtheorem{prop}[thm]{Proposition}
\newtheorem{conj}{Conjecture}
\theoremstyle{definition}
\newtheorem{defn}{Definition}
\theoremstyle{remark}
\newtheorem*{rem}{Remark}
\newtheorem{rems}{Remark}[thm]

\newtheorem{drems}{Remark}[defn]

\numberwithin{equation}{section}
\numberwithin{thm}{section}

\newcommand{\Q}{\mathbb Q}
\newcommand{\C}{\mathbb C}

\newcommand{\Z}{\mathbb Z}

\renewcommand{\P}{\mathbb P}
\newcommand{\A}{\mathbb A}

\DeclareMathOperator\SL{SL}
\DeclareMathOperator\Tr{Tr}
\DeclareMathOperator\Pic{Pic}
\newcommand{\la}{\langle}
\newcommand{\ra}{\rangle}
\newcommand{\binomE}{\genfrac(){0pt}{}}
\newcommand{\obinomE}{\genfrac\la\ra{0pt}{}}
\newcommand{\binomQ}{\genfrac[]{0pt}{}}
\newcommand{\binomI}{\genfrac\{\}{0pt}{}}
\newcommand{\II}{\mathord{I\!I}}
\newcommand{\III}{\mathord{I\!I\!I}}

\begin{document}

\title{$BC_n$-symmetric abelian functions}
\author{Eric M. Rains\footnote{Department of Mathematics, University of California, Davis}}

\date{January 9, 2006}
\maketitle

\begin{abstract}
We construct a family of $BC_n$-symmetric biorthogonal abelian functions
generalizing Koornwinder's orthogonal polynomials, and prove a number of
their properties, most notably analogues of Macdonald's conjectures.  The
construction is based on a direct construction for a special case
generalizing Okounkov's interpolation polynomials.  We show that these
interpolation functions satisfy a collection of generalized hypergeometric
identities, including new multivariate elliptic analogues of Jackson's
summation and Bailey's transformation.

\end{abstract}

\tableofcontents

\section{Introduction}

In \cite{bcpoly}, we considered two families of $BC_n$-symmetric Laurent
polynomials: the multivariate orthogonal polynomials introduced by
Koornwinder \cite{KoornwinderTH:1992}, and the ``interpolation Macdonald
polynomials'' of Okounkov \cite{OkounkovA:1998a}, related by Okounkov's
``binomial formula''.  One of the results of that paper was a multivariate
analogue of Jackson's summation, which had originally been conjectured by
Warnaar \cite{WarnaarSO:2002}.  In fact, that summation was a limiting case
of Warnaar's conjecture, an identity of elliptic hypergeometric series.
This suggested that the theory of \cite{bcpoly} should be extensible to the
elliptic case; this extension is the topic of the present work.

The main difficulty in constructing this extension is the result of
Okounkov \cite{OkounkovA:1998c} which roughly speaking said that the
interpolation Macdonald polynomials could not be generalized; more
precisely, he defined a general notion of interpolation polynomial, and
showed that the Macdonald-type polynomials form a component of the
corresponding moduli space.  Thus in order to extend our results to the
elliptic level, it would be necessary to generalize this notion of
interpolation polynomial.

Another indication that an elliptic extension should exist was the
existence of elliptic analogues
\cite{SpiridonovVP/ZhedanovAS:2000b,SpiridonovVP:2003}, both discrete and
continuous, of the Askey-Wilson orthogonal polynomials
\cite{AskeyR/WilsonJ:1985}.  The difference from the Askey-Wilson theory
(already seen at the trigonometric level in work of Rahman
\cite{RahmanM:1986} for continuous biorthogonality, and later work by
Wilson \cite{WilsonJA:1991} for discrete biorthogonality) was two-fold:
first, rather than being orthogonal, the functions were merely {\em
  biorthogonal} , and second, the constraint of being monic of specified
degree became a constraint on allowable poles.  Extending this prescription
directly to the multivariate level ran into the difficulty that divisors in
multiple variables are much more complicated than in the univariate case.
However, if we first clear the poles, we find that the constraint on poles
becomes a requirement that the function vanish at appropriate points.

We thus arrive at the notion of ``balanced'' interpolation polynomial; this
differs from Okounkov's notion in that, rather than have one collection of
vanishing conditions together with an assumption of triangularity on
monomials, balanced interpolation polynomials satisfy two complementary
collections of vanishing conditions.  (This definition was independently
discovered by Coskun and Gustafson \cite{CoskunH/GustafsonRA:2003}, who also
obtained many of the results of Section \ref{sec:hyperg1}, in particular
Theorems \ref{thm:theta_bde} and \ref{thm:theta_bailey} below.)  This
allowed the affine symmetry of ordinary interpolation polynomials to be
extended to a projective symmetry, and as a result allowed the desired
generalization to elliptic interpolation polynomials.  Moreover, two new
special cases arise: first, by matching up the vanishing conditions, one
obtains a collection of delta functions; more subtly, another special case
appears in which the interpolation polynomials factor completely.  These
special cases then play the role that the Macdonald polynomials did in
\cite{bcpoly}, namely that of instances that can be solved directly, thus
giving boundary conditions for recurrences.

In particular, we have been able to generalize many of the results of
\cite{bcpoly}, especially concerning the Koornwinder polynomials, which we
generalize to a family of biorthogonal abelian functions.  A closely
related family was studied in \cite{xforms} using a certain contour
integral identity and a related integral operator; in particular, it was
shown there that the biorthogonal functions satisfy an analogue of
Macdonald's normalization and evaluation conjectures.  The analogue of the
remaining conjecture (referred to here as ``evaluation symmetry'') requires
a more in depth understanding of interpolation polynomials, and will be
proved in Theorem \ref{thm:biorth_ev_symm} below.  In contrast to
\cite{xforms}, our present approach is essentially algebraic in nature;
although to begin with we use the theory of theta functions over $\C$, we
will eventually see that all of our results can be given purely algebraic
proofs.

We begin in Section 2 by giving the definition and basic properties of
balanced interpolation polynomials, in particular considering four of the
main cases satisfying extra vanishing properties (the ``perfect'' case);
the fifth (``elliptic'') case is the subject of the remainder of the paper.

Section 3 begins our treatment of the elliptic case by defining a family of
$BC_n$-symmetric theta functions, parametrized by partitions contained in a
rectangle, defined via vanishing conditions.  As in \cite{bcpoly}, the key
property of these functions is a difference equation, from which the extra
vanishing property follows.  We also define a corresponding family of
abelian functions, this time indexed by all partitions of length at most
$n$.

Evaluating an interpolation abelian function at a partition gives a
generalized binomial coefficient; we study these coefficients in Section 4.
In particular, we obtain a pair of identities generalizing Jackson's
summation and Bailey's transformation respectively (see
\cite{FrenkelIB/TuraevVG:1997} for the univariate elliptic case); each
identity involves a sum over partitions contained in a skew Young diagram
of a product of binomial coefficients.  Using these, we derive some
unexpected symmetries of binomial coefficients, relating a coefficient to
coefficients with conjugated or complemented partitions.  This then enables
us to prove a number of properties of interpolation functions, including a
connection coefficient formula, a branching rule, a Pieri identity, and a
Cauchy identity.

Section 5 introduces the biorthogonal functions, defined via an expansion
in interpolation functions, and shows how they relate to the functions
considered in \cite{xforms}.  Using the interpolation function identities,
we prove a number of identities for the biorthogonal functions; in addition
to the analogue of Macdonald's conjectures for Koornwinder polynomials, we
obtain connection coefficients, a quasi-branching rule, a quasi-Pieri
identity, and a Cauchy identity.

In Section 6, we depart from the analytic approach via theta functions, and
show that, in fact, the interpolation and biorthogonal functions can be
defined purely algebraically (and thus, over the complex numbers, are
invariant under transformations of the modular parameter).  We sketch a
purely algebraic derivation of the identities of Sections 4 and 5; the key
step being to construct the algebraic analogue of the difference operator
of Section 3.

In Section 7, we use our results on interpolation functions to define a
family of perfect bigrids based on elliptic curves.  Since this family is
an open subset of the space of perfect bigrids, we in particular find that
we cannot add any parameters to our elliptic theory, without further
extending the notion of interpolation polynomial.  We then consider how the
bigrids in this family can degenerate, obtaining four general classes of
degeneration (corresponding to the Kodaira symbols $I_1$, $I_2$, $\II$, and
$\III$; these in turn split somewhat further over a
non-algebraically-closed field), including ordinary hypergeometric series
($\II$ and some instances of $\III$) and basic hypergeometric series ($I_1$
and some instances of $I_2$) as well as some apparently unstudied cases.

Section 8 shows how the present results relate to the results of
\cite{bcpoly}.  The identities of that paper involved three different
functions of skew Young diagrams: binomial coefficients, ``inverse''
binomial coefficients, and principal specializations of skew Macdonald
polynomials.  We show that each of these can be obtained as a limit of
elliptic binomial coefficients, and thus in particular prove that the
Koornwinder polynomials are a limiting case of the biorthogonal functions.

Finally, in Section 9, we consider some open problems relating to
interpolation and biorthogonal functions.

We would like to thank R. Gustafson, A. Okounkov, H. Rosengren, S. Sahi,
and V. Spiridonov for enlightening conversations related to this work.

\vfil\break
\noindent{\bf Notation}

As in \cite{bcpoly}, we use the notations of \cite{MacdonaldIG:1995} for
partitions, with the following additions.  If $\lambda\subset m^n$, then
$m^n-\lambda$ denotes the ``complementary'' partition defined by
\[
(m^n-\lambda)_i = m-\lambda_{n+1-i},\quad 1\le i\le n
\]
In addition, if $\ell(\lambda)\le n$, we define $m^n+\lambda$ to be the
partition
\[
(m^n+\lambda)_i = m+\lambda_i,\quad 1\le i\le n.
\]
Similarly, if $\lambda_1\le m$, we define $m^n\cdot\lambda$ to be the
partition such that
\[
(m^n\cdot\lambda)_i = \begin{cases}m&1\le i\le n\\
     \lambda_{i-n} & i>n.\end{cases}
\]
(Note that this is specifically not defined for $\lambda_1>m$; it is a
concatenation operation, not to be confused with a union operator.)
We also define a relation $\prec_m$ on partitions such that
$\kappa\prec_m\lambda$ if and only if
\[
\kappa \subset \lambda\subset m^n+\kappa
\]
for sufficiently large $n$; similarly, $\kappa\prec'_m\lambda$ iff
$\kappa'\prec_m\lambda'$.  (If $m=1$ in either case, the subscript will be
omitted.)

We will also need the following ``elliptic'' analogues of $q$-symbols and
generalizations.  Let $p$ be a complex number such that $0<|p|<1$.
Then we define
\begin{align}
\theta(x;p)&:= \prod_{0\le k} (1-p^k x)(1-p^{k+1}/x)\\
\theta(x;q;p)_m &:= \prod_{0\le j<m} \theta(q^j x;p)
\end{align}
In each case (and similarly for the $C$ symbols (but {\em not} the $\Delta$
symbols) below), we follow the standard convention that the presence of
multiple arguments before the (first) semicolon indicates a product; thus
for instance $\theta(a x^{\pm 1};p) = \theta(ax;p)\theta(a/x;p)$.  We also
need elliptic analogues of the $C$ symbols of \cite{bcpoly}:
\begin{align}
C^0_\lambda(x;q,t;p)
&:=
\prod_{(i,j)\in \lambda} \theta(q^{j-1}t^{1-i} x;p)\\
C^-_\lambda(x;q,t;p)
&:=
\prod_{(i,j)\in\lambda} \theta(q^{\lambda_i-j}t^{\lambda'_j-i} x;p)\\
C^+_\lambda(x;q,t;p)
&:=
\prod_{(i,j)\in\lambda} \theta(q^{\lambda_i+j-1}t^{2-\lambda'_j-i} x;p).
\end{align}
Note that in the limit $p=0$ we recover the symbols of \cite{bcpoly};
similarly, the analogous symbols of \cite{xforms} can each be expressed as
products of two of our present symbols (with $q$ and $p$ switched in one
symbol).  The transformations of \cite{bcpoly} carry over to the present
case, via the same arguments.

Two combinations of $C$ symbols are of particular importance, and are thus
given their own notations.  For any nonzero complex number $a$ and any
finite sequence of nonzero complex numbers $\dots b_i\dots$, we define two
$\Delta$ symbols as follows.
\begin{align}
\Delta^{0}_{\lambda}(a|\dots b_i\dots;q,t;p)
&:=
\prod_i
\frac{C^0_{\lambda}(b_i;q,t;p)}
     {C^0_{\lambda}(pqa/b_i;q,t;p)}\\
\Delta_{\lambda}(a|\dots b_i\dots;q,t;p)
&:=
\Delta^{0}_{\lambda}(a|\dots b_i\dots;q,t;p)
\frac{C^0_{2\lambda^2}(pqa;q,t;p)}
     {C^-_{\lambda}(pq,t;q,t;p)C^+_{\lambda}(a,pqa/t;q,t;p)}.
\end{align}
(Here, as in \cite{bcpoly}, $2\lambda^2$ denotes the partition such that
$(2\lambda^2)_i = 2(\lambda_{\lceil i/2\rceil}).$) Note that although
$\Delta^0$ essentially follows the standard multiple argument convention,
$\Delta$ very much does not.  We note the following transformations:
\begin{align}
\Delta_{\lambda'}(a|\dots b_i\dots;1/t,1/q;p)
&=
\Delta_\lambda(a/qt|\dots b_i\dots;q,t;p)
\\
\frac{\Delta_{m^n-\lambda}(a|\dots b_i\dots;q,t;p)}
     {\Delta_{m^n}(a|\dots b_i\dots;q,t;p)}
&=
\Delta_{\lambda}(\frac{t^{2n-2}}{q^{2m}a}|\dots \frac{t^{n-1}b_i}{q^m
  a}\dots,t^n,q^{-m},pqt^{n-1},pq/q^mt;q,t;p)
\\
\frac{\Delta_{m^n+\lambda}(a|\dots b_i\dots;q,t;p)}
     {\Delta_{m^n}(a|\dots b_i\dots;q,t;p)}
&=
\Delta_{\lambda}(q^{2m}a|\dots q^m b_i\dots,t^n,pq t^{n-1},
 q^m a/t^{n-1}, pq q^m a/t^n;q,t;p)
\\
\frac{\Delta_{m^n\cdot\lambda}(a|\dots b_i\dots;q,t;p)}
{\Delta_{m^n}(a|\dots b_i\dots;q,t;p)}
&=
\Delta_{\lambda}(a/t^{2n}|\dots b_i/t^n\dots,q^{-m},pq/q^mt,q^m
a/t^{n-1},pq q^m a/t^n;q,t;p),
\end{align}
with corresponding transformations for $\Delta^0$; we also note that
\[
\Delta_{m^n}(a|\dots b_i\dots;q,t;p)
=
\lim_{x\to 1}
\Delta^0_{m^n}(a|\dots b_i\dots,
  q^m pqa/x, pqa/xt^n,pqax,pq t^{n-1}x/q^m;q,t;p)
\]
\medskip
A {\em meromorphic ($p$-)theta function} is a meromorphic function $f$ on
$\C^*$ (the multiplicative group of $\C$; i.e., $f$ is meromorphic away
from 0 and $\infty$) such that
\[
f(x) = a x^n f(px)
\]
for suitable $a\in \C^*$, $n\in \Z$; without the adjective
``meromorphic'', theta functions are assumed holomorphic.  Thus in
particular the function $\theta(x;p)$ is a theta function, satisfying
\[
\theta(x;p) = -x\theta(1/x;p) = -x\theta(px;p).
\]
An {\em elliptic function} is a meromorphic function on $\C^*$ invariant
under multiplication by $p$; in other words, a meromorphic function on the
quotient group $\C^*/\la p\ra$ (which is a compact 1-dimensional complex
Lie group, or in other words an elliptic curve).

We also have the following classes of multivariate functions.

\begin{defn}
A {\em meromorphic $BC_n$-symmetric theta function of degree $m$} is a
meromorphic function $f$ on $(\C^*)^n$ such that
\begin{itemize}
\item[] $f(x_1,\dots, x_n)$ is invariant under permutations of its
  arguments.
\item[] $f(x_1,\dots, x_n)$ is invariant under $x_i\mapsto 1/x_i$ for each
  $i$.
\item[] $f(px_1,x_2,\dots, x_n)=(1/px_1^2)^m f(x_1,x_2,\dots, x_n)$.
\end{itemize}
A {\em $BC_n$-symmetric ($p$-)abelian function} is a meromorphic
$BC_n$-symmetric theta function of degree $0$.
\end{defn}

\begin{rems}
Again, $BC_n$-symmetric theta functions are assumed holomorphic unless
specifically labeled ``meromorphic''.
\end{rems}

\begin{rems}
It is not entirely clear to which root system these functions are truly
attached.  Certainly, the biorthogonal functions considered below
generalize the Koornwinder polynomials, traditionally associated to the
root system $BC_n$, but there is also good reason to associate this class
of theta functions to the affine root system of type $A^{(2)}_n$ (see, for
instance, the work of Looijenga \cite{LooijengaE:1976} and Saito
\cite{SaitoK:1990} on the invariant theory of abelian varieties associated
to (extended) affine root systems.)  In the absence of a general theory, we
have chosen the more familiar label.
\end{rems}

The canonical example of this is the function $\prod_{1\le i\le n}\theta(u
x_i^{\pm 1};p)$ for $u\in \C^*$; indeed, the graded algebra of
$BC_n$-symmetric theta functions is generated by such functions.

We define a difference operator on (meromorphic) $BC_n$-symmetric theta
functions that will play a crucial role in the sequel; this is the
analogue of the difference operator of \cite{bcpoly}, and also appeared
prominently in \cite{xforms}.

\begin{defn}
Let $n$ be a nonnegative integer, and let $q,t,a,b,c,d\in \C^*$ be arbitrary
parameters.  Define a difference operator $D^{(n)}(a,b,c,d;q,t;p)$ acting
on $BC_n$-symmetric meromorphic functions on $(\C^*)^n$ by:
\begin{align}
(D^{(n)}(a,b,c,d;&q,t;p)f)(x_1,\dots x_n)\notag
\\
&=
\sum_{\sigma\in\{\pm 1\}^n}
\prod_{1\le i\le n}
\frac{\theta(a x_i^{\sigma_i},b x_i^{\sigma_i},c x_i^{\sigma_i},d x_i^{\sigma_i};p)}{\theta(x_i^{2\sigma_i};p)}
\prod_{1\le i<j\le n}
\frac{\theta(t x_i^{\sigma_i} x_j^{\sigma_j};p)}{\theta(x_i^{\sigma_i} x_j^{\sigma_j};p)}
f(\dots q^{\sigma_i/2} x_i\dots).
\end{align}
\end{defn}

\begin{prop}
For each nonnegative integer $m$, if $q^m t^{n-1} a b c d = p$, then
$D^{(n)}(a,b,c,d;q,t;p)$ maps the space of $BC_n$-symmetric theta
functions of degree $m$ into itself.
\end{prop}

\begin{proof}
Every term in the sum is a theta function in each variable with the same
functional equation, and thus the sum is itself a theta function in each
variable.  It remains only to show that it is holomorphic, which follows
from the usual symmetry argument.
\end{proof}

\begin{rem}
For $n=1$, the above operators span a four-dimensional space of difference
operators acting on $BC_1$-symmetric theta functions of degree $m$.  This

its action on theta functions \cite{SklyaninEK:1982,SklyaninEK:1983}.
Moreover, the relations of the Sklyanin algebra are simply given by the
quasi-commutation relation of Corollary \ref{cor:quasicommutation} below.
\end{rem}

The case $m=0$ is worth noting:

\begin{cor}
If $t^{n-1}abcd=p$, then
\[
\sum_{\sigma\in \{\pm 1\}^n}
\prod_{1\le i\le n}
\frac{\theta(a x_i^{\sigma_i},b x_i^{\sigma_i},c x_i^{\sigma_i},d x_i^{\sigma_i};p)}{\theta(x_i^{2\sigma_i};p)}
\prod_{1\le i<j\le n}
\frac{\theta(t x_i^{\sigma_i} x_j^{\sigma_j};p)}{\theta(x_i^{\sigma_i} x_j^{\sigma_j};p)}
=
\prod_{1\le i\le n}
\theta(ab t^{n-i},ac t^{n-i},ad t^{n-i};p)
\]
\end{cor}

\begin{proof}
By the proposition, the sum is independent of $x_i$; setting $x_i=a
t^{n-i}$ gives the desired result.
\end{proof}

\begin{rem}
This identity is, in fact, equivalent to the identity used in
\cite{RosengrenH:2001,RosengrenH:2004} to prove $BC_n$-type elliptic
hypergeometric identities, a special case of Theorem 5.1 of
\cite{WarnaarSO:2002}.  In particular, there is a determinantal proof of
this identity, which turns out to generalize to a determinantal formula for
the difference operator, although we will not consider that formula here.
\end{rem}

\section{Balanced interpolation polynomials}

For this section, we fix $R$ to be an arbitrary commutative ring.  A point
on the projective line $\P^1(R)$ is represented by relatively prime
homogeneous coordinates $(z,w)$ (i.e., the ideal generated by $z$ and $w$
contains 1), modulo multiplication by units; a polynomial of degree $m$ on
$\P^1(R)$ (a section of ${\cal O}(m)$) then corresponds to a homogeneous
polynomial of degree $m$ in these coordinates.  Similarly, given an
$n$-tuple of nonnegative integers, one has a notion of a polynomial on
$\P^1(R)^n$ of degree $(m_1,m_2,\dots, m_n)$.  The value of such a
polynomial at a given $n$-tuple of points on $\P^1(R)$ is only defined
modulo units, unless specific homogeneous coordinates are chosen for each
point; however, this freedom leaves invariant the zero-locus of the
polynomial as well as the ratio of two such polynomials of the same degree.
We will thus feel free to write
\[
p(x_1,x_2,\dots, x_n)
\]
for the value of $p$ at the points $x_1,\dots, x_n\in \P^1(R)$, with this
understanding.

Of particular importance is the polynomial $x\cdot y$ of degree $(1,1)$
defined by
\[
(z_1,w_1)\cdot (z_2,w_2) = z_1w_2-w_1z_2.
\]
Note that $x\cdot y$ is invariant under the action of $\SL_2(R)$; since
$(1,0)\cdot (z,w)=w$, we find that $x\cdot y=0$ if and only if $x$ and $y$
represent the same point.  Similarly, $x\cdot y\in R^*$ (the unit group of
$R$) if and only if they represent distinct points in $\P^1(R_0)$ for all
homomorphic images $R_0$ of $R$.

\begin{defn}
A symmetric polynomial on $\P^1(R)$ of degree $m$ in $n$ variables is a
polynomial on $\P^1(R)^n$ of degree $(m,m,\dots, m)$ invariant
under permutations of the $n$ arguments.
\end{defn}

The space of symmetric polynomials of degree $m$ in $n$ variables (denoted
$\Lambda^m_n$) is a free module of rank $\binom{m+n}{n}$; this also counts
the number of partitions $\lambda\subset m^n$.  An explicit basis of this
module is, for instance, given by the homogeneous monomials of degree $m$
in the $n+1$ multilinear polynomials
\[
e_i(z_1,w_1,\dots, z_n,w_n)
=
\sum_{\substack{I\subset\{1,2,\dots, n\}\\ |I|=i}}
\prod_{j\in I} z_j
\prod_{j\notin I} w_j,
\]
$0\le i\le n$.  A polynomial in $\Lambda^m_n$ will be said to be {\em
  primitive} if the ideal generated by its coefficients contains $1$; note
that it suffices to consider its coefficients in the $e_i$ basis.

The following decomposition will prove useful in the sequel.

\begin{lem}\label{lem:sf_decomp}
Fix integers $m$, $n\ge 0$, and a point $x_0\in \P^1(R)$ (with chosen
homogeneous coordinates).  Then we have the following short exact sequence:
\[
0\to \Lambda^{m-1}_n\mathop{\to}\limits^f \Lambda^m_n\mathop{\to}\limits^g \Lambda^m_{n-1}\to 0,
\]
in which the map $f$ is defined by
\[
f(p)(x_1,\dots, x_n) = p(x_1,\dots, x_n)\prod_{1\le i\le n} (x_0\cdot
x_i)
\]
and the map $g$ is defined by
\[
g(p)(x_1,\dots, x_{n-1})=p(x_0,x_1,\dots, x_{n-1}).
\]
Moreover, this exact sequence splits: there exists a map
$h:\Lambda^m_{n-1}\to \Lambda^m_n$ such that $g\circ h=1$.
\end{lem}

\begin{proof}
The claim is clearly invariant under $\SL_2(R)$, so we may assume that
$x_0=(0,1)$.  But then $f(p)=e_n p$, $g(e_i)=e_i$ for $i<n$, and
$g(e_n)=0$; that the sequence is exact follows immediately.  The map $h$ is
then given by $h(e_i)=e_i$.
\end{proof}

\begin{defn}
Let $m$, $n>0$ be positive integers.  A ($R$-valued) {\em extended bigrid
of shape $m^n$} is a function
\[
\gamma:\{0,1\}\times \{1,\dots, n\}\times \{0,1,\dots, m\}\to \P^1(R).
\]
A {\it bigrid of shape $m^n$} is the restriction of an extended bigrid
of shape $m^n$ to the subset of $(\alpha,i,j)$ such that if $\alpha=0$ then
$j<m$, and if $\alpha=1$, then $j>0$.  If the extended bigrid $\gamma^+$
restricts to $\gamma$, then $\gamma^+$ is said to be an {\it extension} of
$\gamma$; in considering the extensions of a given bigrid, we allow
the extended values to lie in any ring containing $R$.
\end{defn}

\begin{defn}
Let $\gamma$ be a bigrid of shape $m^n$, and let $\lambda\subset m^n$
be a partition.  A {\em balanced interpolation polynomial} of index
$\lambda$ for $\gamma$ is a symmetric polynomial $P^*_\lambda(;\gamma)$ of
degree $m$ on $\P^1(R)^n$ such that for all extensions $\gamma^+$ of
$\gamma$, the following vanishing conditions hold: (1) For all partitions
$\mu\subset m^n$, $\lambda\not\subset\mu$,
\[
P^*_\lambda(\dots\gamma^+(0,i,\mu_i)\dots;\gamma)=0,
\]
and (2) for all partitions $\mu\subset m^n$, $\mu\not\subset\lambda$,
\[
P^*_\lambda(\dots\gamma^+(1,i,\mu_i)\dots;\gamma)=0.
\]
The bigrid $\gamma$ is {\em perfect} if for all $\lambda\subset m^n$, there
exists a primitive balanced interpolation polynomial of index $\lambda$ for
$\gamma$.
\end{defn}

In the sequel, we will omit the word ``balanced'' except when necessary to
avoid confusion with Okounkov's interpolation polynomials, which we will
refer to as ``ordinary'' interpolation polynomials throughout.  (The $*$ to
denote interpolation polynomials is inherited from shifted Schur functions
via Okounkov's interpolation polynomials.)  The word ``balanced'' is used
for two reasons.  First, ordinary interpolation polynomials are defined via
vanishing conditions together with a triangularity condition; we have
replaced the latter by another, symmetrical, set of vanishing conditions
(thus obtaining a more ``balanced'' definition).  Second, the corresponding
hypergeometric sums (most notably Theorem \ref{thm:theta_bde} below) are
multivariate analogues of balanced, very-well-poised, hypergeometric
series; the corresponding identity for ordinary interpolation polynomials
is merely very-well-poised.  Moreover, although it is easy to degenerate
the sums to remove very-well-poisedness, there appears to be no
straightforward way to obtain a non-balanced sum without degenerating to
ordinary interpolation polynomials.

Note that it suffices to consider the generic extension of $\gamma$, in
which the extended points are all algebraically independent over $R$.  Thus
the existence of an interpolation polynomial of specified shape is
equivalent to the existence of a primitive solution to a certain finite
system of linear equations, or in other words to the simultaneous vanishing
of an appropriate collection of determinants.  Thus the space of perfect
bigrids is a closed, $R$-rational subscheme of the space of bigrids;
moreover, perfection is preserved under base change.

\begin{prop}
Let $\gamma$ be a perfect bigrid.  Then the {\em complementary} bigrid
$\tilde\gamma$ defined by $\tilde\gamma(\alpha,i,j)=\gamma(1-\alpha,n+1-i,m-j)$ is
perfect.
\end{prop}

\begin{proof}
Indeed, an interpolation polynomial of index $\lambda$ for
$\gamma$ is also an interpolation polynomial of index $m^n-\lambda$ for
$\tilde\gamma$, and vice versa.
\end{proof}

The simplest example of a perfect bigrid is the following.

\begin{prop}
Let $\gamma$ be an arbitrary bigrid of shape $m^1$.  Then $\gamma$ is
perfect, with associated family of (univariate) interpolation polynomials
given by
\[
P^*_j(x;\gamma) = \prod_{0\le i<j} (x\cdot \gamma(0,1,i))
                  \prod_{j<i\le m} (x\cdot \gamma(1,1,i)).
\]
\end{prop}

\begin{proof}
Indeed, the vanishing conditions say precisely that $\gamma(0,1,i)$ is a
zero for $0\le i<j$ and that $\gamma(1,1,i)$ is a zero for $j<i\le m$.
\end{proof}

This generalizes in several different ways to the multivariate case.  The
simplest of these is the case of ``monomial'' bigrids, in which
$\gamma(\alpha,i,j)$ is independent of $i$.  In other words, a monomial bigrid
$\gamma$ is constructed from a bigrid $\gamma_0$ of shape $m^1$ by
$\gamma(\alpha,i,j)=\gamma_0(\alpha,1,j)$.

\begin{prop}
All monomial bigrids are perfect, with associated interpolation
polynomials given by
\[
P^*_\lambda(x_1,\dots, x_n;\gamma)
\propto
\sum_{\pi\in S_n}
\prod_{1\le i\le n}
P^*_{\lambda_{\pi(i)}}(x_i;\gamma_0).
\]
\end{prop}

\begin{proof}
The given polynomial is clearly symmetric and homogeneous, so it remains to
verify the vanishing conditions.  For this, we observe that if $\mu\subset
m^n$ is such that $\mu\not\supset\lambda$, then for any permutation $\pi$,
$\mu_{i}<\lambda_{\pi(i)}$ for some $1\le i\le n$.  Then the $i$th factor
of the term for $\pi$ vanishes as required.  The vanishing conditions for
$\gamma(1,i,j)$ follow by the symmetrical argument.
\end{proof}

Similarly, the Schur case extends to balanced polynomials.  Fix $m$,
$n$, and let $\gamma_0$ be a bigrid of shape $(m+n-1)^1$.  Then we define
a Schur bigrid $\gamma$ by:
\[
\gamma(\alpha,i,j) = \gamma_0(\alpha,1,j+n-i).
\]

\begin{prop}
Schur bigrids are perfect, with associated interpolation polynomials given
by the formula
\[
P^*_\lambda(x_1,\dots, x_n;\gamma)
=
\frac{
\det(P^*_{\lambda_i+n-i}(x_j;\gamma_0))_{1\le i,j\le n}
}{
\det(P^*_{n-i}(x_j;\gamma_0))_{1\le i,j\le n}
}
\propto
\frac{
\det(P^*_{\lambda_i+n-i}(x_j;\gamma_0))_{1\le i,j\le n}
}{
\prod_{1\le i<j\le n} (x_i\cdot x_j)
}
\]
\end{prop}

\begin{proof}
We first observe that the numerator vanishes when $x_j=x_i$, so the result
is a polynomial; similarly, switching $x_i$ and $x_j$ negates both
numerator and denominator, so it is a symmetric polynomial.  Finally, the
numerator is homogeneous of degree $(m+n-1,m+n-1,\dots, m+n-1)$, while the
denominator is homogeneous of degree $(n-1,n-1,\dots, n-1)$.  Thus
$P^*_\lambda(;\gamma)$ as defined is a symmetric polynomial of degree $m$
on $\P^1(R)^n$.

For the vanishing conditions, we observe that, at least in the generic
case, the points at which we must vanish satisfy $x_i\ne x_j$ for all
$i\ne j$.  Thus it suffices to show that the numerator vanishes as
required; this follows by the same argument as in the monomial case.

Thus Schur bigrids are generically perfect with given interpolation
polynomials; the desired result is closed, so holds for all Schur bigrids.
\end{proof}

\begin{cor}
Any bigrid of shape $1^n$ is perfect.
\end{cor}

\begin{proof}
Indeed, any bigrid of shape $1^n$ is a Schur bigrid.
\end{proof}

There is a third important multivariate extension of the univariate case
which is not simply an obvious generalization of a special case for
ordinary interpolation polynomials.  In the monomial and Schur cases, the
polynomials have reasonable formulas, but unlike in the univariate case,
do not in general have nice factorizations.  It turns out that there is a
large class of bigrids for which the corresponding interpolation
polynomials do factor completely.

Let $\eta$ be a function $\eta:\{0,1,\dots, n\}\times\{0,1,\dots, m-1\}\to
\P^1(R)$, and construct a bigrid $\gamma$ by
\begin{align}
\gamma(0,i,j)&=\eta(i,j)\\
\gamma(1,i,j)&=\eta(i-1,j-1)
\end{align}
A bigrid constructed in this manner will be called a ``Cauchy'' bigrid.

\begin{prop}
Cauchy bigrids are perfect; indeed, we can construct corresponding
interpolation polynomials via the product expression:
\[
P^*_\lambda(x_1,\dots, x_n;\gamma)
=
\prod_{1\le i\le n}
\prod_{1\le j\le m}
x_i\cdot \eta(\lambda'_j,j-1).
\]
\end{prop}

\begin{proof}
Again, the only nonobvious fact is that $P^*_\lambda(;\gamma)$ as defined
satisfies the vanishing conditions.

Choose $\mu\subset m^n$ such that $\mu\not\supset\lambda$, and consider
\[
P^*_\lambda(\dots\gamma^+(0,i,\mu_i)\dots;\gamma)
=
\prod_{1\le i\le n}
\prod_{1\le j\le m}
\gamma^+(0,i,\mu_i)\cdot \eta(\lambda'_j,j-1)
\]
Since $\mu\not\supset\lambda$, there exists a position $1\le k\le n$ such
that $\mu_k<\lambda_k$.  If $k'$ is the largest such $k$, then we note that
\[
\lambda_{k'+1}=\mu_{k'+1}\le \mu_k<\lambda_k\le n,
\]
and thus in particular $\lambda'_l=k$ for $\mu_k<l\le \lambda_k$.  But then
\[
\gamma^+(0,k,\mu_k) = \gamma(0,k,\mu_k) = \eta(k,\mu_k)
=
\eta(\lambda'_{\mu_k+1},(\mu_k+1)-1),
\]
so the $i=k$, $j=\mu_k+1$ factor of the above product vanishes.

The other set of vanishing conditions follow by symmetry.
\end{proof}

\begin{rems}
Compare the proof of Lemma 6.3 of \cite{OkounkovA:1998a}.  Similarly,
interpolation polynomials for Cauchy bigrids will be related to the Cauchy
identity for interpolation theta functions, Theorem \ref{thm:interp_cauchy}
below.  In addition, when the Cauchy bigrid is also elliptic (see Section
\ref{sec:ell_bigrid}), the corresponding interpolation polynomials satisfy
a number of symmetries, which we will use to prove corresponding symmetries
for the general elliptic case.
\end{rems}

\begin{rems}
This is the only case (at least for a regular (Definition
\ref{def:regular}) bigrid) for which the associated interpolation
polynomials factor completely.  Indeed, this holds by direct computation
for shape $1^2$, and every equality required of a Cauchy bigrid is
supported on some truncation (Definition \ref{def:truncated}) of shape
$1^2$.
\end{rems}

The fourth ``elementary'' case is that of a ``delta'' bigrid, for which
$\gamma(0,i,j)=\gamma(1,i,j)$ whenever both sides are defined.  If such a
bigrid is perfect, then the corresponding interpolation polynomials must
satisfy
\[
P^*_\lambda(\dots \gamma_0(i,\mu_i)\dots;\gamma)
\propto \delta_{\lambda\mu}
\]
(here $\gamma_0(i,j)=\gamma(0,i,j)$ or $\gamma(1,i,j)$, whichever is
well-defined).  Indeed, aside from the fact that the extended bigrid
$\gamma^+$ need not satisfy $\gamma^+(0,i,j)=\gamma^+(1,i,j)$, the first
set of vanishing conditions says that this vanishes unless $\lambda\subset
\mu$, while the second requires that it vanish unless $\mu\subset\lambda$.
That all delta bigrids are perfect requires some additional machinery, and
is shown in Proposition \ref{prop:delta_is_perfect} below.

\begin{defn}\label{def:truncated}
Let $\gamma$ be a perfect bigrid.  We define four {\em truncated} bigrids
as follows.  The bigrids $\gamma^-$, of shape $m^{n-1}$, and $\gamma_-$, of
shape $(m-1)^n$, are simply the restrictions of $\gamma$ to the appropriate
domains; the bigrids ${}^-\gamma$ and ${}_-\gamma$ are then defined by
complementing the restrictions of $\tilde\gamma$; thus
\begin{align}
{}^-\gamma(\alpha,i,j) &= \gamma(\alpha,i+1,j)\\
{}_-\gamma(\alpha,i,j) &= \gamma(\alpha,i,j+1)
\end{align}
\end{defn}

The relation of truncation to interpolation polynomials is described by
the following propositions.

\begin{prop}
Let $\gamma$ be a bigrid of shape $m^n$, let $\lambda\subset m^{n-1}$ be a
partition, and let $P^*_\lambda(;\gamma)$ be a corresponding interpolation
polynomial.  Then the specialization
\[
P^*_\lambda(x_1,x_2,\dots, x_{n-1},\gamma(0,n,0);\gamma)
\]
is an interpolation polynomial for $\gamma^-$ of index $\lambda$.
Similarly, if $P^*_{m\cdot\lambda}(;\gamma)$ is an interpolation polynomial
for $\gamma$ of index $m\cdot\lambda$, then
\[
P^*_{m\cdot\lambda}(\gamma(1,1,m),x_1,x_2,\dots, x_{n-1};\gamma)
\]
is an interpolation polynomial for ${}^-\gamma$ of index $\lambda$.
\end{prop}

We cannot quite conclude, however, that if $\gamma$ is perfect, then so
are $\gamma^-$ and ${}^-\gamma$; the constructed interpolation polynomials
could be trivial even if the original polynomials were primitive.

\begin{prop}
Let $\gamma$ be a bigrid of shape $m^n$, let $\lambda\subset (m-1)^n$ be a
partition, and let $P^*_\lambda(;\gamma_-)$ be a corresponding interpolation
polynomial for $\gamma_-$.  Then
\[
\prod_{1\le i\le n} x_i\cdot \gamma(1,1,m)
P^*_\lambda(x_1,\dots, x_n;\gamma_-)
\]
is an interpolation polynomial of index $\lambda$ for $\gamma$.
Similarly, if $P^*_\lambda(;{}_-\gamma)$ is an interpolation polynomial,
then 
\[
\prod_{1\le i\le n} x_i\cdot \gamma(0,n,0)
P^*_\lambda(x_1,\dots, x_n;{}_-\gamma)
\]
is an interpolation polynomial of index $1^n+\lambda$ for $\gamma$.
\end{prop}

With a slight additional hypothesis, we obtain a converse:

\begin{prop}\label{prop:trunc2}
Let $\gamma$ be a bigrid of shape $m^n$, and let $\lambda\subset
(m-1)^n$. If the point $\gamma(1,1,m)$ has unit inner product with all
points in the image of $\gamma_-$, then all interpolation polynomials of
index $\lambda$ for $\gamma$ are obtained from interpolation polynomials
for $\gamma_-$ via the above proposition.  Similarly, if the point
$\gamma(0,n,0)$ has unit inner product with all points in the image of
${}_-\gamma$, then all interpolation polynomials of index $1^n+\lambda$ for
$\gamma$ are obtained via the above proposition.
\end{prop}

\begin{proof}
Passing to the complementary bigrid reduces the second claim to the first.
Thus consider an interpolation polynomial $P^*_\lambda(;\gamma)$.  The
partition $m$ is not contained in $\lambda$, and thus the polynomial
\[
P^*_\lambda(\gamma(1,1,m),x_2,\dots, x_n;\gamma)=0.
\]
It follows that $P^*_\lambda(x_1,\dots, x_n;\gamma)$ is a multiple of
$x_1\cdot \gamma(1,1,m)$, and thus, by symmetry, that it is a multiple of
\[
\prod_{1\le i\le n} x_i\cdot\gamma(1,1,m).
\]
We need to show that
\[
\frac{P^*_\lambda(\gamma(1,1,m),x_2,\dots, x_n;\gamma)}
{\prod_{1\le i\le n} x_i\cdot\gamma(1,1,m)}
\]
satisfies the vanishing identities required of an interpolation polynomial
for $\gamma_-$; by the hypothesis, it suffices to show that
$P^*_\lambda(;\gamma)$ itself satisfies the identities.  But the required
identities are a subset of the vanishing identities that $P^*_\lambda$
already satisfies.
\end{proof}

\begin{cor}
Let $\gamma$ be a perfect grid of shape $m^n$.  If $\gamma(1,1,m)$ has unit
inner product with the image of $\gamma_-$, then $\gamma_-$ is perfect; if
$\gamma(0,n,0)$ has unit inner product with the image of ${}_-\gamma$, then
${}_-\gamma$ is perfect.
\end{cor}

Naturally, the interesting cases of perfect bigrids are those for which the
interpolation polynomials are uniquely determined.  This is an open
condition (the complement of a closed subscheme again cut out by
determinants); unfortunately, the ideal for the complement appears to be
prohibitively complicated.  We thus need to find a simpler set of
sufficient conditions for uniqueness.  As in \cite[Definition
  2.1]{OkounkovA:1998c}, the idea is to choose a subset of the equations
and add a normalization condition in order to obtain a square system of
linear equations with a simple determinant.

\begin{defn}\label{def:regular}
A bigrid is {\em regular} if for all $i\ge i'$, $j<j'$,
$\alpha,\beta\in \{0,1\}$,
\[
\gamma(\alpha,i,j)\cdot\gamma(\beta,i',j') \in R^*.
\]
\end{defn}

Note that regularity of a bigrid $\gamma$ implies that the extra hypotheses of
Proposition \ref{prop:trunc2} hold for all truncations of $\gamma$;
regularity is slightly stronger in that the additional conditions
\[
\gamma(1,i,j)\cdot\gamma(0,i',j-1)\in R^*
\]
for $i\ge i'$ must also hold.  (It is unclear whether these extra
conditions are necessary for primitive interpolation polynomials to be unique
and form a basis; they are needed for our proof below, however.)  It
follows therefore that if $\gamma$ is a regular perfect bigrid, then so are
$\gamma_-$ and ${}_-\gamma$.

Given a bigrid $\gamma$ of shape $m^n$ ($m>0$), a partition
$\lambda\subset m^n$, and an integer $0\le l\le n$, we define a point
$\gamma_l(\lambda)\in
\P^1(R)^n$ by
\[
\gamma_l(\lambda) =
(\gamma([1\le l],1,\lambda_1),\gamma([2\le l],2,\lambda_2),\dots,
\gamma([n\le l],n,\lambda_n)),
\]
where $[i\le l]$ is 1 if $i\le l$ and $0$ otherwise.  For a pair of
partitions $\lambda,\mu\subset m^n$, we define
$\gamma_\lambda(\mu)=\gamma_l(\mu)$ with $l$ determined as follows.  Let
$l_0$ be the largest index such that $\mu_l\ne \lambda_l$ ($l_0=0$ if no
such index exists), and let $l_1$ similarly be the largest index such that
$\mu_l=m$.  If $l_0=0$ or $\mu_{l_0}<\lambda_{l_0}$ then take $l=l_1$;
otherwise, take $l=l_0$.  Thus, for instance, for $m=n=3$, $\lambda=21$, we
have:
\begin{align}
\gamma_{21}(0) &= (\gamma(0,1,0),\gamma(0,2,0),\gamma(0,3,0))\\
\gamma_{21}(1) &= (\gamma(0,1,1),\gamma(0,2,0),\gamma(0,3,0))\\
\gamma_{21}(2) &= (\gamma(1,1,2),\gamma(0,2,0),\gamma(0,3,0))\\
\gamma_{21}(11) &= (\gamma(0,1,1),\gamma(0,2,1),\gamma(0,3,0))\\
\gamma_{21}(21) &= (\gamma(1,1,2),\gamma(0,2,1),\gamma(0,3,0))\\
\gamma_{21}(22) &= (\gamma(1,1,2),\gamma(1,2,2),\gamma(0,3,0))\\
\gamma_{21}(111) &= (\gamma(1,1,1),\gamma(1,2,1),\gamma(1,3,1))\\
\gamma_{21}(211) &= (\gamma(1,1,2),\gamma(1,2,1),\gamma(1,3,1))\\
\gamma_{21}(221) &= (\gamma(1,1,2),\gamma(1,2,2),\gamma(1,3,1))\\
\gamma_{21}(222) &= (\gamma(1,1,2),\gamma(1,2,2),\gamma(1,3,2))
\end{align}

\begin{defn}
Fix a ring $R$, positive integers $m$, $n$, a $R$-valued bigrid $\gamma$
of shape $m^n$, and a partition $\lambda\subset m^n$.  A primitive symmetric
polynomial $p$ of degree $m$ and $n$ variables is a {\em
  quasi-interpolation polynomial} of index $\lambda$ and bigrid $\gamma$ if
for all partitions $\mu\subset m^n$, $\mu\ne\lambda$,
\[
p(\gamma_\lambda(\mu))=0.
\]
\end{defn}

\begin{lem}\label{lem:quasi1}
If the bigrid $\gamma$ of shape $m^n$ is regular, then for all
$\lambda\subset m^n$, there exists a quasi-interpolation polynomial of
index $\lambda$, unique up to scale.  The quasi-interpolation polynomials
for $\gamma$ span the space of symmetric polynomials of degree $m$ on
$\P^1(R)^n$.
\end{lem}

\begin{proof}
Introduce a total ordering $<_l$ on partitions by first sorting the parts
into increasing order, then lexicographically ordering the sorted
partitions.  We construct a corresponding filtration of $\Lambda^m_n$ as
follows.  (That is, a sequence of submodules $F_\mu(\lambda;\gamma)$ such
that $F_\mu(\lambda;\gamma)\supset F_\nu(\lambda;\gamma)$ whenever
$\mu<_l\nu$.)  Recall from Lemma \ref{lem:sf_decomp} the short exact
sequence:
\[
0\to \Lambda^{m-1}_n\mathop{\to}\limits^f \Lambda^m_n\mathop{\to}\limits^g \Lambda^m_{n-1}\to 0,
\]
in which the map $f$ is defined by
\[
f(p)(x_1,\dots, x_n) = p(x_1,\dots, x_n)
\prod_{1\le i\le n} (\gamma(0,n,0)\cdot x_i)
\]
and the map $g$ is defined by
\[
g(p)(x_1,\dots, x_{n-1})=p(x_1,\dots, x_{n-1},\gamma(0,n,0)).
\]
The submodules $F_\mu(\lambda;\gamma)\subset \Lambda^m_n$ are then defined
as follows.  First, as a base case, $F_0(\lambda;\gamma)$ is always taken
to be $\Lambda^m_n$; this in particular covers the case $m=0$ or $n=0$.
If $\mu_n=\lambda_n=0$, $m>0$, then
\[
F_\mu(\lambda;\gamma)=g^{-1}(F_\mu(\lambda;\gamma_-));
\]
if $\mu_n=0$ but $\lambda_n>0$, then
\[
F_\mu(\lambda;\gamma)=g^{-1}(F_\mu(\lambda;\gamma')),
\]
where $\gamma'(\alpha,i,j):=\gamma(0,i,j)$ is a delta bigrid of shape $m^{n-1}$.
If $\mu_n>0$, $\lambda_n>0$, then
\[
F_\mu(\lambda;\gamma)=f(F_{\mu-1^n}(\lambda-1^n;{}_-\gamma)).
\]
Finally, if $\mu_n>0$ but $\lambda_n=0$, then
\[
F_\mu(\lambda;\gamma)=f(F_{\mu-1^n}(0;\gamma'')),
\]
where $\gamma''(\alpha,i,j):=\gamma(1,i,j-1)$ is a delta bigrid of shape
$(m-1)^n$.

Moreover, since the short exact sequences used to construct the filtration
are all split, we can construct a corresponding sequence of polynomials
$p_\mu(;\lambda;\gamma)$ such that
\[
F_\mu = \la p_\nu:\mu\le_l \nu\ra.
\]
(When $m=0$ or $n=0$, we take $p_0=1$)
From the construction of this filtration, it follows easily that
\[
p_\mu(\gamma_\lambda(\nu);\lambda;\gamma)=0
\]
if $\nu<_l\mu$, while
\[
p_\mu(\gamma_\lambda(\mu);\lambda;\gamma)\in R^*,
\]
using the fact that $\gamma$ is regular.  (Note that if $\gamma$ is
regular, then so are the delta bigrids $\gamma'$, $\gamma''$ used above)
In particular, the determinant of the matrix
\[
(p_\mu(\gamma_\lambda(\nu)))_{\mu,\nu\subset m^n}
\]
is a unit, so for any collection of values $f_\mu\in R$, there is
a unique polynomial $p$ such that
\[
p(\gamma_\lambda(\mu))=f_\mu.
\]
Taking $f_\mu=\delta_{\lambda\mu}$ gives the desired existence and
uniqueness claims for the quasi-interpolation polynomials.

For $\nu\le_l\lambda$, $\gamma_\lambda(\nu)$ is independent of $\lambda$;
it follows that
\[
P^*_\lambda(\gamma_{m^n}(\nu))=0
\]
when $\nu<_l\lambda$, while
\[
P^*_\lambda(\gamma_{m^n}(\lambda))\in R^*.
\]
Thus the corresponding matrix has unit determinant, which immediately
implies that the $P^*_\lambda$ form a basis of $\Lambda^m_n$.
\end{proof}

\begin{rem}
Compare the proof of Proposition 2.6 of \cite{OkounkovA:1998c}.
\end{rem}

\begin{thm}
Let $\gamma$ be a regular, perfect bigrid of shape $m^n$.  Then the
corresponding primitive interpolation polynomials are uniquely determined up
to scale, and form a basis of the space of symmetric polynomials of degree
$m$ on $\P^1(R)^n$.  Moreover, for all $\lambda$,
\[
P^*_\lambda(\dots \gamma^+(0,i,\lambda_i)\dots;\gamma)
\]
is a primitive polynomial in the indeterminates $\gamma^+(0,i,m)$, $1\le
i\le n$.
\end{thm}

\begin{proof}
We find that the vanishing conditions defining the quasi-interpolation
polynomials form a subset of the conditions defining the interpolation
polynomials, so each interpolation polynomial must be a unit multiple of
the corresponding quasi-interpolation polynomial.  We then find that
\[
P^*_\lambda(\dots \gamma^+(0,i,\lambda_i)\dots;\gamma)\in R^*
\]
in the special case
\[
\gamma^+(0,i,m)=\gamma(1,i,m),
\]
which implies primitivity in general.
\end{proof}

\begin{cor}
If $\gamma$ is a regular, perfect bigrid, then the same is true for any
truncation of $\gamma$.
\end{cor}

\begin{proof}
The only thing remaining to prove is that if $\lambda_n=0$ and
$P^*_\lambda(;\gamma)$ is a primitive interpolation polynomial for $\gamma$,
then substituting $\gamma(0,n,0)$ for one argument gives another primitive
polynomial.  But the corresponding statement is true for
quasi-interpolation polynomials.
\end{proof}

In the sequel, we will restrict our attention to regular, perfect bigrids.

\begin{prop}\label{prop:delta_is_perfect}
Delta bigrids are perfect.
\end{prop}

\begin{proof}
Since regular delta bigrids are Zariski-dense in the space of all delta
bigrids, it suffices to consider the case of a regular delta bigrid
$\gamma$.  Let $P^*_\lambda(;\gamma)$ denote the corresponding
quasi-interpolation polynomials; thus
\[
P^*_\lambda(\dots \gamma(0,i,\mu_i)\dots;\gamma)=\delta_{\lambda\mu}.
\]
Taking a ring extension as necessary, we
can also choose a regular Cauchy bigrid $\gamma'$ such that
$\gamma'(0,i,j)=\gamma(0,i,j)$, with corresponding interpolation
polynomials.  But then
\[
P^*_\lambda(\dots \gamma^+(0,i,\mu_i)\dots;\gamma')=0
\]
unless $\lambda\subseteq\mu$; thus the polynomials corresponding to
$\gamma$ are triangular with respect to the basis $P^*_\lambda(;\gamma')$
under the inclusion ordering.  In particular,
\[
P^*_\lambda(;\gamma') = \sum_{\mu\supseteq\lambda} c_{\lambda\mu}
P^*_\mu(;\gamma)
\]
for appropriate coefficients $c_{\lambda\mu}$, and therefore
\[
P^*_\lambda(;\gamma) = \sum_{\mu\supseteq\lambda} d_{\lambda\mu}
P^*_\mu(;\gamma')
\]
for another set of coefficients $d_{\lambda\mu}$.  But then if
$\mu\not\supseteq\lambda$,
\[
P^*_\lambda(\dots \gamma^+(0,i,\mu_i)\dots;\gamma)
=
\sum_{\lambda\supseteq\nu\supseteq\mu} d_{\lambda\nu} P^*_\nu(\dots
\gamma^+(0,i,\mu_i)\dots;\gamma)
= 0.
\]
The other set of vanishing conditions follow symmetrically; the result
follows.
\end{proof}

Given a pair $(C,\phi)$, where $C$ is a curve of genus 1 and $\phi$ is a
quadratic map $C\to\P^1$, the line bundle ${\cal O}(1)$ on $\P^1$ pulls
back to a line bundle on $C$, on which the Galois group $\Z/2\Z$ acts.
Thus symmetric polynomials on $P^1$ pull back to sections of the
corresponding line bundle on $C^n$, invariant under the action of $BC_n$.
When $C$ is a complex elliptic curve and $\phi$ preserves the identity,
symmetric polynomials of degree $m$ in $n$ variables thus pull back to
$BC_n$-symmetric theta functions of degree $m$.  In the next section, we
will construct certain ``interpolation theta functions'', pull-backs of
interpolation polynomials associated to ``elliptic'' bigrids.  We will
consider the case of general curves $C$ in Section \ref{sec:modular} below.

\section{Interpolation theta functions}\label{sec:int_theta}

\begin{defn}
Fix integers $m,n\ge 0$, a partition $\lambda\subset m^n$, and generic
points $a,b,q,t\in \C^*$; also choose another generic point $v\in \C^*$.
The {\em interpolation theta function}
\[
P^{*(m,n)}_\lambda(x_1,\dots, x_n;a,b;q,t;p)
\]
is the unique function on $(\C^*)^{n}\times (\C^*)^4$ with the following
properties:
\begin{itemize}
\item[(1)] $P^{*(m,n)}_\lambda$ is a $BC_n$-symmetric theta function of
  degree $m$ in the variables $x_1,\dots, x_n$.
\item[(2)] For any partition $\mu\subset m^n$, $\mu\ne \lambda$, let
$l_0$ be the largest index such that $\mu_l\ne \lambda_l$, let $l_1$
  be the largest index such that $\mu_l=m$ ($0$ if none exists).  If
  $\mu_{l_0}<\lambda_{l_0}$, then set $l=l_1$, otherwise set $l=l_0$.
Then
\[
P^{*(m,n)}_\lambda(b q^{m-\mu_1},\dots, b q^{m-\mu_l} t^{l-1},
                   a q^{\mu_{l+1}} t^{n-l-1}, \dots, a q^{\mu_n};a,b;q,t;p)
=
0.
\]
\item[(3)] $P^{*(m,n)}_\lambda$ is normalized by the specialization
\[
P^{*(m,n)}_\lambda(\dots v t^{n-i}\dots;a,b;q,t;p)
=
C^0_\lambda(t^{n-1} av,a/v;q,t;p)
C^0_{m^n-\lambda}(t^{n-1} bv,b/v;q,t;p)
\]
\end{itemize}
\end{defn}

Indeed, aside from the normalization, this is simply the pull-back of the
definition of quasi-interpolation polynomials relative to the bigrid
\[
\gamma(0,i,j):=\phi(a q^j t^{n-i})\quad
\gamma(1,i,j):=\phi(b q^{m-j} t^{i-1}),
\]
where $\phi$ maps the elliptic curve $\C^*/\la p\ra$ to $\P^1(\C)$
via a pair of $BC_1$-symmetric theta functions of degree 1; i.e., for some
$c,d\in \C^*$,
\[
\phi(z) = (\theta(c z^{\pm 1};p),\theta(d z^{\pm 1};p)).
\]
Thus by Lemma \ref{lem:quasi1}, interpolation theta functions exist and are
unique, at least for generic values of $a,b,q,t,v$.  We note in particular
that the associated bigrid is regular precisely when
\begin{align}
q^j t^i&\notin ab q^{m-1} t^{n-1} \la p\ra & 0\le i<n,0\le j<m\\
(b/a) q^j t^i&\notin ab q^{m-1} t^{n-1}\la p\ra & |i|\le n-1,|j|\le m-2\\
q^j t^i&\notin \la p\ra & 0\le i<n,1\le j<m\\
(b/a) q^j t^i&\notin \la p\ra & |i|\le n-1,|j|\le m-1\\
q^j t^i&\notin (ab)^{-1} q^{1-m} t^{1-n} \la p\ra & 0\le i<n,2\le j<m\\
(b/a) q^j t^i&\notin (ab)^{-1} q^{1-m} t^{1-n}\la p\ra & |i|\le n-1,|j|\le m-2
\end{align}
Moreover, we can also conclude that $P^{*(m,n)}_\lambda(;a,b;q,t;p)$ is a
meromorphic theta function in each of $a,b,q,t,v$.  (We will show below
that it is independent of $v$, explaining its omission from the notation)

As one might expect from the fact that we are considering these functions
at all, interpolation theta functions indeed satisfy extra vanishing
conditions corresponding to the fact that the associated bigrid is perfect.
The normalization can thus be justified as follows: viewed as a function of
$v$, the left-hand side must (by extra vanishing) vanish wherever the
right-hand side vanishes.  By degree considerations and the theta function
condition, it follows that the ratio must be independent of $v$.

As in \cite{bcpoly}, we first show that the interpolation theta functions
satisfy a difference equation, and in the process prove independence of $v$.

\begin{thm}\label{thm:theta_diff}
Let $c,d\in \C^*$ be chosen so that $q^m t^{n-1} abcd=p$.  Then
\begin{align}
D^{(n)}(a,b,c,d&;q,t;p)
P^{*(m,n)}_{\lambda}(;\sqrt{q}a,\sqrt{q}b;q,t;p)
\notag\\
&=
\prod_{1\le i\le n} 
\theta(ab q^m t^{n-i},ac q^{\lambda_i} t^{n-i},bc q^{m-\lambda_i} t^{i-1};p)
P^{*(m,n)}_{\lambda}(;a,b;q,t;p)
\label{eq:diff_eq}
\end{align}
\end{thm}

\begin{proof}
To begin with, we take the $v$ normalizing the interpolation theta function
in the left-hand side of the equation to be $\sqrt{q}$ times the $v$
normalizing the interpolation theta function on the right.

We first need to show that the left-hand side satisfies the relevant
vanishing conditions.  Suppose more generally that we evaluate the
left-hand side at a point of the form
\[
b q^{m-\mu_1},\dots, b q^{m-\mu_l} t^{l-1},
                   a q^{\mu_{l+1}} t^{n-l-1}, \dots, a q^{\mu_n}
\]
for an arbitrary choice of $l$, where $\mu$ simply satisfies the conditions
\begin{align}
m\ge \mu_1\ge \mu_2\ge\cdots\ge\mu_l,\\
\mu_{l+1}\ge \mu_{l+2}\ge\cdots\ge \mu_n\ge 0,
\end{align}
but need not satisfy $\mu_l\ge \mu_{l+1}$.  Carrying this point through the
difference operator and using the $BC_n$ symmetry of $P^{*(m,n)}$, we find
that in each term, $P^{*(m,n)}(;\sqrt{q}a,\sqrt{q}b;q,t;p)$ is evaluated at
a point of the same form, with $a,b$ replaced by $q^{1/2}a$, $q^{1/2} b$,
and $\mu$ replaced by $\nu$ satisfying
\begin{align}
\mu_i\le \nu_i\le \mu_i+1,&\quad 1\le i\le l\\
\mu_i-1\le \nu_i\le \mu_i,&\quad l+1\le i\le n
\end{align}
(Indeed, $\nu_i=\mu_i+1/2-\sigma_i$ if $1\le i\le l$, and
$\nu_i=\mu_i-1/2-\sigma_i$ if $l+1\le i\le n$.)  We furthermore observe
that if $m\ge \nu_1$, $\nu_i<\nu_{i+1}$ for $i\ne l$, or $\nu_n<0$, then
that term of the expansion necessarily vanishes.  Thus the only surviving
$\nu$ satisfy
\begin{align}
m\ge \nu_1\ge\nu_2\ge\cdots\ge \nu_l,\\
\nu_{l+1}\ge \nu_{l+2}\ge\cdots\ge \nu_n\ge 0.
\end{align}
In particular, if $\mu$ and $l$ are chosen to make the right-hand side
vanish, then all terms on the left-hand side will also vanish, as required.

It follows that
\[
D^{(n)}(a,b,c,d;q,t;p)
P^{*(m,n)}_{\lambda}(;\sqrt{q}a,\sqrt{q}b;q,t;p)
\propto
P^{*(m,n)}_{\lambda}(;a,b;q,t;p),
\]
for some scale factor independent of $x_1,\dots, x_n$.  Taking
$\mu=\lambda$, and $l$ to be the largest index such that $\lambda_l=m$, we
find that the interpolation theta function on the right has a nonzero value
(by Lemma \ref{lem:quasi1}), and that only one term on the left survives.
The desired dependence of the scale factor on $c$ follows immediately.

If we now take $c=v$ and set $x_i=v t^{n-i}$, again only one term survives
on the left, and the difference equation follows.  By symmetry, the same
difference equation would hold if we had normalized the interpolation
theta function on the left using $q^{-1/2} v$ instead, and thus these
interpolation theta functions agree.  In other words, interpolation theta
functions are invariant under $v\mapsto qv$; since $q$ is generic, they
must indeed be independent of $v$.
\end{proof}

\begin{rem}
In particular, we find that the interpolation theta functions are solutions
to a ``generalized eigenvalue problem'' (ala
\cite{SpiridonovVP/ZhedanovAS:2001}), to wit:
\begin{align}
D^{(n)}&(a,b,c,\frac{p}{q^m t^{n-1}abc};q,t;p)
P^{*(m,n)}_{\lambda}(;\sqrt{q}a,\sqrt{q}b;q,t;p)\\
&=
\left(\prod_{1\le i\le n} 
\frac{\theta(ac q^{\lambda_i} t^{n-i},bc q^{m-\lambda_i} t^{i-1};p)}
     {\theta(ac' q^{\lambda_i} t^{n-i},bc' q^{m-\lambda_i} t^{i-1};p)}\right)
D^{(n)}(a,b,c',\frac{p}{q^m t^{n-1}abc'};q,t;p)
P^{*(m,n)}_{\lambda}(;\sqrt{q}a,\sqrt{q}b;q,t;p)
\notag
\end{align}
\end{rem}

\begin{cor}\label{cor:theta_vanishing}
Fix a partition $\lambda\subset m^n$, and let $1\le l\le n$ and a sequence
$0\le \mu_n\le\mu_{n-1}\le \cdots\le \mu_{l+1}$ be chosen such that
$\mu_{l+1}<\lambda_{l+1}$.  Then
\[
P^{*(m,n)}_\lambda(x_1,\dots, x_l,a q^{\mu_{l+1}} t^{n-l-1}, \dots, a q^{\mu_n} ;a,b;q,t;p)
=0.
\]
Similarly, if $m\ge \mu_1\ge \mu_2\ge\cdots\ge \mu_l$ is chosen such that
$\mu_l>\lambda_l$, then
\[
P^{*(m,n)}_\lambda(
b q^{m-\mu_1},\dots, b q^{m-\mu_l} t^{l-1},
x_{l+1},\dots, x_n
;a,b;q,t;p)
=0.
\]
\end{cor}

\begin{proof}
By analytic continuation, it suffices to prove that if $\mu_1,\dots, \mu_n$
are as above, such that either $\mu_l>\lambda_l$ or
$\mu_{l+1}<\lambda_{l+1}$, then
\[
P^{*(m,n)}_\lambda(
b q^{m-\mu_1},\dots, b q^{m-\mu_l} t^{l-1},
a q^{\mu_{l+1}} t^{n-l-1}, \dots, a q^{\mu_n}
;a,b;q,t;p)=0.
\]
We proceed by induction on
\[
\sum_{1\le i\le l} (m-\mu_i)+
\sum_{l+1\le i\le n} \mu_i.
\]
Indeed, if we evaluate both sides of the difference equation
\eqref{eq:diff_eq} at this point, we find by the inductive assumption that
at most one term on the left survives (that for which $\sigma_i=1$, $1\le
i\le n$).  If $\mu_l>\lambda_l$, then setting $c=b q^{m-\mu_l} t^{l-1}$
makes this term vanish without annihilating the factor on the right;
similarly, if $\mu_{l+1}<\lambda_{l+1}$, we set $c=a q^{\mu_{l+1}}
t^{n-l-1}$.  The result follows.
\end{proof}

\begin{cor}\label{cor:comp_symm_P}
We have the symmetry property
\[
P^{*(m,n)}_{m^n-\lambda}(;a,b;q,t;p)
=
P^{*(m,n)}_{\lambda}(;b,a;q,t;p)
\]
\end{cor}

\begin{proof}
This follows immediately from Corollary \ref{cor:theta_vanishing} and the
fact that the normalization has this symmetry.
\end{proof}

This also gives the following commutation relation for the difference
operators:

\begin{cor}\label{cor:quasicommutation}
If $cd = c'd'$, then
\begin{align}
D^{(n)}(a,b,c',d';q,t;p)
D^{(n)}(q^{1/2}a,q^{1/2}b,&q^{-1/2}c,q^{-1/2}d;q,t;p)\\
&=
D^{(n)}(a,b,c,d;q,t;p)
D^{(n)}(q^{1/2}a,q^{1/2}b,q^{-1/2}c',q^{-1/2}d';q,t;p)\notag
\end{align}
\end{cor}

\begin{proof}
Suppose $q^m t^{n-1}abcd = p$.  Applying both sides to interpolation theta
functions, we find that the relation holds on the full space of
$BC_n$-symmetric theta functions of degree $m$.  Now, applied to a
$BC_n$-symmetric meromorphic function $f$, both sides can be expressed in the
form
\[
\sum_{\epsilon\in \{-1,0,1\}^n}
C_\epsilon(\dots x_i\dots;a,b,c,d,c',d')
f(\dots q^{\epsilon_i} x_i\dots)
\]
for suitable functions $C_\epsilon$; the claim is that the two resulting
families of functions agree.  For $m$ sufficiently large, this follows
from the action on theta functions; since each $C_\epsilon$ is a meromorphic
theta function, the result follows by analytic continuation.
\end{proof}

When $q^m t^n ab=pq$, the associated bigrid is a Cauchy bigrid; we thus
obtain the following formula for interpolation theta functions in that
case:

\begin{prop}
If $q^m t^n ab=pq$, then
\[
P^{*(m,n)}_\lambda(\dots x_i\dots;a,b;q,t;p)
=
\prod_{1\le i\le n,1\le j\le m}
\theta(a t^{n-\lambda'_j} q^{j-1} x_i,a t^{n-\lambda'_j} q^{j-1}/x_i;p),
\]
\end{prop}

We will be using this Cauchy specialization quite frequently; the point is
that by using the difference equation, one can shift $ab$ arbitrarily, but
cannot change $a/b$.  Thus the Cauchy specialization, which includes cases
having all values of $a/b$, gives a suitable base case for inductive
arguments.  (It thus plays much the same role that Macdonald polynomials
did in \cite{bcpoly}.)  For instance:

\begin{prop}
Interpolation theta functions satisfy the identity
\begin{align}
P^{*(m,n)}_\lambda(\dots a q^{\lambda_i} t^{n-i}\dots;&a,b;q,t;p)\\
&=
C^0_{m^n-\lambda}(\frac{pq}{q^m ab};q,t;p)
C^+_{m^n-\lambda}(\frac{t^{n-1}b}{q^m a};q,t;p)
\frac{C^-_\lambda(pq,t;q,t;p) C^+_\lambda(t^{2n-2}a^2;q,t;p)}
     {C^0_\lambda(t^n;q,t;p)}\notag
\end{align}
\end{prop}

\begin{proof}
Using the difference equation, it follows immediately that this formula
holds, up to a scale factor depending only on $a/b$ (more precisely,
independent under $(a,b)\mapsto (\sqrt{q}a,\sqrt{q}b)$, but by analytic
continuation, this amounts to the same thing); taking the Cauchy
specialization shows that this scale factor is 1.
\end{proof}

The other important base case is the delta case:

\begin{prop}
If $q^m t^{n-1} ab=1$, then
\begin{align}
P^{*(m,n)}_\lambda(&
\dots a q^{\mu_i} t^{n-i}\dots
;a,b;q,t;p)\notag\\
&=
\delta_{\lambda\mu}
C^0_{m^n-\lambda}(pqt^{n-1};q,t;p)
C^+_{m^n-\lambda}(\frac{1}{q^{2m} a^2};q,t;p)
\frac{C^-_\lambda(pq,t;q,t;p) C^+_\lambda(t^{2n-2}a^2;q,t;p)}
     {C^0_\lambda(t^n;q,t;p)}\notag
\end{align}
\end{prop}

The dependence of interpolation theta functions on $m$ is fairly simple;
we have the following identity.

\begin{prop}
If $\lambda\subset m^n$, then
\[
P^{*(m+k,n)}_\lambda(\dots x_i\dots;a,b;q,t;p)
=
\prod_{1\le i\le n} \theta(b x_i,b/x_i;q;p)_k
P^{*(m,n)}_\lambda(\dots x_i\dots;a,bq^k;q,t;p).
\]
In particular,
\[
P^{*(m,n)}_0(\dots x_i\dots;a,b;q,t;p)
=
\prod_{1\le i\le n} \theta(b x_i,b/x_i;q;p)_m
\]
\end{prop}

In particular, the dependence on $m$ is independent of $\lambda$,
and thus the following abelian functions are well-defined:

\begin{defn}
The interpolation abelian function $R^{*(n)}_\lambda(;a,b;q,t;p)$
is defined by
\[
R^{*(n)}_\lambda(;a,b;q,t;p)
=
\frac{P^{*(m,n)}_\lambda(;a,q^{-m}b;q,t;p)}
     {P^{*(m,n)}_0(;a,q^{-m}b;q,t;p)},
\]
for any $m\ge \lambda_1$.
\end{defn}

Of course, the price we pay for dealing with these functions is that we
cannot evaluate them at points $b q^{\mu_i} t^{n-i}$, and as a result the
complementation symmetry becomes more complicated.  However, the
abelianness, the lack of dependence on $m$, and the resulting ability to
have an infinite family of functions, is often an overall win, and as a
result we will tend to use these functions in preference to the
interpolation theta functions.

Translating the above identities to the interpolation abelian functions,
we obtain:

\begin{prop}
If $t^{n-1}abcd=p$, then
\begin{align}
D^{(n)}(a,b,c,d&;q,t;p)
R^{*(n)}_{\lambda}(;\sqrt{q}a,\sqrt{q}b;q,t;p)\notag\\
&=
\prod_{1\le i\le n} 
\theta(ab t^{n-i},ac q^{\lambda_i} t^{n-i},bc q^{-\lambda_i} t^{i-1};p)
R^{*(n)}_{\lambda}(;a,b;q,t;p)
\label{eq:diff_eq_Rs}
\end{align}
As an abelian function in $v$,
\[
R^{*(n)}_\lambda(v t^{n-i};a,b;q,t;p)
=
\Delta^0_\lambda(t^{n-1}a/b|t^{n-1}av,a/v;q,t;p)
\]
If $t^n ab=pq$, then
\[
R^{*(n)}_\lambda(x;a,b;q,t;p)
=
\prod_{1\le i\le n,1\le j\le \lambda_1}
\frac{\theta(a t^{n-\lambda'_j} q^{j-1} x_i,a t^{n-\lambda'_j} q^{j-1}/x_i;p)}
     {\theta(a t^{n} q^{j-1} x_i,a t^{n} q^{j-1}/x_i;p)}
\]
Finally,
\[
R^{*(n)}_\lambda(a q^{\lambda_i} t^{n-i};a,b;q,t;p)
=
\frac{C^+_\lambda(t^{2n-2}a^2;q,t;p)C^0_\lambda(pqt^{n-1}a/b;q,t;p)}
     {C^+_\lambda(t^{n-1}a/b;q,t;p)C^0_\lambda(t^{n-1}ab;q,t;p)}
\Delta_\lambda(t^{n-1}a/b|t^n;q,t;p)^{-1}
\]
\end{prop}

We also note the following symmetries that follow from symmetries of the
definition of interpolation theta functions.

\begin{prop}
The interpolation abelian functions satisfy the following symmetry identities:
\begin{align}
R^{*(n)}_\lambda(\dots z_i\dots;a,b;q,t;p)
&=
(\frac{q^2 t^{2n-2} a^2}{b^2})^{|\lambda|}
t^{-4n(\lambda)}
q^{4n(\lambda')}
R^{*(n)}_\lambda(\dots z_i\dots;1/a,1/b;1/q,1/t;p)\\
R^{*(n)}_\lambda(\dots z_i\dots;a,b;q,t;p)
&=
R^{*(n)}_\lambda(\dots -z_i\dots;-a,-b;q,t;p)\label{eq:interp_negate}\\
R^{*(n)}_\lambda(\dots \sqrt{p} z_i\dots;\sqrt{p}a,\sqrt{p}b;q,t;p)
&=
(bt^{1-n}/aq)^{|\lambda|}
t^{2n(\lambda)} q^{-2n(\lambda')}
R^{*(n)}_\lambda(\dots z_i\dots;a,b;q,t;p)\\
R^{*(n)}_\lambda(\dots \sqrt{p} z_i\dots;\sqrt{p}a,b/\sqrt{p};q,t;p)
&=
(pq/ab)^{|\lambda|}
R^{*(n)}_\lambda(\dots z_i\dots;a,b;q,t;p)\\
R^{*(n)}_{m^n+\lambda}(\dots z_i\dots;a,b;q,t;p)
&=
\prod_{1\le i\le n}
\frac{\theta(a z_i^{\pm 1};q;p)_m}
     {\theta(\frac{pq}{b} z_i^{\pm 1};q;p)_m}
R^{*(n)}_\lambda(\dots z_i\dots;aq^m,b/q^m;q,t;p)\\
R^{*(n+k)}_\lambda(\dots z_i\dots,a,at,\dots, at^{k-1};a,b;q,t;p)
&=
\frac{C^0_{\lambda}(t^n,pq a/bt;q,t;p)}
     {C^0_{\lambda}(t^{n+k},pq t^k a/bt;q,t;p)}
R^{*(n)}_\lambda(\dots z_i\dots;t^k a,b;q,t;p)\label{eq:specbranch}
\end{align}
\end{prop}

\section{Binomial coefficients and hypergeometric identities}
\label{sec:hyperg1}

\begin{defn}
The {\em generalized binomial coefficients} are meromorphic functions of
$a,b,q,t$ defined by
\begin{align}
\binomE{\lambda}{\mu}_{[a,b];q,t;p} &:=
\Delta_\mu(\frac{a}{b}|t^n,1/b;q,t;p) R^{*(n)}_\mu(\dots \sqrt{a}
q^{\lambda_i} t^{1-i}\dots;t^{1-n} \sqrt{a},b/\sqrt{a};q,t;p),
\end{align}
for any integer $n\ge \ell(\lambda),\ell(\mu)$.
\end{defn}

Note that by equation \eqref{eq:interp_negate} this is indeed independent of
the choice of $\sqrt{a}$.  While for many purposes the given normalization
of the binomial coefficients is the nicest, it does have the significant
drawback of being singular at points of the form $b=q^{-k} t^l$ for
integers $k,l\ge 0$.  With this in mind, we introduce a second
normalization:
\[
\obinomE{\lambda}{\mu}_{[a,b](v_1,\dots, v_k);q,t;p}
:=
\frac{\Delta^0_\lambda(a|b,v_1,\dots, v_k;q,t;p)}
{\Delta^0_\mu(a/b|1/b,v_1,\dots, v_k;q,t;p)}
\binomE{\lambda}{\mu}_{[a,b];q,t;p}.
\]
(We will omit the parentheses when $k=0$.)
Thus for instance for $b=1$ we find
\[
\obinomE{\lambda}{\mu}_{[a,1](v_1,\dots, v_k);q,t;p} = \delta_{\lambda\mu}.
\]
The standard normalization satisfies particularly nice transformation laws:
\begin{align}
\binomE{\lambda}{\mu}_{[pa,b];q,t;p}
&=
\binomE{\lambda}{\mu}_{[a,b];q,t;p}
\\
\binomE{\lambda}{\mu}_{[a,pb];q,t;p}
&=
\binomE{\lambda}{\mu}_{[a,b];q,t;p}\\
\binomE{\lambda}{\mu}_{[1/a,1/b];1/q,1/t;p}
&=
\binomE{\lambda}{\mu}_{[a,b];q,t;p}.
\end{align}
We also note the following special values:
\begin{align}
\binomE{\lambda}{0}_{[a,b];q,t;p}
&=
1\\
\binomE{\lambda}{\lambda}_{[a,b];q,t;p}
&=
\frac{C^+_\lambda(a;q,t;p)\Delta^0_\lambda(a/b|1/b;q,t;p)}{C^+_\lambda(\frac{a}{b};q,t;p)\Delta^0_\lambda(a|b;q,t;p)}\\
\binomE{m^n}{\lambda}_{[a,b];q,t;p}
&=
\Delta_\lambda(\frac{a}{b}|t^n,q^{-m},t^{1-n}q^m a,1/b;q,t;p)
\end{align}
Furthermore, we readily obtain the following symmetry, from the
corresponding symmetry of interpolation functions:
\[
\frac{\binomE{m^n+\lambda}{m^n+\mu}_{[a,b];q,t;p}}
{\binomE{m^n}{m^n}_{[a,b];q,t;p}}
=
\frac{\Delta^0_\lambda(q^{2m}a|b,pqaq^m/b,pqt^{n-1}q^m,t^{1-n}q^{2m}a;q,t;p)}
     {\Delta^0_\mu(q^{2m}a/b|1/b,pqaq^m,pqt^{n-1}q^m,t^{1-n}q^{2m}a/b;q,t;p)}
\binomE{\lambda}{\mu}_{[q^{2m}a,b];q,t;p}.
\]

As announced in \cite{xforms}, the main summation identity for elliptic
generalized binomial coefficients is the following identity; as the proof
below suggests, it can be thought of as a ``bulk'' difference equation, ala
the bulk Pieri identities and branching rule of \cite{bcpoly}.

\begin{thm}\label{thm:theta_bde}
For any partitions $\kappa\subset\lambda$, and generic parameters
$a,b,c,d,e,q,t\in \C^*$ such that $bcde=apq$,
\[
\binomE{\lambda}{\kappa}_{[a,c];q,t;p}
=
\frac{\Delta^0_\kappa(a/c|1/c,bd,be,pqa/b;q,t;p)}
     {\Delta^0_\lambda(a|c,bd,be,pqa/b;q,t;p)}
\sum_{\kappa\subset\mu\subset\lambda}
\Delta^0_\mu(a/b|c/b,pqa,d,e;q,t;p)
\binomE{\lambda}{\mu}_{[a,b];q,t;p}
\binomE{\mu}{\kappa}_{[a/b,c/b];q,t;p}.
\label{eq:theta_bde}
\]
\end{thm}

\begin{proof}
If we take the difference equation \eqref{eq:diff_eq} for interpolation
abelian functions and rewrite it in terms of binomial coefficients, we
obtain an identity of the form:
\begin{align}
{\binomE{\lambda}{\kappa}\!}_{[a^2,ab];q,t;p}
={}&
\frac{\Delta^0_\kappa(a/b|1/ab,ac,ad,pq^2a^2;q,t;p)}
     {\Delta^0_\lambda(a^2|ab,ac,ad;q,t;p)}\\
&\sum_{\mu\prec\lambda}\!
\Delta^0_\mu(qa^2|qab,qac,qad;q,t;p)
F_{\lambda/\mu}(a;q,t;p)
{\binomE{\mu}{\kappa}\!}_{[a^2q,abq];q,t;p}\notag
\end{align}
with $abcd=p$; the coefficients $F_{\lambda/\mu}(a;q,t;p)$ are
complicated, but are in principle explicit.  Extending $F$ to be 0 if
$\mu\not\prec\lambda$, we can replace the range of summation with
$\kappa\subset\mu\subset\lambda$ without changing the sum.  Now, in the
special case $b=1/aq$ (and thus $cd=pq$), the binomial coefficient on the
right-hand side becomes a delta function, and we thus obtain the identity
\[
F_{\lambda/\kappa}(a;q,t;p)
=
\obinomE{\lambda}{\kappa}_{[a^2,1/q];q,t;p}.
\]
Substituting this in gives the special case $b=1/q$ of the theorem.

We next observe that the values of $b$ for which the theorem holds is
closed under multiplication.  Indeed, if $b=b_1b_2$ such that the theorem
holds for $b=b_1$ and $b=b_2$, then we in particular find using the case
$b=b_1$ of the theorem that
\begin{align}
\binomE{\lambda}{\mu}_{[a,b_1b_2];q,t;p}
&=
\frac{\Delta^0_\kappa(a/b_1b_2|1/b_1b_2,pqa/b_1,b_1d,ce;q,t;p)}
     {\Delta^0_\lambda(a|b_1b_2,pqa/b_1,b_1d,ce;q,t;p)}\\
&\phantom{{}={}}\sum_{\kappa\subset\mu\subset\lambda}
\Delta^0_\mu(a/b_1|b_2,pqa,d,ce/b_1;q,t;p)
\binomE{\lambda}{\mu}_{[a,b_1];q,t;p}
\binomE{\mu}{\kappa}_{[a/b_1,b_2];q,t;p}.\notag
\end{align}
Substituting this in to the right-hand side, we find that the resulting
double sum can be summed over $\mu$ using the case $b=b_2$ of the theorem,
and then summed over $\nu$ using the case $b=b_1$ again.

We thus conclude that the theorem holds whenever $b=q^{-k}$ for $k$ a
positive integer.  Since for $a,c,d$ fixed, both sides are meromorphic
theta functions in $b$ with the same multiplier, the result for general
$b$ follows by analytic continuation.
\end{proof}

\begin{rems}
This identity includes a number of known elliptic hypergeometric identities
as special cases.  Taking both $\lambda$ and $\kappa$ to consist of a
single part gives the Frenkel-Turaev summation for univariate elliptic
hypergeometric series \cite{FrenkelIB/TuraevVG:1997} (the elliptic analogue
of the ${}_8\phi_7$ summation due to Jackson).  Taking $\lambda$ to be a
rectangle and $\kappa=0$ gives an identity conjectured by Warnaar
\cite[Corollary 6.2]{WarnaarSO:2002} (also proved in
\cite{RosengrenH:2001}), which in our notation reads
\[
\sum_{\mu\subset m^n}
\Delta_\mu(a|t^n,q^{-m},b_0,b_1,b_2,b_3;q,t;p)
=
\frac{
C^0_{m^n}(pqa,pqa/b_0b_1,pqa/b_0b_2,pqa/b_1b_2;q,t;p)
}{
C^0_{m^n}(pqa/b_0,pqa/b_1,pqa/b_2,pqa/b_0b_1b_2;q,t;p)
},\label{eq:ip_normalization}
\]
where $q^{-m} t^n b_0b_1b_2b_3 = a^2 pqt$.
\end{rems}

\begin{rems}
When $t=q$, the binomial coefficient is related to Schur-type
interpolation functions, and can thus be expressed as a determinant; to be
precise, we have:
\[
\obinomE{\lambda}{\mu}_{[a,b];q,q;p}
=
\prod_{1\le i<j\le n}
\frac{q^{-\mu_i}\theta(q^{\mu_i-i-\mu_j+j},q^{\mu_i+\mu_j+2-i-j}a/b;p)}
{q^{-\lambda_i}\theta(q^{\lambda_i-i-\lambda_j+j},q^{\lambda_i+\lambda_j+2-i-j}
  a;p)}
\det_{1\le i,j\le n}\left(
  \obinomE{\lambda_i+n-i}{\mu_j+n-j}_{[\frac{a}{q^{2n-2}},b];q,q;p}
\right)
\]
(The given scale factors can be derived via Warnaar's determinant identity,
Lemma 5.3 of \cite{WarnaarSO:2002}; also note that setting $\lambda=\mu$
makes the determinant triangular, so trivial to evaluate.) This gives rise
to Warnaar's Schlosser-type identity (Theorem 5.1 of
\cite{WarnaarSO:2002}), via the case $t=q$, $b=q^{-k}$, $\kappa=0$ above,
in the following way.  When thus specialized, the binomial coefficients
with lower index $\kappa$ are ratios of $C^0$ symbols, while the remaining
binomial coefficient can be expressed as a determinant.  Now, for this
coefficient not to vanish, we must have $\lambda_i-k\le \mu_i\le \lambda_i$
for each $i$; if in addition $\lambda_i-\lambda_{i+1}\ge k$, the relevant
matrix becomes diagonal.  Since the sum is terminating for generic
$\lambda$, we can analytically continue in $q^{\lambda_i}$, $1\le i\le n$
(summing over $\lambda_i-\mu_i$); the desired identity results.
\end{rems}

\begin{cor}
The binomial coefficient $\obinomE{\lambda}{\mu}_{[a,1/q];q,t;p}$ vanishes
unless $\mu\prec\lambda$, in which case
\begin{align}
\obinomE{\lambda}{\mu}_{[a,1/q];q,t;p}
&=
\prod_{\substack{(i,j)\in \lambda\\\lambda_i=\mu_i}}
\frac{\theta(q^{\lambda_i+j-1} t^{2-\lambda'_j-i} a;p)}
     {\theta(q^{\mu_i-j} t^{\mu'_j-i} t;p)}
\prod_{\substack{(i,j)\in \lambda\\\lambda_i\ne \mu_i}}
\frac{\theta(q^{\lambda_i-j} t^{\lambda'_j-i} pq;p)}
     {\theta(q^{\mu_i+j-1} t^{2-\mu'_j-i} apq^2/t;p)}\\
&\phantom{{}={}}\prod_{\substack{(i,j)\in \mu\\\lambda_i=\mu_i}}
\frac{\theta(q^{\lambda_i-j} t^{\lambda'_j-i} t;p)}
     {\theta(q^{\mu_i+j-1} t^{2-\mu'_j-i} aq;p)}
\prod_{\substack{(i,j)\in \mu\\\lambda_i\ne \mu_i}}
\frac{\theta(q^{\lambda_i+j-1} t^{2-\lambda'_j-i} pqa/t;p)}
     {\theta(q^{\mu_i-j} t^{\mu'_j-i} pq;p)}\notag
\end{align}
\end{cor}

\begin{proof}
The proof of Theorem \ref{thm:theta_bde} gives a formula for the binomial
coefficient, which in particular implies that it vanishes when required.
Simplifying this as in the proof of Theorem 4.12 of \cite{bcpoly} gives the
desired result.
\end{proof}

\begin{rem}
Similarly, for any integer $k\ge 0$,
$\obinomE{\lambda}{\mu}_{[a,q^{-k}];q,t;p}$ vanishes unless
$\mu\prec_k\lambda$.
\end{rem}

When $c=1$ in the bulk difference equation \eqref{eq:theta_bde}, the left-hand
side becomes a delta function, and we thus conclude:

\begin{cor}
Elliptic generalized binomial coefficients satisfy the inversion identity
\[
\sum_{\kappa\subset\mu\subset\lambda}
\binomE{\lambda}{\mu}_{[a,b];q,t;p}
\binomE{\mu}{\kappa}_{[a/b,1/b];q,t;p}
=
\delta_{\lambda\kappa}
\]
\end{cor}

From the form of Cauchy interpolation polynomials, the following symmetry
of interpolation theta functions is essentially automatic:
\[
P^{*(m,n)}_\mu(\lambda;a,pq/q^mt^na;q,t;p)
=
P^{*(n,m)}_{n^m-\lambda'}(n^m-\mu';\sqrt{p}/a,\sqrt{p}t^{n-1} q^m a;1/t,1/q;p);
\]
note in particular that the roles of argument and index are interchanged
via this symmetry.  The special form of the bigrid in our case gives rise
to another symmetry as well, which follows from a computation along the
lines of Lemma 2.1 of \cite{bcpoly}:
\begin{align}
P^{*(m,n)}_\mu(\lambda;a,pq/q^mt^na;q,t;p)
&=
\frac{C^0_{m^n}(t^n;q,t;p)}{C^0_{m^n}(q^{-m};q,t;p)}
\frac{C^0_\lambda(q^m t^{2n-1} a^2,(pq/t) t^n;q,t;p)}
     {C^0_\lambda(t^{n-1} a^2,(pq/t) q^{-m};q,t;p)}
\frac{C^0_\mu(t^{n-1} a^2,q^{-m};q,t;p)}
     {C^0_\mu(q^m t^{2n-1} a^2,t^n;q,t;p)}\notag\\
&\phantom{{}={}}
P^{*(n,m)}_{\mu'}(\lambda';t^{n-1/2} q^{m-1/2} a,pq^{3/2} t^{1/2}/a;1/t,1/q;p)
\end{align}
We thus obtain an action of the group $Z_2^2$ on Cauchy interpolation theta
functions; this symmetry, it turns out, actually extends to all
interpolation theta functions.  In terms of binomial coefficients, we
obtain the following two symmetries: ``duality'' and ``complementation''.

\begin{cor} (Duality)
\[
\binomE{\lambda}{\mu}_{[a,b];q,t;p}
=
\binomE{\lambda'}{\mu'}_{[aqt,b];1/t,1/q;p}
\]
\end{cor}

\begin{proof}
When $b=pq/t$, the relevant interpolation functions are of Cauchy type.
We thus find that
\[
\binomE{\lambda}{\mu}_{[a,pq/t];q,t;p}
=
F_\lambda(a;q,t;p) G_\mu(a;q,t;p)
\binomE{\lambda'}{\mu'}_{[aqt,pq/t];1/t,1/q;p}
\]
for appropriate functions $F_\lambda$, $G_\mu$, with $G_0=1$.  Setting
$\mu=0$ shows that $F_\lambda=1$; setting $\mu=\lambda$ then shows that
$G_\mu=1$.  (These could, of course, also be computed directly.)

The bulk difference equation \eqref{eq:theta_bde} is preserved by this
symmetry; it follows that the set of $b$ for which this duality holds is
closed under multiplication, and thus the result follows by analytic
continuation.
\end{proof}

\begin{cor}
The binomial coefficient $\obinomE{\lambda}{\mu}_{[a,t];q,t;p}$ vanishes
unless $\mu\prec'\lambda$, in which case
\begin{align}
\obinomE{\lambda}{\mu}_{[a,t];q,t;p}
&=
\prod_{\substack{(i,j)\in \lambda\\\lambda'_j=\mu'_j}}
\frac{\theta(q^{\lambda_i+j-1} t^{2-\lambda'_j-i} a;p)}
     {\theta(q^{\mu_i-j} t^{\mu'_j-i} pq;p)}
\prod_{\substack{(i,j)\in \lambda\\\lambda'_j\ne \mu'_j}}
\frac{\theta(q^{\lambda_i-j} t^{\lambda'_j-i} t;p)}
     {\theta(q^{\mu_i+j-1} t^{2-\mu'_j-i} apq/t^2;p)}\\
&\phantom{{}={}}\prod_{\substack{(i,j)\in \mu\\\lambda'_j=\mu'_j}}
\frac{\theta(q^{\lambda_i-j}t^{\lambda'_j-i} pq;p)}
     {\theta(q^{\mu_i+j-1}t^{2-\mu'_j-i} a/t;p)}
\prod_{\substack{(i,j)\in \mu\\\lambda'_j\ne \mu'_j}}
\frac{\theta(q^{\lambda_i+j-1} t^{2-\lambda'_j-i} pqa/t;p)}
     {\theta(q^{\mu_i-j} t^{\mu'_j-i} t;p)}\notag
\end{align}
\end{cor}

\begin{rem}
Similarly, $\obinomE{\lambda}{\mu}_{[a,t^k];q,t;p}$ vanishes unless
$\mu\prec'_k\lambda$.
\end{rem}

\begin{cor}
If $\lambda_1,\mu_1\le m$, the binomial coefficients satisfy the identity
\[
\frac{\binomE{m^n\cdot\lambda}{m^n\cdot\mu}_{[a,b];q,t;p}}
{\binomE{m^n}{m^n}_{[a,b];q,t;p}}
=
\frac{\Delta^0_\lambda(\frac{a}{t^{2n}}|b,\frac{pq a}{t^n b},\frac{pq}{q^m t^{n+1}},\frac{q^m a}{t^{2n-1}};q,t;p)}
     {\Delta^0_\mu(\frac{a}{t^{2n}b}|\frac{1}{b},\frac{pq
         a}{t^n},\frac{pq}{q^m t^{n+1}},\frac{q^m a}{t^{2n-1} b};q,t;p)}
\binomE{\lambda}{\mu}_{[t^{-2n}a,b];q,t;p}.
\]
\end{cor}

Similarly,

\begin{cor} (Complementation)
\[
\frac{\binomE{\lambda}{\mu}_{[a,b];q,t;p}}
{\Delta_{m^n}(a|t^n,q^{-m},b,q^m t^{1-n}a/b;q,t;p)}
=
\frac{\Delta_\mu(a/b|t^n,q^{-m},1/b,q^m t^{1-n}a;q,t;p)}
     {\Delta_\lambda(a|t^n,q^{-m},b,q^m t^{1-n}a/b;q,t;p)}
\binomE{m^n-\mu}{m^n-\lambda}_{[q^{-2m} t^{2n-2}b/a,b];q,t;p}.
\]
\end{cor}

If $b=q^{-k}$ or $b=t^k$ there is a further symmetry:

\begin{cor}\label{cor:binom_spec_recur}
If $\ell(\lambda)\le n$,
\[
\obinomE{k^n+\lambda}{\mu}_{[a,q^{-k}];q,t;p}
=
\frac{C^0_{k^n}(q^{-k};q,t;p)}{C^0_{k^n}(pqa q^k;q,t;p)}
\frac{\Delta_\mu(q^k a|t^n,t^{1-n}q^k a;q,t;p)}
     {\Delta_\lambda(q^{2k} a|t^n,t^{1-n}q^k a;q,t;p)}
\obinomE{\mu}{\lambda}_{[q^k a,q^{-k}];q,t;p}.
\]
If $\lambda_1\le m$,
\[
\obinomE{m^k\cdot\lambda}{\mu}_{[a,t^k];q,t;p}
=
\frac{C^0_{m^k}(t^k;q,t;p)}{C^0_{m^k}(pq a t^{-k};q,t;p)}
\frac{\Delta_\mu(t^{-k} a|q^{-m},q^{m}t^{1-k} a;q,t;p)}
     {\Delta_\lambda(t^{-2k} a|q^{-m},q^{m}t^{1-k} a;q,t;p)}
\obinomE{\mu}{\lambda}_{[t^{-k} a,t^k];q,t;p}.
\]
\end{cor}

\begin{proof}
If we express the binomial coefficient on the left of the first equation in
terms of interpolation theta functions, we find that complementing $\mu$
gives an interpolation theta function that is still evaluated at a
partition, namely $m^n-\lambda$.  Re-expressing the result as a binomial
coefficient and applying complementation symmetry gives the desired result.
The second equation follows by duality.
\end{proof}

Since the bulk difference equation \eqref{eq:theta_bde} can be viewed as a
multivariate Jackson summation indexed over a skew Young diagram, we also
expect there to be a corresponding analogue of the Bailey transformation.

\begin{thm}\label{thm:theta_bailey} The sum
\[
\frac{\Delta^0_\lambda(a|b,apq/bf;q,t;p)}
     {\Delta^0_\kappa(a/c|b/c,apq/bd;q,t;p)}
\sum_{\kappa\subset\mu\subset\lambda}
\frac{\Delta^0_\mu(a/b|c/b,f,g;q,t;p)}
     {\Delta^0_\mu(a/b|1/b,d,e;q,t;p)}
\binomE{\lambda}{\mu}_{[a,b];q,t;p}
\binomE{\mu}{\kappa}_{[a/b,c/b];q,t;p}
\label{eq:theta_bailey}
\]
is symmetric in $b$ and $b'$, where $bb'de=capq$, $bb'fg=apq$.
\end{thm}

\begin{proof}
Summing
\begin{align}
\frac{\Delta^0_\lambda(a|b',apq/b'f;q,t;p)}
{\Delta^0_\kappa(a/c|b/c,apq/bd;q,t;p)}
\sum_{\mu,\nu}
\frac{\binomE{\lambda}{\nu}_{[a,b'];q,t;p}
\binomE{\nu}{\mu}_{[a/b',b/b'];q,t;p}
\binomE{\mu}{\kappa}_{[a/b,c/b];q,t;p}
}
{\Delta^0_\nu(a/b'|1/b',apq/b,apq/b'f,apq/b'g;q,t;p)
 \Delta^0_\mu(a/b|b'/b,pqa/c,d,e;q,t;p)}
\end{align}
over $\nu$ gives \eqref{eq:theta_bailey}; summing over
$\mu$ gives the same sum with $b$ replaced by $b'$.
\end{proof}

\begin{rems}
The bulk difference equation \eqref{eq:theta_bde} is the special case $c=b$
above.
\end{rems}

\begin{rems}
The rectangular case $\lambda=m^n$, $\mu=0$ gives Warnaar's conjectured
multivariate elliptic Bailey transformation, Conjecture 6.1 of
\cite{WarnaarSO:2002}, as a corollary.  One can also obtain a
transformation of Schlosser-type sums by setting $t=q$, $b=q^{-N_1}$,
$b'=q^{-N_2}$, and $\kappa=0$.
\end{rems}

This induces a symmetry under the Weyl group $D_4$.

\begin{thm}
Define a function $\Omega_{\lambda/\kappa}(a,b{:}v_0,v_1,v_2,v_3;q,t;p)$ by
\begin{align}
\Omega_{\lambda/\kappa}&(a,b{:}v_0,v_1,v_2,v_3;q,t;p)\\
&:=
\sum_{\kappa\subset\mu\subset\kappa}
\frac{
C^0_{\lambda/\mu}(pqa/v_0,pqa/v_1,pqa/v_2,pqa/v_3;q,t;p)}
{C^0_{\mu/\kappa}(v_0/c,v_1/c,v_2/c,v_3/c;q,t;p)^{-1}}
\obinomE{\lambda}{\mu}_{[pqa^2,pqac];q,t;p}
\obinomE{\mu}{\kappa}_{[a/c,b/c];q,t;p},\notag
\end{align}
where $c=\sqrt{bv_0v_1v_2v_3/apq}$.  Then $\Omega_{\lambda/\kappa}$ is
invariant under the natural action of $D_4$ on $(v_0,v_1,v_2,v_3)$; in
particular,
\[
\Omega_{\lambda/\kappa}(a,b{:}v_0,v_1,v_2,v_3;q,t;p)
=
\Omega_{\lambda/\kappa}(a,b{:}v_0,v_1,1/v_2,1/v_3;q,t;p).
\]
\end{thm}

\begin{proof}
Indeed, the symmetry
\[
\Omega_{\lambda/\kappa}(a,b{:}v_0,v_1,v_2,v_3;q,t;p)
=
\Omega_{\lambda/\kappa}(a,b{:}v_0,v_1,1/v_2,1/v_3;q,t;p)
\]
is simply Theorem \ref{thm:theta_bailey}, up to a change of variables;
since $\Omega_{\lambda/\kappa}$ has a manifest $S_4$ symmetry, the $D_4$
symmetry follows immediately.
\end{proof}

From the center of $D_4$, we obtain the following transformation,
generalizing the commutation relation of Corollary
\ref{cor:quasicommutation} (which corresponds to the case $b=b'=1/q$ when
evaluated at a partition).

\begin{cor}\label{cor:centralD4}
If $bcd=b'c'd'$, then
\[
\frac{\Delta^0_\lambda(a|b,c,d;q,t;p)}
{\Delta^0_\kappa(\frac{a}{bb'}|\frac{pqa}{b'},\frac{pqa}{bc},\frac{pqa}{bd};q,t;p)}
\sum_\mu
\frac{
\binomE{\lambda}{\mu}_{[a,b];q,t;p}
\binomE{\mu}{\kappa}_{[a/b,b'];q,t;p}}
{
\Delta^0_{\mu}(a/b|\frac{pqa}{bb'},\frac{pqa}{c'},\frac{pqa}{d'},1/b,c,d;q,t;p)
}
\]
is invariant under $(b,c,d,b',c',d')\mapsto (b',c',d',b,c,d)$.
\end{cor}

The case $b=q^{-1}$ gives rise to the following difference equation:

\begin{thm}\label{thm:raise_op_interp}
Define a difference operator $D^{+(n)}_q(u_0{:}u_1{:}u_2,u_3,u_4;t;p)$ by
\begin{align}
(D^{+(n)}_q(u_0{:}u_1{:}u_2,u_3,u_4;t,p)f)&(\dots z_i\dots)\notag\\
&:=
\prod_{1\le i\le n}
\frac{\theta(pq t^{n-i} u_1/u_0;p)}
{\prod_{2\le r\le 5} \theta(u_r t^{n-i} u_1;p)}\\
&\phantom{{}:={}}
\prod_{1\le i\le n} (1+R_{z_i})
\prod_{1\le i\le n}
\frac{\prod_{1\le r\le 5} \theta(u_r z_i;p)}
     {\theta(pq z_i/u_0;p)\theta(z_i^2;p)}
\prod_{1\le i<j\le n}
\frac{\theta(t z_i z_j;p)}
     {\theta(z_i z_j;p)}
f(\dots q^{1/2} z_i\dots),\notag
\end{align}
where $u_5=p^2q/t^{n-1}u_0u_1u_2u_3u_4$, and $R_{z_i}$ acts on functions by
replacing $z_i$ by its reciprocal.  Then
\[
D^{+(n)}_q(u_0{:}u_1{:}u_2,u_3,u_4;t;p)
R^{*(n)}_\lambda(;q^{1/2}u_1,q^{-1/2}u_0;q,t;p)
=
\sum_{\lambda\subset\mu\subset \lambda+1^n}
c_{\lambda\mu}
R^{*(n)}_{\mu}(;u_1,u_0;q,t;p)
\]
where
\[
c_{\lambda\mu}
=
\frac{\Delta_\mu    (  t^{n-1}u_1/u_0|t^n,\frac{pq}{u_0u_2},\frac{pq}{u_0u_3},\frac{pq}{u_0u_4},\frac{pq}{u_0u_5};q,t;p)}
     {\Delta_\lambda(q t^{n-1}u_1/u_0|t^n,\frac{pq}{u_0u_2},\frac{pq}{u_0u_3},\frac{pq}{u_0u_4},\frac{pq}{u_0u_5};q,t;p)}
\obinomE{\mu}{\lambda}_{[t^{n-1}u_1/u_0,1/q];q,t;p}
\]
\end{thm}

One can obtain other identities from our main identities
\eqref{eq:theta_bde} and Theorem \eqref{thm:theta_bailey} by first applying
some combination of duality and complementation symmetry, specializing the
result so that the sum terminates regardless of $\lambda$ and/or $\kappa$,
then analytically continuing.

For instance:

\begin{thm}
The interpolation theta functions satisfy the following connection
coefficient identity:
\begin{align}
[P^{*(m,n)}_\mu(;a,b';q,t;p)]&P^{*(m,n)}_\lambda(;a,b;q,t;p)\\
&=
\frac{C^0_{m^n}(t^{n-1}ab,\frac{b}{a};q,t;p)}
     {C^0_{m^n}(t^{n-1}ab',\frac{b'}{a};q,t;p)}
\frac{\Delta_\mu(\frac{t^{n-1}a}{q^m b'}|\frac{b}{b'},t^n,q^{-m},\frac{pq}{q^m b b'};q,t;p)}
     {\Delta_\lambda(\frac{t^{n-1}a}{q^m b}|\frac{b'}{b},t^n,q^{-m},\frac{pq}{q^m b b'};q,t;p)}
\binomE{\mu}{\lambda}_{[\frac{t^{n-1}a}{q^m b'},\frac{b}{b'}];q,t;p}\notag
\end{align}
\end{thm}

\begin{proof}
If we set $e=q^{-m}$ in the bulk difference equation \eqref{eq:theta_bde},
the resulting sum will vanish for $\mu\not\subset m^n$, regardless of
$\lambda$, and we may thus replace the range of summation with $\mu\subset
m^n$.  If we express the binomial coefficients with upper index $\lambda$
in terms of interpolation theta functions, we find that the remaining
factors involving $\lambda$ are the same on both sides; dividing by those
factors and analytically continuing gives the above result.
\end{proof}

Applying Corollary \ref{cor:comp_symm_P} gives the following.

\begin{cor}
We have the following connection coefficient identity:
\[
[P^{*(m,n)}_\mu(;a',b;q,t;p)]
P^{*(m,n)}_{\lambda}(;a,b;q,t;p)
=
\obinomE{\lambda}{\mu}_{[\frac{t^{n-1}a}{q^m b},\frac{a}{a'}](t^{n-1}aa');q,t;p}
\]
Similarly,
\[
[R^{*(n)}_\mu(;a',b;q,t;p)]
R^{*(n)}_{\lambda}(;a,b;q,t;p)
=
\obinomE{\lambda}{\mu}_{[\frac{t^{n-1}a}{b},\frac{a}{a'}](t^{n-1}aa');q,t;p}
\label{eq:conn_Rs}
\]
\end{cor}

Combining the two identities for interpolation theta functions gives the
following result.

\begin{thm}\label{thm:multi6j}
General connection coefficients for interpolation theta functions are
given by the sum
\begin{align}
[P^{*(m,n)}_\kappa(;a',b';q,t;p)]P^{*(m,n)}_\lambda(;a,b;q,t;p)
&=
\frac{C^0_{m^n}(t^{n-1}a'b,\frac{b}{a'};q,t;p)}
     {C^0_{m^n}(t^{n-1}a'b',\frac{b'}{a'};q,t;p)}
\frac{
\Delta_\kappa(\frac{t^{n-1}a'}{q^m b'}|t^n,q^{-m},\frac{pq}{q^m b
  b'};q,t;p)}
{\Delta^{0}_\lambda(\frac{t^{n-1}a}{q^m b}|\frac{pq}{q^ma'b};q,t;p)}
\notag
\\
&\phantom{{}={}}
\sum_{\mu\subset\lambda,\kappa}
\frac{
\obinomE{\lambda}{\mu}_{[\frac{t^{n-1}a}{q^m b},\frac{a}{a'}];q,t;p}
\obinomE{\kappa}{\mu}_{[\frac{t^{n-1}a'}{q^m b'},\frac{b}{b'}];q,t;p}}
{\Delta_\mu(\frac{t^{n-1}a'}{q^m b}|t^n,q^{-m},t^{n-1}aa',\frac{pq}{q^m b b'};q,t;p)}
\notag
\end{align}
\end{thm}

\begin{rem}
These connection coefficients generalize Rosengren's construction of
elliptic 6-j symbols in \cite{RosengrenH:2003b}.
\end{rem}

Using the connection coefficient formula \eqref{eq:conn_Rs}, we can extend
the special branching rule \eqref{eq:specbranch} to a general branching
rule.

\begin{thm}\label{thm:Rs_branch}
Interpolation functions satisfy the following ``bulk'' branching rule.
\[
R^{*(n+k)}_{\lambda}(\dots z_i\dots,v,vt,\dots vt^{k-1};a,b;q,t;p)
=
\sum_\kappa c_{\lambda\kappa} R^{*(n)}_\kappa(\dots z_i\dots;a,b;q,t;p),
\]
where
\[
c_{\lambda\kappa}=
\obinomE{\lambda}{\kappa}_{[t^{n+k-1}a/b,t^k](t^{n+k-1}av,t^na/v,pqa/tb);q,t;p}.
\]
\end{thm}

\begin{proof}
Expand the left-hand side in interpolation functions with parameters
$(v,b)$, rewrite those as $n$-variable interpolation functions, and
convert the parameters back.  This gives the coefficients of the branching
rule as a sum over partitions, which can be summed via the bulk difference
equation \eqref{eq:theta_bde}.
\end{proof}

Similarly:

\begin{thm}
The interpolation functions satisfy the following generalized Pieri identity:
\[
\prod_{1\le i\le n}
\frac{\theta(v z_i^{\pm 1};q;p)_m}
     {\theta((pq/b) z_i^{\pm 1};q;p)_m}
R^{*(n)}_{\lambda}(\dots z_i\dots;a,q^{-m} b;q,t;p)
=
\sum_\kappa c_{\kappa\lambda} R^{*(n)}_\kappa(\dots z_i\dots;a,b;q,t;p),
\]
where
\[
c_{\kappa\lambda}
=
\Delta^0_{m^n}(t^{n-1}v/b|t^{n-1}va,v/a;q,t;p)
\frac{\Delta_\kappa(t^{n-1}a/b|t^n,pq/vb,q^mv/b;q,t;p)}
     {\Delta_\lambda(q^m t^{n-1}a/b|t^n,pq/vb,q^mv/b;q,t;p)}
\obinomE{\kappa}{\lambda}_{[t^{n-1}a/b,q^{-m}];q,t;p}.
\]
\end{thm}

If we expand a Cauchy-type interpolation theta function in other
interpolation theta functions, the result can be analytically continued
again; we obtain the following result.

\begin{thm}\label{thm:interp_cauchy}
Interpolation theta functions satisfy the following ``Cauchy'' identity.
\begin{align}
\prod_{1\le i\le n,1\le j\le m} \theta(y_j x_i,y_j/x_i;p)
=
\sum_{\mu\subset m^n}
\frac{\Delta_\mu(\frac{t^{n-1}a}{q^m b}|t^n,q^{-m};q,t;p)}
     {C^0_{m^n}(t^{n-1}ab,b/a;q,t;p)}
&P^{*(m,n)}_\mu(\dots x_i\dots;a,b;q,t;p)\\
&P^{*(n,m)}_{n^m-\mu'}(\dots \sqrt{p}/y_i\dots;\sqrt{p}/a,\sqrt{p}/b;1/t,1/q;p)\notag
\end{align}
\end{thm}

\section{Biorthogonal functions}

\begin{defn}
Let $t_0$, $t_1$, $t_2$, $t_3$, $u_0$, $u_1$, $q$, $t$ be parameters such
that $t^{2n-2}t_0t_1t_2t_3u_0u_1=pq$.  Then define
\[
\tilde{R}^{(n)}_{\lambda}(;t_0{:}t_1,t_2,t_3;u_0,u_1;q,t;p)
:=
\sum_{\mu\subset\lambda}
\frac{
\binomE{\lambda}{\mu}_{[\frac{1}{u_0u_1},\frac{1}{t^{n-1}t_0u_1}];q,t;p}
R^{*(n)}_\mu(;t_0,u_0;q,t;p)}
{\Delta^{0}_\mu(\frac{t^{n-1}t_0}{u_0}|t^{n-1}t_0t_1,t^{n-1}t_0t_2,t^{n-1}t_0t_3,t^{n-1}t_0u_1;q,t;p)}.
\label{eq:binomial_formula}
\]
\end{defn}

\begin{rem}
This is an analogue of the ``binomial formula'' of \cite{OkounkovA:1998a}.
\end{rem}

The specific form of the above expansion is effectively determined by the
requirement that $\tilde{R}^{(n)}_\lambda$ be symmetrical in $t_0$ through
$t_3$; more precisely:

\begin{thm}
The function $\tilde{R}^{(n)}_\lambda$ satisfies the symmetry
\[
\tilde{R}^{(n)}_{\lambda}(;t_1{:}t_0,t_2,t_3;u_0,u_1;q,t;p)
=
\frac{\tilde{R}^{(n)}_{\lambda}(;t_0{:}t_1,t_2,t_3;u_0,u_1;q,t;p)}
     {\tilde{R}^{(n)}_{\lambda}(\dots t^{n-i}
       t_1\dots;t_0{:}t_1,t_2,t_3;u_0,u_1;q,t;p)},
\]
where
\[
\tilde{R}^{(n)}_{\lambda}(\dots t^{n-i}
       t_1\dots;t_0{:}t_1,t_2,t_3;u_0,u_1;q,t;p)
=
\Delta^0_\lambda(1/u_0u_1|t^{n-1}t_1t_2,t^{n-1}t_1t_3,t^{1-n}/t_1u_1,pqt^{n-1}t_0/u_0;q,t;p).
\]
\end{thm}

\begin{proof}
If we expand 
\[
\tilde{R}^{(n)}_{\lambda}(;t_0{:}t_1,t_2,t_3;u_0,u_1;q,t;p)
\]
in terms of $R^{*(n)}_\lambda(;t_1,u_0;q,t;p)$ using the connection
coefficient formula \eqref{eq:conn_Rs}, the result simplifies via the bulk
difference equation \eqref{eq:theta_bde} to give the desired identity.
\end{proof}

\begin{rem}
Alternatively, one can apply the same steps used above to derive the
connection coefficient identity, but starting with Theorem
\ref{thm:theta_bailey} instead of \eqref{eq:theta_bde}.  The resulting
computation is somewhat more complicated, but does have the merit of {\em
  deriving} the expansion above.
\end{rem}

In particular, we have:
\[
\tilde{R}^{(n)}_\lambda(;t_0{:}t_1,t_2,t_3;u_0,\frac{1}{t^{n-1}t_1};q,t;p)
=
\frac{R^{*(n)}_\lambda(;t_1,u_0;q,t;p)}
{\Delta^0_\lambda(t^{n-1}t_1/u_0|t^{n-1}t_0t_1,t_1/t_0;q,t;p)}.
\]

In addition, from the difference equation \eqref{eq:diff_eq_Rs} for
interpolation abelian functions, we directly obtain the following
difference equation for $\tilde{R}^{(n)}_\lambda$.

\begin{lem}
The function $\tilde{R}^{(n)}_\lambda$ satisfies the difference equation
\begin{align}
D^{(n)}(u_0,t_0,t_1,t^{1-n}p/u_0t_0t_1;q,t;p)
&\tilde{R}^{(n)}_{\lambda}(;q^{1/2}t_0{:}q^{1/2}t_1,q^{-1/2}t_2,q^{-1/2}t_3;q^{1/2}u_0,q^{-1/2}u_1;q,t;p)\\
&=
\prod_{1\le i\le n} \theta(t^{n-i} u_0t_0,t^{n-i} u_0t_1,t^{n-i} t_0t_1;p)
\tilde{R}^{(n)}_{\lambda}(;t_0{:}t_1,t_2,t_3;u_0,u_1;q,t;p)
\notag
\end{align}
\end{lem}

\begin{thm}
The functions $\tilde{R}^{(n)}_\lambda$ agree with the biorthogonal
functions of \cite{xforms}.  To be precise, in the notation of that paper,
\[
\tilde{\cal R}^{(n)}_{\lambda,\mu}(;t_0{:}t_1,t_2,t_3;u_0,u_1;t;p,q)
=
\tilde{R}^{(n)}_\lambda(;t_0{:}t_1,t_2,t_3;u_0,u_1;p,t;q)
\tilde{R}^{(n)}_\mu(;t_0{:}t_1,t_2,t_3;u_0,u_1;q,t;p).
\]
\end{thm}

\begin{proof}
We need to show
\[
\tilde{\cal R}^{(n)}_{0,\lambda}(;t_0{:}t_1,t_2,t_3;u_0,u_1;t;p,q)
=
\tilde{R}^{(n)}_\lambda(;t_0{:}t_1,t_2,t_3;u_0,u_1;q,t;p).
\]
This equation is preserved by the action of the difference operators
\[
D^{(n)}(u_0,t_r,t_s,t^{1-n}p/u_0t_rt_s;q,t;p),\quad 0\le r<s\le 3;
\]
in particular, both sides satisfy the same generalized eigenvalue equation
with respect to the compositions
\begin{align}
&D^{(n)}(u_0,t_0,t_1;q,t;p)
D^{(n)}(q^{1/2} u_0,q^{-1/2} t_2,q^{-1/2} t_3;q,t;p)
,\\
&D^{(n)}(u_0,t_0,t_2;q,t;p)
D^{(n)}(q^{1/2} u_0,q^{-1/2} t_1,q^{-1/2} t_3;q,t;p)
,\dots
\end{align}
and thus agree up to scalar multiples.  Since both sides evaluate to 1 at
$\dots t_0 t^{n-i}\dots$, they agree everywhere.
\end{proof}

\begin{rem}
Similarly, one can show using Theorem \ref{thm:raise_op_interp} that the
difference operator defined there acts as a raising operator on the
biorthogonal functions, just as in \cite{xforms}.
\end{rem}

In \cite{xforms}, it was shown that these functions are biorthogonal with
respect to an appropriate contour integral, and a number of other
properties were given.  A few remaining properties were outside the scope
of that paper; we are now in a position to prove these.  Most striking of
these is evaluation symmetry, generalizing the analogous property of
Macdonald \cite{MacdonaldIG:1995} and Koornwinder
\cite{KoornwinderTH:1992,vanDiejenJF:1996,SahiS:1999} polynomials (see also
\cite{bcpoly} for a proof for Koornwinder polynomials along the present
lines).

\begin{thm}\label{thm:biorth_ev_symm}
For otherwise generic parameters satisfying $t^{2n-2}t_0t_1t_2t_3u_0u_1=pq$,
\[
\tilde{R}^{(n)}_{\lambda}(\dots t_0 t^{n-i} q^{\kappa_i}\dots;t_0{:}t_1,t_2,t_3;u_0,u_1;q,t;p)
=
\tilde{R}^{(n)}_{\kappa}(\dots \hat{t}_0 t^{n-i} q^{\lambda_i}\dots;\hat{t}_0{:}\hat{t}_1,\hat{t}_2,\hat{t}_3;\hat{u}_0,\hat{u}_1;q,t;p),
\]
where
\[
\hat{t}_0 = \sqrt{t_0t_1t_2t_3/pq}\quad
\hat{t}_0\hat{t}_1=t_0t_1\quad
\hat{t}_0\hat{t}_2=t_0t_2\quad
\hat{t}_0\hat{t}_3=t_0t_3\quad
\frac{\hat{u}_0}{\hat{t}_0} = \frac{u_0}{t_0}\quad
\frac{\hat{u}_1}{\hat{t}_0} = \frac{u_1}{t_0}.
\]
\end{thm}

\begin{proof}
Upon specializing the variables as required, the result can be expressed
as a sum over binomial coefficients; changing from the original
parameters to the ``hatted'' parameters interchanges the binomial coefficients.
\end{proof}

\begin{rems}
Aside from some simple factors, this sum over binomial coefficients is the
analytic continuation of the ``multivariate 6-j symbol'' of Theorem
\ref{thm:multi6j}.
\end{rems}

\begin{rems}
When the biorthogonal function is specialized to an interpolation
function, we obtain the identity
\[
\frac{R^{*(n)}_\mu(\dots (v/a) t^{n-i} q^{\lambda_i}\dots;a,b/a;q,t;p)}
     {R^{*(n)}_\mu(\dots (v/a) t^{n-i}\dots;a,b/a;q,t;p)}
=
\frac{R^{*(n)}_\lambda(\dots (v/a') t^{n-i} q^{\mu_i}\dots;a',b/a';q,t;p)}
     {R^{*(n)}_\lambda(\dots (v/a') t^{n-i}\dots;a',b/a';q,t;p)},
\]
where
\[
a' = \sqrt{t^{n-1}b}v/a.
\]
If $v=q^{-m}t^{1-n}$ here, we recover the complementation symmetry of
binomial coefficients; by analytic continuation, the two results are
equivalent.
\end{rems}

\begin{rems}
In the univariate case, this was proved in Section 9 of
\cite{SpiridonovVP/ZhedanovAS:2000b}.
\end{rems}

Using inversion of binomial coefficients, one can expand interpolation
functions in terms of biorthogonal functions.

\begin{thm}
Biorthogonal functions satisfy the following ``inverse binomial formula''.
\begin{align}
[\tilde{R}^{(n)}_{\mu}(;t_0{:}t_1,t_2,t_3;u_0,u_1;q,t;p)]R^{*(n)}_\lambda(;t_0,u_0;q,t;p)
=
\frac{\binomE{\lambda}{\mu}_{[t^{n-1}t_0/u_0,t^{n-1}t_0u_1];q,t;p}}
{\Delta^{0}_\lambda(\frac{t^{n-1}t_0}{u_0}|\frac{pq}{u_0t_1},\frac{pq}{u_0t_2},\frac{pq}{u_0t_3},\frac{pq}{u_0u_1};q,t;p)}
\end{align}
\end{thm}

Combined with the binomial formula \eqref{eq:binomial_formula}, we obtain
connection coefficient formulas, analogous to connection coefficients for
Askey-Wilson polynomials \cite{AskeyR/WilsonJ:1985}.

\begin{thm}\label{thm:conn_biorth1}
If $t^{2n-2}t_0t_1t_2t_3u_0u_1=pq$ and $t'_1t'_2t'_3u'_1=t_1t_2t_3u_1$,
then
\begin{align}
[\tilde{R}^{(n)}_{\kappa}(;t_0{:}t'_1,t'_2,t'_3;u_0,u'_1;q,t;p)]
\tilde{R}^{(n)}_{\lambda}(;t_0{:}t_1,t_2,t_3;&u_0,u_1;q,t;p)\\
=
\sum_{\kappa\subset \mu\subset\lambda}&
\frac{\Delta^{0}_\mu(\frac{t^{n-1}t_0}{u_0}|t^{n-1}t_0t'_1,t^{n-1}t_0t'_2,t^{n-1}t_0t'_3,t^{n-1}t_0u'_1;q,t;p)}
{\Delta^{0}_\mu(\frac{t^{n-1}t_0}{u_0}|t^{n-1}t_0t_1,t^{n-1}t_0t_2,t^{n-1}t_0t_3,t^{n-1}t_0u_1;q,t;p)}\notag\\
&
\quad\binomE{\lambda}{\mu}_{[\frac{1}{u_0u_1},\frac{1}{t^{n-1}t_0u_1}];q,t;p}
\binomE{\mu}{\kappa}_{[t^{n-1}t_0/u_0,t^{n-1}t_0u'_1];q,t;p}.\notag
\end{align}
\end{thm}

If $t'_3=t_3$, the same connection coefficients can be computed via
$R^{*(n)}(;t_3,u_0;q,t;p)$; the result is precisely our generalized Bailey
transformation, Theorem \ref{thm:theta_bailey}.  If also $t'_2=t_2$, the bulk
difference equation \eqref{eq:theta_bde} applies, giving the following
result announced in \cite{xforms}.

\begin{cor}
\begin{align}
[\tilde{R}^{(n)}_{\mu}(;t_0{:}t_1v,t_2,t_3;u_0,u_1/v;q,t;p)]
\tilde{R}^{(n)}_{\lambda}(;t_0{:}t_1,t_2,t_3;&u_0,u_1;q,t;p)\\
&=
\obinomE{\lambda}{\mu}_{[1/u_0u_1,1/v](t^{n-1}t_2t_3,pqt^{n-1}t_0/u_0,t_1v/u_1);q,t;p}\notag
\end{align}
\end{cor}

We also have discrete biorthogonality, which was derived in \cite{xforms}
via residue calculus, but can also be derived via the present
``hypergeometric'' methods.

\begin{thm}
For any partitions $\lambda,\kappa\subset m^n$, and
for otherwise generic parameters satisfying $t_0t_1=q^{-m}t^{1-n}$,
$t^{n-1}t_2t_3u_0u_1 = pq^{m+1}$
\begin{align}
\sum_{\mu\subset m^n}
\tilde{R}^{(n)}_{\lambda}(\dots t_0 t^{n-i} q^{\mu_i}\dots;t_0{:}t_1,t_2,t_3;u_0,u_1;q,t;p)
\tilde{R}^{(n)}_{\kappa}(\dots t_0 t^{n-i} q^{\mu_i}\dots;t_0{:}t_1,t_2,t_3;u_1,u_0;q,t;p)\\
\Delta_{\mu}(t^{2n-2} t_0^2|t^n,t^{n-1}t_0t_1,t^{n-1}t_0t_2,t^{n-1}t_0t_3,t^{n-1}t_0u_0,t^{n-1}t_0u_1;q,t;&p)
=
0\notag
\end{align}
unless $\lambda=\kappa$, when the sum is
\[
\frac{\Delta^0_{m^n}(\frac{t^{n-1}t_1}{u_0}|\frac{t_1}{t_0},\frac{pq}{u_0t_2},\frac{pq}{u_0t_3},\frac{pq}{u_0u_1};q,t;p)}
{\Delta_\lambda(1/u_0u_1|t^n,t^{n-1}t_0t_1,t^{n-1}t_0t_2,t^{n-1}t_0t_3,\frac{1}{t^{n-1}t_0u_0},\frac{1}{t^{n-1}t_0u_1};q,t;p)}
\]
\end{thm}

\begin{proof}
The argument of \cite{bcpoly} carries over essentially verbatim;
alternatively, the argument of \cite{RosengrenH:2003b} generalizes, using 
Theorem \ref{thm:multi6j}.
\end{proof}

\begin{rem}
Note that the above inner product is normalized by equation
\eqref{eq:ip_normalization}.
\end{rem}

We also have a special quasi-branching rule (having the branching rule for
interpolation functions (Theorem \ref{thm:Rs_branch}) as a special case):

\begin{thm}\label{thm:quasi-branch}
The biorthogonal functions satisfy the expansion
\begin{align}
\tilde{R}^{(n+k)}_{\lambda}(\dots z_i\dots,t_0,t_0t,\dots t_0t^{k-1};t_0{:}t_1,t_2,t_3;u_0,u_1;&q,t;p)\\
&=
\sum_{\kappa\subset\lambda}
c_{\lambda\kappa}
\tilde{R}^{(n)}_{\kappa}(\dots
z_i\dots;t_0t^k{:}t_1,t_2,t_3;u_0,u_1t^k;q,t;p)\notag
\end{align}
where
\[
c_{\lambda\kappa}
=
\obinomE{\lambda}{\kappa}_{[\frac{1}{u_0u_1},t^k](\frac{pq}{t^{n+k}u_0u_1},\frac{1}{t^{k-1}t_0u_1},\frac{pqt^{n+k-1}t_0}{u_0});q,t;p}
\]
\end{thm}

\begin{proof}
Expand the left-hand side in interpolation functions, apply equation
\eqref{eq:specbranch}, and expand back into biorthogonal functions; again
the resulting sum can be simplified via the bulk difference equation
\eqref{eq:theta_bde}.
\end{proof}

Similarly, we obtain a special quasi-Pieri identity.

\begin{thm}\label{thm:quasi-Pieri}
The biorthogonal functions satisfy the expansion
\begin{align}
\prod_{1\le i\le n}
\frac{\theta(t_0 z_i,t_0/z_i;q;p)_m}
     {\theta((pq/u_0) z_i,(pq/u_0)/z_i;q;p)_m}
\tilde{R}^{(n)}_{\lambda}(\dots z_i\dots;&q^m t_0{:}t_1,t_2,t_3;q^{-m} u_0,u_1;q,t;p)\\
&=
\sum_\kappa
c_{\kappa\lambda}
\tilde{R}^{(n)}_{\kappa}(\dots z_i\dots;t_0{:}t_1,t_2,t_3;u_0,u_1;q,t;p),\notag
\end{align}
where
\begin{align}
c_{\kappa\lambda}
&=
\Delta^0_{m^n}(\frac{t^{n-1}t_0}{u_0}|t^{n-1}t_0t_1,t^{n-1}t_0t_2,t^{n-1}t_0t_3,t^{n-1}t_0u_1;q,t;p)\notag\\
&\phantom{{}={}}
\obinomE{\kappa}{\lambda}_{[\frac{1}{u_0u_1},q^{-m}](t^n,\frac{1}{t^{n-1}t_0u_1},\frac{q^mt_0}{u_0});q,t;p}.
\end{align}
\end{thm}

Finally, we have a Cauchy identity for biorthogonal functions.

\begin{thm}\label{thm:Cauchy_biorth}
The function
\[
F(x_1,\dots x_n;y_1,\dots y_m):=
\frac{\prod_{1\le i\le n,1\le j\le m} \theta(y_j x_i,y_j/x_i;p)}
{\prod_{1\le i\le n} \theta(q^{-m}u_0 x_i,q^{-m}u_0/x_i;q;p)_m
\prod_{1\le j\le m} \theta(pq^m/u_0y_j,q^m y_j/u_0;1/t;p)_n}
\]
admits an expansion
\begin{align}
\sum_{\mu\subset m^n}
c_\mu
\tilde{R}^{(n)}_{\mu}(\dots x_i\dots;t_0{:}t_1,t_2,t_3;&u_0,u_1;q,t;p)
\tilde{R}^{(m)}_{n^m-\mu'}(\dots y_i\dots;t_0{:}t_1,t_2,t_3;
\frac{t^n u_0}{q^m},\frac{t^{n-1}u_1}{q^{m-1}};t,q;p),\notag\\
&=\frac{C^0_{m^n}(pq/u_0t_0,q^m t_0/u_0;q,t;p)F(x_1,\dots x_n;y_1,\dots y_m)}
{\Delta^0_{m^n}(t^{n-1}t_0/u_0|t^{n-1}t_0t_1,t^{n-1}t_0t_2,t^{n-1}t_0t_3,t^{n-1}t_0u_1;q,t;p)}
\end{align}
where
\[
c_\mu=
\Delta_\mu(\frac{1}{u_0u_1}|t^n,q^{-m},\frac{1}{t^{n-1}t_0u_1},\frac{q^m t_0}{u_0};q,t;p).
\]
\end{thm}

\section{Algebraic modularity and rationality}\label{sec:modular}

The purpose of the present section is to give a purely algebraic definition
of the biorthogonal functions; as a consequence, it will follow that the
biorthogonal abelian functions are in addition modular functions.  This
could of course be shown directly by determining how the interpolation
functions and $C^*$ symbols behave under modular transformation; the
approach used here has the advantage of being more conceptual in nature.
In addition, we also obtain analogues of interpolation functions, binomial
coefficients, and biorthogonal functions over arbitrary fields (including
those of positive characteristic).

In addition to making sure the construction is sensible geometrically,
we also want things to be reasonable {\em arithmetically}; that
is, in such a way that the functions depend rationally on the
parameters.  In particular, the difference equation as given above
is problematical in this respect, as it requires us to choose a
square root of $q$, despite the fact that the interpolation functions
themselves are independent of that choice.  The simplest way to
fix that is to relax the notion of elliptic curve slightly, by
forgetting which point represents the identity; this allows us to
absorb the freedom in choosing the square root.  We thus obtain
the following fundamental definition.  Recall that for a (smooth)
curve $C$, the Picard group $\Pic(C)$ is the group of divisors
modulo principal divisors; for an integer $n$, $\Pic^n(C)$ is then
the preimage of $n$ under the natural degree map.  For each $n$,
$\Pic^n(C)$ has a natural structure of algebraic variety, a principal
homogeneous space over the algebraic group $\Pic^0(C)$.  We will
use multiplicative notation for divisors (and the group law in
$\Pic(C)$), and will denote the divisor associated to a point $p\in
C$ by $\la p\ra$, and to a function $f$ on $C$ by $\la f\ra$.

\begin{defn}
A {\em genus 1 hyperelliptic curve} is a pair $(C,\tau)$, where $C$ is a
genus 1 curve, and $\tau\in \Pic^2(C)$ is a divisor class of degree 2
(which thus induces an involutory automorphism $\iota_\tau:x\mapsto \tau/x$
of $C$ with genus 0 quotient).  A $BC_n(\tau)$-symmetric function on $C$ is
a function $f$ on $C^n$ invariant under permutations of the variables and
replacements $x_i\mapsto \tau/x_i$.
\end{defn}

\begin{drems}
Note in particular that $\Pic^0(C)$ is an elliptic curve, and $\Pic^1(C)$
can be canonically identified with $C$.  Each point $x\in C$ thus induces a
map $\Pic(C)\to \Pic^0(C)$ defined by $y\mapsto y x^{-\deg(y)}$, and in
particular gives an identification of $C$ with the elliptic curve
$\Pic^0(C)$.
\end{drems}

\begin{drems}
The action of $\iota_\tau$ on $C$ extends to an automorphism of the Picard
group $\Pic(C)$, which we will denote by $\iota_\tau(x):=\tau^{\deg(x)}/x$.
Similarly, its action on the $i$th copy of $C$ in $C^n$ will be denoted by
$\iota_{i,\tau}$.
\end{drems}

\begin{drems}
For $x\in \Pic^0(C)$, there is a natural map from the genus 1 hyperelliptic
curve $(C,\tau)$ to the curve $(C,x^2 \tau)$ given by $p\mapsto xp$.
Since all of the functions we will define below are canonically defined,
they in particular will transform nicely under this isomorphism; in
general, replacing each element $y\in \Pic(\tau)$ by $x^{\deg(y)}y$ will
leave the function invariant.
\end{drems}

While there is in general no canonical choice of divisor representing
$\tau$, there is in fact a canonical choice of representative for $\tau^2$,
namely the ramification divisor of the quotient map $C\mapsto
C/\iota_\tau$.  In odd characteristic, this is simply the product of the
fixed points of $\iota_\tau$ in $\Pic^1(C)$; in even characteristic, it is
the square or fourth power of that product (depending on whether the curve
is ordinary or supersingular).  We will denote this divisor by
$\la\tau^2\ra$.

Given a divisor $D$, let ${\cal L}(D)$ denote the space of functions $f$ on
$C$ such that $\la f\ra D$ is effective; by the Riemann-Roch theorem, this
space has dimension $\max(\deg(D),0)$, unless $\deg(D)=0$ and $D$ is
principal (when the dimension is 1).  The following proposition is key to
proving the difference equation in arbitrary characteristic.

\begin{prop}\label{prop:th_symm1_alg}
Let $D$ be a divisor of degree $m$.  Then $\iota_\tau^*{\cal L}(D\iota_\tau(D))={\cal
  L}(D\iota_\tau(D))$, and the dimension of the space of $\iota_\tau^*$-invariants is
$\max(m+1,0)$.  Moreover,
\begin{align}
(1-\iota_\tau^*){\cal L}(D\iota_\tau(D))&\subset {\cal
  L}(D\iota_\tau(D)/\la\tau^2\ra),\\
(1+\iota_\tau^*){\cal L}(D\iota_\tau(D)\la\tau^2\ra)&\subset {\cal
  L}(D\iota_\tau(D)),
\end{align}
with images of dimension $\max(m-1,0)$, $\max(m+1,0)$, respectively.
\end{prop}

\begin{proof}
The statement about the space of $\iota_\tau^*$-invariants follows
immediately from Riemann-Roch, since $\iota_\tau^*$-invariants are simply
functions on the quotient $\P^1$.  For the remaining statements,
multiplication by a $\iota_\tau^*$-invariant function allows us to replace
$D$ by any other divisor of degree $m$; in particular, we can assume that
$D$ is supported away from the support of $\la \tau^2\ra$.  That
\[
(1-\iota_\tau^*){\cal L}(D\iota_\tau(D))\subset {\cal
  L}(D\iota_\tau(D)/\la\tau^2\ra)
\]
then follows from the definition of ramification divisor; the dimension of
the image can be computed from the dimension of the kernel.  The remaining
claims follow upon multiplication by a $\iota_\tau^*$-anti-invariant function.
\end{proof}

Fix a genus 1 hyperelliptic curve, as well as two points $q,t\in
\Pic^0(C)$.  Given four generic points $u_0$, $u_1$, $u_2$, $u_3\in C$ such
that $t^{n-1}u_0u_1u_2u_3=\tau^2$, we define a difference operator ${\tilde D}^{(n)}(;u_0,u_1,u_2;q,t;C,\tau)$ acting on $BC_n(\tau/q)$-symmetric
functions as follows:
\[
({\tilde D}^{(n)}(;u_0,u_1,u_2;q,t;C,\tau)f)(\dots x_i\dots)
=
\prod_{1\le i\le n} (1+\iota_{i,\tau}^*)
\delta^{(n)}(\dots x_i\dots;u_0,u_1,u_2,u_3;t;C,\tau)
f(\dots x_i\dots)
\]
for a function $\delta^{(n)}$ defined inductively as follows:
\begin{itemize}
\item[(0)] $\delta^{0}(u_0,u_1,u_2,u_3;t;C,\tau)=1$.
\item[(1)] As a function of $x_n$,
$\delta^{(n)}(\dots x_i\dots;u_0,u_1,u_2,u_3;t;C,\tau)$ has divisor
\[
\frac{
\la\frac{\tau}{u_0}\ra
\la\frac{\tau}{u_1}\ra
\la\frac{\tau}{u_2}\ra
\la\frac{\tau}{u_3}\ra
\prod_{1\le i<n} \la\frac{\tau}{t x_i}\ra}
{\la\tau^2\ra\prod_{1\le i<n} \la\frac{\tau}{x_i}\ra}
\]
\item[(2)] Setting $x_n=u_0$ gives
\[
\delta^{(n)}(x_1,\dots, x_{n-1},u_0;u_0,u_1,u_2,u_3;t)
=
\delta^{(n-1)}(x_1,\dots, x_{n-1};tu_0,u_1,u_2,u_3;t)
\]
\end{itemize}

\begin{prop}
For fixed $(C,\tau)$ and $t$, $\delta^{(n)}$ is a well-defined rational
function of all variables and parameters.  It is invariant under
permutations of $x_1$ through $x_n$ and under permutations of $u_0$ through
$u_3$.
\end{prop}

\begin{proof}
For $n=1$, we find from Proposition \ref{prop:th_symm1_alg} that
\[
(1+\iota_\tau^*)\delta^{(1)}(x;u_0,u_1,u_2,u_3;t;C,\tau)
\]
is a constant; evaluating this constant at $x=u_0$ shows that it is 1.
Evaluating at $x=u_1$ shows that
\[
\delta^{(1)}(x;u_0,u_1,u_2,u_3;t;C,\tau)
=
\delta^{(1)}(x;u_1,u_0,u_2,u_3;t;C,\tau)
\]
as required, and thus the proposition holds when $n=1$.

For $n>1$, we first observe by induction that $\delta^{(n)}$ is symmetric
in $x_1$ through $x_{n-1}$ and in $u_1$, $u_2$, $u_3$.  Now, if we set
$x_n=u_0$, the divisor in $x_{n-1}$ will be
\[
\frac{
\la\frac{\tau}{tu_0}\ra\la\frac{\tau}{u_1}\ra\la\frac{\tau}{u_2}\ra\la\frac{\tau}{u_3}\ra
\prod_{1\le i<n-1} \la\frac{\tau}{t x_i}\ra}
{\la\tau^2\ra\prod_{1\le i<n-1} \la\frac{\tau}{x_i}\ra};
\]
from the dependence on $x_{n-1}$ of the divisor in $x_n$, we deduce the
dependence on $x_n$ of the divisor in $x_{n-1}$, and thus conclude that
the divisor in $x_{n-1}$ of $\delta^{(n)}$ is
\[
\frac{
\la\frac{\tau}{u_0}\ra\la\frac{\tau}{u_1}\ra\la\frac{\tau}{u_2}\ra\la\frac{\tau}{u_3}\ra
\prod_{i\ne n-1} \la\frac{\tau}{t x_i}\ra}
{\la\tau^2\ra\prod_{1\le i<n} \la\frac{\tau}{x_i}\ra}
\]
We claim that setting $x_{n-1}=u_1$ gives
$\delta^{(n-1)}(;u_0,tu_1,u_2,u_3)=\delta^{(n-1)}(;tu_1,u_0,u_2,u_3)$;
indeed, setting $x_n=u_0$, $x_{n-1}=u_1$ gives
$\delta^{(n-2)}(;tu_0,tu_1,u_2,u_3)$, and the divisor in $x_n$ is also
correct.  It follows that
\[
\delta^{(n)}(x_1,\dots, x_n;u_0,u_1,u_2,u_3;t;C,\tau)
=
\delta^{(n)}(x_1,\dots, x_{n-2},x_n,x_{n-1};u_1,u_0,u_2,u_3;t;C,\tau),
\]
at which point the proposition follows.
\end{proof}

\begin{rem}
If $C$ is a complex elliptic curve of the form $\C^*/\la p\ra$,
then $\delta^{(n)}$ can be expressed in theta functions as
\begin{align}
\delta^{(n)}(\dots x_i\dots;u_0,u_1,u_2,u_3;&q,t;C,\tau)\\
&=
\prod_{1\le i\le n}
\frac{\theta(u_0 x_i/\tau,u_1 x_i/\tau,u_2 x_i/\tau,u_3 x_i/\tau;p)}
     {\theta(x_i^2/\tau) \theta(u_0 u_1 t^{n-i}/\tau,u_0u_2
       t^{n-i}/\tau,u_0u_3 t^{n-i}/\tau;p)}
\prod_{1\le i<j\le n}
\frac{\theta(t x_ix_j/\tau;p)}{\theta(x_ix_j/\tau;p)},\notag
\end{align}
where the parameters are lifted to $\C^*$ so that $u_0u_1u_2u_3=p\tau^2$;
thus the difference operator defined above agrees (up to a scale factor)
with the earlier definition.
\end{rem}

For an integer $m\ge 0$, let $A_m(u_0;q;C,\tau)\subset {\cal L}(\prod_{1\le
  i\le n} \la q^{-i} u_0\ra\la \frac{\tau}{q^{-i}u_0}\ra)$ denote the
subspace of $\iota_\tau$-invariant functions, and let
$A^{(n)}_m(u_0;q;C,\tau)$ be the $n$-th symmetric power of that space,
viewed as a space of functions on $C^n$; in particular,
$A^{(n)}_m(u_0;q;C,\tau)$ consists of $BC_n(\tau)$-symmetric functions.  We
also write $A^{(n)}(u_0;q;C,\tau)$ for the union $\lim_{m\to\infty}
A^{(n)}_m(u_0;q;C,\tau)$.

\begin{thm}\label{thm:modular_diff_op}
For any integer $m\ge 0$,
\[
{\tilde D}^{(n)}(u_0,u_1,u_2;q,t;C,\tau)
A^{(n)}_m(u_0;q;C,\tau/q)
\subset
A^{(n)}_m(u_0;q;C,\tau).
\]
Moreover,
\begin{align}
({\tilde D}^{(n+k)}(u_0,u_1,u_2;q,t;C,\tau)f)&(z_1,\dots,z_n,u_0,\dots, t^{k-1}u_0)\notag\\
&=
({\tilde D}^{(n)}(u_0,u_1,u_2;q,t;C,\tau)f(\underline{\quad},u_0,\dots,t^{k-1}u_0))
(z_1,\dots,z_n).
\end{align}
\end{thm}

\begin{proof}
The second claim follows readily from the definition of the difference
operator and of $\delta^{(n)}$, so we need only consider the first claim.
In particular, it suffices to show that as a function of $z_n$, the only
poles of
\[
{\tilde D}^{(n)}(u_0,u_1,u_2;q,t;C,\tau)f
\]
for $f\in A^{(n)}_m(u_0;q;C,\tau/q)$ are of the form $q^{-i} u_0$ or $q^i
\tau/u_0$ for $1\le i\le m$, of multiplicity at most 1.  Now, as a
function of $z_n$,
\[
\delta^{(n)}(\dots z_i\dots;u_0,u_1,u_2;q,t;C,\tau)f(\dots z_i\dots)
\]
has polar divisor at most
\[
\la \tau^2\ra
\prod_{1\le i<n} \la \tau/z_i\ra
\prod_{1\le i\le m} \la q^{-i} u_0\ra\la q^{i-1}\tau/u_0\ra
/\la \tau/u_0\ra,
\]
with the last factor coming from the numerator of $\delta^{(n)}$.
In particular, upon symmetrization by $\iota_\tau$, the result clearly has
polar divisor at most
\[
\la \tau^2\ra
\prod_{1\le i<n} \la z_i\ra\la \tau/z_i\ra
\prod_{1\le i\le m} \la q^{-i} u_0\ra\la q^i \tau/u_0\ra.
\]
We thus need only show that the potential poles corresponding to the first
two sets of factors disappear upon symmetrization; we may assume (since it
is true generically) that the three sets of factors have disjoint support.

That the singularities corresponding to $\la\tau^2\ra$ disappear follows
immediately from Proposition \ref{prop:th_symm1_alg}, so it remains by
symmetry to consider the potential singularity at $z_n=z_{n-1}$.  Now, in
the sum
\begin{align}
\prod_{1\le i\le n} (1+\iota_{i,\tau}^*)&
\delta^{(n)}(\dots z_i\dots;u_0,u_1,u_2;q,t;C,\tau)f(\dots z_i\dots)\\
=&
\prod_{1\le i\le n-2} (1+\iota_{i,\tau}^*)
[
(1+\iota_{n-1,\tau}^*+\iota_{n,\tau}^*+\iota_{n-1,\tau}^*\iota_{n,\tau}^*)
\delta^{(n)}(\dots z_i\dots;u_0,u_1,u_2;q,t;C,\tau)f(\dots z_i\dots)
]\notag
\end{align}
the factors corresponding to $1$ and $\iota_{n-1,\tau}^*\iota_{n,\tau}^*$
are already nonsingular at $z_n=z_{n-1}$, so it suffices to show that
\[
(\iota_{n-1,\tau}^*+\iota_{n,\tau}^*)
\delta^{(n)}(\dots z_i\dots;u_0,u_1,u_2;q,t;C,\tau)f(\dots z_i\dots)
\]
has no pole at $z_{n-1}=z_n$.  Now, since $\delta^{(n)}$ and $f$ are
symmetric under permutations of $z_1,\dots, z_n$, this sum can be written as
the symmetrization of
\[
\iota_{n,\tau}^*
\delta^{(n)}(\dots z_i\dots;u_0,u_1,u_2;q,t;C,\tau)f(\dots z_i\dots)
\]
under exchanging $z_{n-1}$ and $z_n$; the disappearance of the pole then
follows from local considerations.
\end{proof}

We will also need to construct algebraic analogues of certain products of
$C^*$-symbols.  First, we need the following family of univariate
functions.  For a nonnegative integer $m$, let $u_0,\dots,u_{2m-1}\in C$
be a collection of $2m$ points such that $u_0u_1\cdots u_{2m-1}=\tau^m$, and
let $D$ be a divisor of degree $m$.  Then there exists up to scale a unique
function $f$ with divisor
\[
\frac{\prod_{0\le i<2m}\la u_i\ra}
     {(D\iota_\tau(D))}.
\]
Thus the function
\[
\omega(x{:}u_0,\dots,u_{2m-1};C,\tau)
:=
\frac{f(x)}{f(\tau/x)}
\]
is well-defined; multiplication of $f$ by an $\iota^*_\tau$-invariant
function shows that $\omega$ is independent of the choice of $D$.
Note that 

\begin{lem}
The function $\omega$ has divisor
\[
\frac{\prod_{0\le i<2m}\la u_i\ra}
     {\prod_{0\le i<2m}\la\tau/u_i\ra},
\]
and satisfies the normalization that for any point $x_0\in C$ such that $x_0^2=\tau$,
\[
\omega(x_0{:}u_0,\dots,u_{2m-1};C,\tau)=1.
\]
Moreover,
\[
\omega(x{:}u_0,\dots,u_{2m-1},v,\tau/v;C,\tau)=
\omega(x{:}u_0,\dots,u_{2m-1};C,\tau),
\]
and
\[
\omega(x{:}u_0,\dots,u_{2m-1};C,\tau)
\omega(x{:}v_0,\dots,v_{2n-1};C,\tau)
=
\omega(x{:}u_0,\dots,u_{2m-1},v_0,\dots,v_{2n-1};C,\tau).
\]
\end{lem}

\begin{proof}
The first claim is straightforward; for the second claim, we find that
\[
\omega(x{:}u_0,\dots,u_{2m-1};C,\tau)-1
=
\frac{f(x)-f(\tau/x)}{f(\tau/x)}
\in
{\cal L}(\la \tau^2\ra^{-1}),
\]
and thus it vanishes as required.  The final claims follow immediately
from the original definition of $\omega$.
\end{proof}

Another univariate function can be obtained from the observation that
since $\tau$ induces a map to $\P^1$, any four points induce a uniquely
defined scalar, namely the cross-ratio of their images.  We denote this by
\[
\chi(a,b,c,d;C,\tau),
\]
which can also be defined as the unique $\iota_\tau^*$-invariant function
of $a$ such that
\[
\chi(b,b,c,d;C,\tau)=0,\quad
\chi(c,b,c,d;C,\tau)=1,\quad
\chi(d,b,c,d;C,\tau)=\infty.
\]
We can also define this by
\[
\chi(a,b,c,d;C,\tau)
=
\omega(x{:}bax/\tau,bx/a,dx/c,dcx/\tau;C,x^2bd/\tau),
\]
for any $x\in C$.  Note also that since $\chi$ is defined as a cross-ratio,
all of the usual transformations of cross ratios apply; note in particular
the symmetry
\[
\chi(a,b,c,d;C,\tau)
\chi(a',b,c',d;C,\tau)
=
\chi(a,b,c',d;C,\tau)
\chi(a',b,c,d;C,\tau).
\]
In addition, the identification with $\omega$, with its different natural
symmetry, gives rise to the transformations
\begin{align}
\chi(a,b,c,d;C,\tau)&=\chi(a,b,c,d;C,abcd/\tau)\\
\omega(x{:}v_0,v_1,v_2,v_3;C,\tau)
&=
\omega(x{:}v_0,v_1,x^2/v_2,x^2/v_3;C,\tau x^2/v_2v_3).
\end{align}

\begin{defn}
Let $a$, $v_0,\dots,v_{2m-1}\in \Pic^0(C)$ be such that
\[
v_0v_1\cdots v_{2m-1} = q^m a^m.
\]
Then for any partition $\lambda$, we define
\[
\Delta^0_\lambda(a|v_0,\dots,v_{2m-1};q,t;C)
=
\prod_{(i,j)\in \lambda}
\omega(q^{1-j}t^{i-1}x{:}v_0 x,\dots,v_{2m-1} x;C,qax^2),
\]
for a generic point $x\in C$ (on which the value does not depend).
\end{defn}

\begin{defn}
Let $a$, $b_0,\dots,b_{2l-1}\in \Pic^0(C)$ be such that $\prod_{0\le
  r<2l} b_r = (t/q) (qa)^{l-1}$.  Then we define a symbol
$\Delta_\lambda(a|b_0,\dots,b_{2l-1};q,t;C)$ for all partitions $\lambda$
via the recurrence
\begin{align}
\Delta_{m\cdot\lambda}(a|b_0,\dots,b_{2l-1};q,t;C)
&=
\lim_{x\to 1}
\Delta^0_{m}(a|b_0,\dots, b_{2l-1},q^{m+1}a/x,qa/xt,qax,qx/q^m;q,t;C)\\
&\phantom{{}={}}\Delta_{\lambda}(t^{-2}a|b_0/t,\dots,b_{2l-1}/t,q^{1-m}/t,q^{-m},q^{m+1}a/t,q^ma;q,t;C)\notag
\end{align}
with
\[
\Delta_0(a|b_0,\dots,b_{2l-1};q,t;C)=1;
\]
here, of course, the limit is to be taken in the algebraic sense: the
evaluation at the appropriate point of the stated rational function in the
limit variable.  Similarly, we define
$\binomE{\lambda}{\lambda}_{[a,b];q,t;C}$ via the recurrence
\begin{align}
\binomE{m\cdot\lambda}{m\cdot\lambda}_{[a,b];q,t;C}
&=
\lim_{x\to 1} \Delta^0_{m}(a/b|q^m a/x,1/bx,qax/b,qx/q^m;q,t;C)\\
&\phantom{{}={}}
\frac{\Delta^0_\lambda(a/t^2|b,qa/tb,q^{1-m}/t^2,q^ma/t;q,t;C)}
     {\Delta^0_\lambda(a/bt^2|1/b,qa/t,q^{1-m}/t^2,q^ma/tb;q,t;C)}
\binomE{\lambda}{\lambda}_{[a/t^2,b];q,t;C},\notag
\end{align}
with
\[
\binomE{0}{0}_{[a,b];q,t;C}=1.
\]
\end{defn}

\begin{rem}
For a complex elliptic curve, the symbols $\Delta^0_\lambda$ and
$\Delta_\lambda$ are precisely those defined above; the parameters must be
lifted so that the constraints on $\prod_i b_i$ are satisfied as given,
except with $q$ replaced by $pq$.  The remaining symbol corresponds to the
product
\[
\frac{C^+_\lambda(a  ;q,t;p)C^0_\lambda(1/b,pqa/b;q,t;p)}
     {C^+_\lambda(a/b;q,t;p)C^0_\lambda( b ,pqa  ;q,t;p)}
=
\binomE{\lambda}{\lambda}_{[a,b];q,t;p}.
\]
\end{rem}

The algebraic analogue of interpolation functions are defined as follows.

\begin{defn}
For generic parameters $a,b,v\in C$, $q,t\in \Pic^0(C)$, we define the
algebraic interpolation function $R^{*(n)}_\lambda(;a,b(v);q,t;C,\tau)$ as
follows.
\begin{itemize}
\item[(1)] $R^{*(n)}_\lambda(;a,b(v);q,t;C,\tau)\in
  A^{(n)}_{\lambda_1}(b;q;C,\tau)$.
\item[(2)]
For each integer $m>\lambda_1$, and every partition $\mu\subset m^n$
with $\mu\ne \lambda$, let $l$ be as in the definition of $P^{*(m,n)}$
above.  Then for generic $c\in C$,
\[
\!\!\!\!\!\!\!\!\!
\Bigl(
\prod_{\substack{1\le i\le n\\ 1\le j\le m}}
\chi(z_i,q^{-j} b,t^{n-i} v,c;C,\tau)
R^{*(n)}_\lambda(;a,b(v);q,t;C,\tau)
\Bigr)(b q^{-\mu_1},\dots,b q^{-\mu_l} t^{l-1},
  a q^{\mu_{l+1}} t^{n-l-1},\dots,a q^{\mu_n})
=
0.
\]
\item[(3)]
$R^{*(n)}_\lambda$ satisfies the normalization
\[
R^{*(n)}_\lambda(\dots v t^{n-i}\dots;a,b(v);q,t;C,\tau)=1.
\]
\end{itemize}
\end{defn}

As before, these functions are generically uniquely determined by the
conditions; moreover, they satisfy a difference equation, the proof of
which is a direct analogue of the corresponding proof for theta functions.

\begin{thm}
The functions $R^{*(n)}_\lambda$ satisfy the difference equation
\begin{align}
{\tilde D}^{(n)}(a,b,c;q,t;C,\tau)
R^{*(n)}_{\lambda}(;a,b(v);&q,t;C,\tau/q)\\
&=
\Delta^0_\lambda(\frac{t^{n-1}a}{b}|\frac{avt^{n-1}}{\tau},\frac{acqt^{n-1}}{\tau},\frac{\tau}{bv},\frac{q\tau}{bc};q,t;C)
R^{*(n)}_{\lambda}(;a,b(v);q,t;C,\tau)\notag
\end{align}
and the normalization condition
\[
R^{*(n)}_\lambda(\dots v' t^{n-i}\dots;a,b(v);q,t;C,\tau)
=
\Delta^0_\lambda(t^{n-1}a/b|t^{n-1}av'/\tau,a/v',q\tau/vb,t^{n-1}qv/b;q,t;C).
\]
\end{thm}

\begin{proof}
As mentioned, the proof follows that of Theorem \ref{thm:theta_diff}; the
only difference is the treatment of the normalization.  We find as before that
\[
{\tilde D}^{(n)}(a,b,c;q,t;C,\tau)
R^{*(n)}_{\lambda}(;a,b(v);q,t;C,\tau/q)
\propto
R^{*(n)}_{\lambda}(;a,b(v);q,t;C,\tau);
\]
for $c=v$, we find immediately that
\[
{\tilde D}^{(n)}(a,b,v;q,t;C,\tau)
R^{*(n)}_{\lambda}(;a,b(v);q,t;C,\tau/q)
=
R^{*(n)}_{\lambda}(;a,b(v);q,t;C,\tau).
\]
Now, if we evaluate at the point $(\dots a t^{n-i} q^{\lambda_i}\dots)$, 
we find that only one term survives, and thus the constant of
proportionality in the difference equation can be computed as
\[
\frac{
\delta^{(n)}(\dots a t^{n-i} q^{\lambda_i}\dots;a,b,c,\tau^2/abc;t;C,\tau)}
{
\delta^{(n)}(\dots a t^{n-i} q^{\lambda_i}\dots;a,b,v,\tau^2/abv;t;C,\tau)}
.
\]
We readily verify from divisor conditions that the function
\[
\frac{\delta^{(n)}(\dots x_i\dots;a,b,c,t^{1-n}\tau^2/abc;t;C,\tau)}
     {\delta^{(n)}(\dots x_i\dots;a,b,v,t^{1-n}\tau^2/abv;t;C,\tau)}
=
\prod_{1\le i\le n}
\chi(x_i,\tau/c,a t^{n-i},\tau/v;C,t^{n-1}ab),
\]
and we thus obtain
\[
\frac{
\delta^{(n)}(\dots a t^{n-i} q^{\lambda_i}\dots;a,b,c,\tau^2/abc;t;C,\tau)}
{
\delta^{(n)}(\dots a t^{n-i} q^{\lambda_i}\dots;a,b,v,\tau^2/abv;t;C,\tau)}
=
\prod_{1\le i\le n}
\chi(a t^{n-i} q^{\lambda_i},\tau/c,a t^{n-i},\tau/v;C,t^{n-1}ab)
\]
On the other hand, we find
\begin{align}
\Delta^0_\lambda(t^{n-1}a/b|avt^{n-1}/\tau,acqt^{n-1}/\tau,\tau/bv,q\tau/bc;q,t;C)
&=
\prod_{1\le i\le n}
\prod_{1\le j\le \lambda_i}
\chi(a t^{n-i}q^j,\tau/c,a t^{n-i} q^{j-1},\tau/v;C,t^{n-1}ab)\\
&=
\prod_{1\le i\le n}
\chi(a t^{n-i} q^{\lambda_i},\tau/c,a t^{n-i},\tau/v;C,t^{n-1}ab)
\end{align}
as required.

The formula for changing the normalization follows by applying the
difference equation to both sides of the equation
\[
R^{*(n)}_\lambda(;a,b(v);q,t;C,\tau)
=
R^{*(n)}_\lambda(;a,b(t^{1-n}\tau/v);q,t;C,\tau),
\]
which gives a recurrence for the relevant scale factor, with the stated
solution.
\end{proof}

In particular, we obtain the extra vanishing conditions just as before.
The various identities of Section \ref{sec:int_theta} for the analytic
interpolation functions all carry over to the algebraic interpolation
functions.  In particular, we note the Cauchy case
\[
R^{*(n)}_\lambda(\dots z_i\dots;a,b(v);q,t;C,t^n ab/q)
=
\prod_{\substack{1\le i\le n\\1\le j}}
\chi(z_i,a t^{n-\lambda'_j} q^{j-1},v t^{n-i},a t^n q^{j-1};C,t^n ab/q)
\]

\begin{prop}
We have the identities
\begin{align}
R^{*(n)}_{m^n-\lambda}(\dots z_i\dots;a,b(v);q,t;C,\tau)
&=
\frac{R^{*(n)}_\lambda(\dots z_i\dots;q^{-m}b,q^m a(v);q,t;C,\tau)}
{\prod_{\substack{1\le i\le n\\1\le j\le m}}
\chi(z_i,q^{-j} b,v t^{n-i},q^{j-1} a;C,\tau)}
,
\\
R^{*(n)}_{m^n+\lambda}(\dots z_i\dots;a,b(v);q,t;C,\tau)
&=
\frac{R^{*(n)}_\lambda(\dots z_i\dots;aq^m,b/q^m(v);q,t;C,\tau)}
{\prod_{\substack{1\le i\le n\\1\le j\le m}}
\chi(z_i,q^{-j} b,v t^{n-i},q^{j-1} a;C,\tau)}
,\\
R^{*(n+k)}_\lambda(\dots z_i\dots,a,at,\dots,at^{k-1};a,b(v);q,t;C,\tau)
&=
\Delta^0_\lambda(t^{n+k-1}a/b|t^n,qa/bt,t^k a/v,qt^{n+k-1}v/b;q,t;C)\notag\\
&\phantom{{}={}}
R^{*(n)}_\lambda(\dots z_i\dots;t^k a,b(v);q,t;C,\tau),
\end{align}
and thus
\[
R^{*(n)}_\lambda(\dots a t^{n-i} q^{\lambda_i}\dots;a,b(v);q,t;C,\tau)
=
\frac{\binomE{\lambda}{\lambda}_{[t^{2n-2} a^2/\tau,t^{n-1}ab/\tau];q,t;C}}
{\Delta_\lambda(\frac{t^{n-1}a}{b}|t^n,\tau/t^{n-1}ab,t^{n-1}av/\tau,a/v;q,t;C)}.
\]
\end{prop}

\begin{proof}
The first two identities are straightforward, as both sides satisfy the same
vanishing conditions and normalization.  For the second identity, the
vanishing conditions are trivial to verify, and thus both sides are
proportional.  From the difference equation, it follows that the scale
factor is independent of $\tau$, and we may thus reduce to the Cauchy case,
for which the verification is straightforward.  The final equation follows
immediately.
\end{proof}

We also have the following symmetries which are trivial consequences of
the isomorphism invariance of our definitions:
\begin{align}
R^{*(n)}_\lambda(\dots z_i\dots;a,b(v);q,t;C,\tau)
&=
R^{*(n)}_\lambda(\dots (x/a)z_i\dots;x,x b/a(xv/a);q,t;C,x^2 \tau/a^2)\\
&=
R^{*(n)}_\lambda(\dots z_i\dots;\tau/a,\tau/b(\tau/v);1/q,1/t;C,\tau)\\
&=
R^{*(n)}_\lambda(\dots z_i\dots;a,b(\tau/t^{n-1}v);q,t;C,\tau)
\end{align}

With the above in mind, we define the algebraic binomial coefficients as
follows.  For a genus 1 curve $C$ and $a,b,q,t\in \Pic^0(C)$, we define
\[
\binomE{\lambda}{\mu}_{[a,b];q,t;C}
:=
\Delta_\mu(\frac{a}{b}|t^n,1/b,v,t^{1-n}a/v;q,t;C)
R^{*(n)}_\mu(x q^{\lambda_i} t^{1-i};t^{1-n}x,xb/a(xv/a);q,t;C,x^2/a),
\]
where the right-hand side is independent of $n$, $x\in C$, $v\in
\Pic^0(C)$.  (In particular, the binomial coefficients really are elliptic
in nature.)  When $\mu=\lambda$, this is of course consistent with our
previous notation $\binomE{\lambda}{\lambda}_{[a,b];q,t;C}$.  If $C$ is a
complex elliptic curve, then we obtain the analytic binomial
coefficients:
\[
\binomE{\lambda}{\mu}_{[a,b];q,t;\C/\la p\ra}
=
\binomE{\lambda}{\mu}_{[a,b];q,t;p}.
\]
Similarly, if $v_0v_1v_2=q^2a^2/b$, we define
\[
\obinomE{\lambda}{\mu}_{[a,b](v_0,v_1,v_2);q,t;C}
:=
\frac{
\Delta^0_\lambda(a|b,v_0,v_1,v_2;q,t;C)}
{
\Delta^0_\mu(a/b|1/b,v_0,v_1,v_2;q,t;C)}
\binomE{\lambda}{\mu}_{[a,b];q,t;C}.
\]

In particular, the same proof as in the analytic case gives the
algebraic bulk difference equation.

\begin{thm}
For otherwise generic parameters on $\Pic^0(C)$ satisfying $bcde=aq$,
\[
\binomE{\lambda}{\kappa}_{[a,c];q,t;C}
=
\frac{\Delta^0_\kappa(a/c|1/c,bd,be,aq/b;q,t;C)}
     {\Delta^0_\lambda(a|c,bd,be,aq/b;q,t;C)}
\sum_{\kappa\subset\mu\subset\lambda}
\Delta^0_\mu(a/b|c/b,qa,d,e;q,t;C)
\binomE{\lambda}{\mu}_{[a,b];q,t;C}
\binomE{\mu}{\kappa}_{[a/b,c/b];q,t;C}
\]
In particular,
\[
\sum_\mu
\binomE{\lambda}{\mu}_{[a,b];q,t;C}
\binomE{\mu}{\kappa}_{[a/b,1/b];q,t;C}
=
\delta_{\lambda\kappa}.
\]
\end{thm}

\begin{rem}
One immediate consequence is Spiridonov's observation
\cite{SpiridonovVP:2002} that the various more traditionally
hypergeometric special cases discussed following Theorem
\ref{thm:theta_bde} above are modular and abelian in all parameters.  (This
includes the $BC_n$-type sum of \cite{RosengrenH:2004}, since as observed
there, that sum can be obtained by specializing Warnaar's Schlosser-type
identity.)
\end{rem}

\begin{thm}
The algebraic interpolation functions satisfy the connection coefficient
identity
\[
[R^{*(n)}_\mu(;a',b(v);q,t;C,\tau)]
R^{*(n)}_{\lambda}(;a,b(v);q,t;C,\tau)
=
\obinomE{\lambda}{\mu}_{[\frac{t^{n-1}a}{b},\frac{a}{a'}](\frac{t^{n-1}aa'}{\tau},\frac{q\tau}{bv},\frac{qt^{n-1}v}{b});q,t;C}.
\]
\end{thm}

We define algebraic biorthogonal functions by:
\[
\tilde{R}^{(n)}_{\lambda}(;t_0{:}t_1,t_2,t_3;u_0,u_1;q,t;C,\tau)
:=
\sum_\mu
\frac{
\binomE{\lambda}{\mu}_{[\frac{\tau}{u_0u_1},\frac{\tau}{t^{n-1}t_0u_1}];q,t;C}
R^{*(n)}_\mu(;t_0,u_0(t_1);q,t;C,\tau)}
{
\Delta^0_\mu(\frac{t^{n-1}t_0}{u_0}|\frac{t^{n-1}t_0t_2}{\tau},\frac{t^{n-1}t_0t_3}{\tau},\frac{t^{n-1}t_0u_1}{\tau},
\frac{q t^{n-1}t_1}{u_0};q,t;C)
},
\]
where $t_0$, $t_1$, $t_2$, $t_3$, $u_0$, $u_1\in C$ such that
$t^{2n-2}t_0t_1t_2t_3u_0u_1 = q\tau^3$.  In particular, when $C$ is
complex, we exactly recover the analytic biorthogonal functions, which
must therefore be modular and abelian in all parameters.

\begin{thm}
The algebraic biorthogonal functions satisfy the symmetry identity
\[
\tilde{R}^{(n)}_{\lambda}(;t_1{:}t_0,t_2,t_3;u_0,u_1;q,t;C,\tau)
=
\frac{\tilde{R}^{(n)}_{\lambda}(;t_0{:}t_1,t_2,t_3;u_0,u_1;q,t;C,\tau)}
     {\tilde{R}^{(n)}_{\lambda}(\dots t^{n-i}
       t_1\dots;t_0{:}t_1,t_2,t_3;u_0,u_1;q,t;C,\tau)},
\]
where
\[
\tilde{R}^{(n)}_{\lambda}(\dots t^{n-i} t_1\dots;t_0{:}t_1,t_2,t_3;u_0,u_1;q,t;C,\tau)
=
\Delta^0_\lambda(\frac{\tau}{u_0u_1}|\frac{t^{n-1}t_1t_2}{\tau},\frac{t^{n-1}t_1t_3}{\tau},\frac{\tau}{t^{n-1}t_1u_1},\frac{qt^{n-1}t_0}{u_0};q,t;C),
\]
and the difference equation
\[
{\tilde D}^{(n)}(u_0,t_0,t_1;q,t;C,\tau)
\tilde{R}^{(n)}_{\lambda}(;t_0{:}t_1,t_2/q,t_3/q;u_0,u_1/q;q,t;C,\tau/q)
=
\tilde{R}^{(n)}_{\lambda}(;t_0{:}t_1,t_2,t_3;u_0,u_1;q,t;C,\tau)
\]
\end{thm}

From isomorphism invariance, we obtain the identities
\begin{align}
\tilde{R}^{(n)}_{\lambda}(\dots
z_i\dots;t_0{:}t_1,t_2,t_3;u_0,u_1;q,t;C,\tau)
&=
\tilde{R}^{(n)}_{\lambda}(\dots x z_i\dots;x t_0{:}x t_1,x t_2,x t_3;x
u_0,x u_1;q,t;C,x^2\tau)\\
&=
\tilde{R}^{(n)}_{\lambda}(\dots
z_i\dots;\tau/t_0{:}\tau/t_1,\tau/t_2,\tau/t_3;\tau/u_0,\tau/u_1;q,t;C,\tau).
\end{align}

The statement of evaluation symmetry for the algebraic biorthogonal
functions requires a certain amount of care, since the usual statement
involves a square root; again, this can be absorbed into the hyperelliptic
structure, but there does not seem to be as natural a way of doing so.
What we find is that the hatted parameters (including a $\hat\tau$) must be
related to the original parameters by
\begin{align}
\frac{\hat{t}_1}{\hat{t}_0} = \frac{q\tau}{t_2t_3},\quad
\frac{\hat{t}_2}{\hat{t}_0} = \frac{q\tau}{t_1t_3},\quad
\frac{\hat{t}_3}{\hat{t}_0} = \frac{q\tau}{t_1t_2},\quad
\frac{\hat{u}_0}{\hat{t}_0} = \frac{u_0}{t_0},\quad
\frac{\hat{u}_1}{\hat{t}_0} = \frac{u_1}{t_0},\quad
\frac{\hat{t}_0^2}{\hat{\tau}} = \frac{t_0t_1t_2t_3}{q\tau^2}.
\end{align}
In particular, we can choose $\hat{t}_0$ arbitrarily, at which point the
remaining hatted parameters are determined.  This freedom of $\hat{t}_0$
comes from the translation symmetries of $C$, and corresponds to the fact
that the evaluation of a biorthogonal function at a partition is really a
function on $\Pic^0(C)^7$.

The remaining identities satisfied by the analytic interpolation and
biorthogonal functions all carry over to the algebraic case
straightforwardly; for instance, we have the connection coefficient
\begin{align}
[\tilde{R}^{(n)}_{\mu}(;t_0{:}t_1v,t_2,t_3;u_0,u_1/v;q,t;C,\tau)]&
\tilde{R}^{(n)}_{\lambda}(;t_0{:}t_1,t_2,t_3;u_0,u_1;q,t;C,\tau)\\
&=
\obinomE{\lambda}{\mu}_{[\frac{\tau}{u_0u_1},\frac{1}{v}](\frac{t^{n-1}t_2t_3}{\tau},\frac{qt^{n-1}t_0}{u_0},\frac{vt_1}{u_1});q,t;C},\notag
\end{align}
and so forth.

\section{Elliptic bigrids and degenerations}\label{sec:ell_bigrid}

As we mentioned above, there is a fifth class of perfect bigrids, namely
elliptic bigrids.  In contrast to the earlier classes, the space of
elliptic bigrids is of bounded dimension (8, that is) as $m,n\to\infty$,
but as we have seen, in return for the decreased freedom, we gain a large
number of new properties.

Elliptic bigrids of a given shape are parametrized by sextuples
$(a,b;q,t;C,\phi)$, where $C$ is a genus 1 curve over a field $k$,
$\phi:C\to \P^1(k)$ is a degree 2 function (with associated divisor class
$\tau\in \Pic^2(C)$), $q,t\in \Pic^0(C)$, and $a,b\in C$.  The
corresponding elliptic bigrid of shape $m^n$ is then defined by the formula
\begin{align}
\gamma(0,i,j) &= \phi(a q^jt^{n-i})\\
\gamma(1,i,j) &= \phi(b q^{-j}t^{i-1}).
\end{align}

\begin{thm}\label{thm:elliptic_are_perfect}
Elliptic bigrids are perfect.
\end{thm}

\begin{proof}
Fix a pair $(C,\phi)$, and let the other parameters be generic; also
choose a generic point $v\in C$.  For each partition $\lambda\subset m^n$,
define a $BC_n(\tau)$-symmetric function $f_\lambda$ by:
\[
f_\lambda(\dots x_i\dots)
=
\prod_{\substack{1\le i\le n\\1\le j\le m}}
\frac{\phi(x_i)-\phi(q^j\tau/b)}{\phi(v t^{n-i})-\phi(q^j\tau/b)}
R^{*(n)}_\lambda(\dots x_i\dots;a,b(v);q,t;C,\tau).
\]
Since $f_\lambda$ is $BC_n(\tau)$-symmetric, it factors through a
symmetric rational function on $\P^1(k)$; consideration of poles shows
that this rational function is in fact a polynomial.  Now, from the
vanishing properties of $R^{*(n)}_\lambda$, it follows that
\[
f_\lambda(\dots a t^{n-i} q^{\mu_i}\dots)=0
\]
for $\mu\not\subset\lambda$.  In particular, for $\nu\subset m^n$ with
$\nu\not\subset\lambda$, we can replace the parts equal to $m$ by any
partition with parts $\ge m$.  Thus
\[
f_\lambda(\dots \gamma^+(0,i,\nu_i)\dots)=0
\]
for a Zariski-dense collection of extensions of $\gamma$.  Complementation
symmetry of $R^{*(n)}$ gives the complementary vanishing conditions.  In
other words, the polynomial associated to $f_\lambda$ is an interpolation
polynomial, and thus the generic elliptic bigrid is perfect.  Since
perfection is closed, the theorem follows.
\end{proof}

If two sextuples $(a,b;q,t;C,\phi)$, $(a',b';q',t';C',\phi')$ are
equivalent in the sense that there exists an isomorphism $\psi:C\to C'$
such that
\[
\psi(a)=a',\quad \psi(b)=b',\quad \psi(q)=q',\quad \psi(t)=t',\quad
\phi'\circ\psi=\phi,
\]
the associated elliptic bigrid is the same.  In other words, we obtain a
map from ${\cal M}$ to the space of perfect bigrids, where ${\cal M}$ is
the moduli problem classifying sextuples modulo equivalence.

\begin{thm}\label{thm:birational_bigrid}
If $m,n\ge 2$, $(m,n)\ne (2,2)$, the map from ${\cal M}$ to the space of
bigrids is birational; that is, an elliptic bigrid generically determines a
unique (up to equivalence) sextuple $(a,b;q,t;C,\tau)$.  Moreover, the
closure of the space of elliptic bigrids of shape $m^n$ is a rational
variety of dimension 8.
\end{thm}

\begin{proof}
Let $\gamma$ be a generic elliptic bigrid; we need to prove that it
determines a unique sextuple.  It suffices to consider the cases $(m,n)\in
\{(3,2),(2,3)\}$, since a larger elliptic bigrid necessarily contains one
of the two smaller bigrids.  Consider the first case.  We claim that the
sextuple, and thus the full bigrid, is uniquely determined by the eight
points
\[
\gamma(0,1,0),\gamma(0,1,1),\gamma(0,1,2),
\gamma(0,2,0),\gamma(0,2,1),\gamma(0,2,2),
\gamma(0,1,1),\gamma(0,1,2).
\]
To see this, consider the five pairs of points
\[
(\gamma(0,1,0),\gamma(0,1,1)),
(\gamma(0,1,1),\gamma(0,1,2)),
(\gamma(0,2,0),\gamma(0,2,1)),
(\gamma(0,2,1),\gamma(0,2,2)),
(\gamma(1,1,1),\gamma(1,1,2)).
\]
Each of these is of the form $(\phi(x),\phi(qx))$ for $x\in C$.

\begin{lem}
Let $C$ be a genus 1 curve, with $\phi:C\to\P^1$ a degree 2 function, and
$q\in \Pic^0$ a generic point.  Then there exists a unique (up to scale)
homogeneous polynomial $p$ on $\P^1\times \P^1$ of bidegree $(2,2)$ that
vanishes on precisely those points of the form $(\phi(x),\phi(qx))$; the
resulting subvariety of $\P^1\times \P^1$ is isomorphic to $C$.  Moreover,
this polynomial is symmetric under exchanging the factors $\P^1$.
\end{lem}

\begin{proof}
Let $\phi$ have polar divisor $\la z_1\ra\la z_2\ra$.  Then evaluation at
$(\phi(x),\phi(qx))$ maps the space of polynomials of bidegree $(2,2)$, of
dimension 9, to the space ${\cal L}((\la z_1\ra\la z_2\ra\la z_1/q\ra\la
z_2/q\ra)^2)$, of dimension 8; existence of $p$ follows.  Any such
polynomial cuts out a one-dimensional subscheme of arithmetic genus
$(2-1)(2-1)=1$, which is therefore the full collection of points
$(\phi(x),\phi(qx))$, a curve isomorphic to $C$.  In particular, uniqueness
of $p$ follows from the connectedness of the subscheme.  Since
\[
(\phi(x),\phi(qx))=(\phi(\tau/x),\phi(\tau/qx)),
\]
the subscheme is symmetric, and thus the symmetry of $p$ follows.
\end{proof}

\begin{rem}
Such ``symmetric biquadratic relations'' appear in the work of Baxter on
the eight-vertex model, and other related statistical mechanical models
\cite{BaxterRJ:1982}; they have also been mentioned in the context of
univariate elliptic biorthogonal functions
\cite{SpiridonovVP/ZhedanovAS:2000b,SpiridonovVP/ZhedanovAS:2001}.  In
particular, the latter reference derives such a relation from the existence
of an appropriate difference operator.
\end{rem}

Now, the space of symmetric polynomials of bidegree $(2,2)$ is
six-dimensional, and thus five points in $\P^1\times \P^1$ generically
determine a unique (up to scale) such polynomial; a simple computation of
minors shows that this remains true for our somewhat constrained points.
We thus obtain a genus 1 curve $C_0$, together with a degree 2 map $\phi_0$
given by projection onto the first factor.  The given five points are then
of the form $a_0t_0$, $a_0q_0t_0$, $a_0$, $a_0q_0$, $b_0/q_0^2$
for uniquely determined points $a_0$, $b_0\in C_0$, $q_0$, $t_0\in
\Pic^0(C_0)$.  But then $\gamma$ is precisely the elliptic bigrid
corresponding to the sextuple $(a_0,b_0;q_0,t_0;C_0,\phi_0)$, and the
theorem follows for $(m,n)=(3,2)$; the other case follows symmetrically.
\end{proof}

\begin{rem}
In particular, it follows that an elliptic bigrid takes values over a field
$k$ if and only if its parameters can be defined over that field.
\end{rem}

A similar argument shows that for $(m,n)\in\{(l,1),(1,l)\}$ with $l\ge 4$,
the corresponding map from quintuples (forgetting $t$ or $q$ as
appropriate) is birational, with rational image of dimension 7; for $l<4$,
the map is surjective, and the range too small for birationality.  In the
remaining case $(m,n)=(2,2)$, although the space of bigrids is
8-dimensional, the map fails to be birational.

\begin{lem}
The closure of the space of elliptic bigrids of shape $2^2$ is given by the
equation
\[
\chi(\gamma(0,1,0),\gamma(0,2,1),\gamma(1,1,2),\gamma(1,2,1))
=
\chi(\gamma(0,1,1),\gamma(0,2,0),\gamma(1,1,1),\gamma(1,2,2))
\]
where $\chi$ is the cross-ratio function on $\P^1$.  In particular, this
closure is a seven-dimensional rational variety.
\end{lem}

\begin{proof}
The four points $at$, $\tau/aq$, $b/q^2$, $q\tau/bt\in C$ multiply to
$\tau^2/q^2$; it follows that the four points
\[
(\gamma(0,1,0),\gamma(0,1,1)),
(\gamma(0,2,1),\gamma(0,2,0)),
(\gamma(1,1,2),\gamma(1,1,1)),
(\gamma(1,2,1),\gamma(1,2,2))
\in \P^1\times\P^1
\]
satisfy a common bilinear identity (the space of bilinear polynomials and
the space ${\cal L}(D)$ with $\la D\ra=\tau^2/q^2$ are both 4-dimensional).
The given identity follows.  That the image of the space of elliptic
bigrids is seven-dimensional (and thus agrees with the subscheme cut out
by the equation) follows from a direct computation for the $I_1$
degeneration (see below).
\end{proof}

While we have so far been unable to derive a full classification result for
perfect bigrids--that is, a description of its decomposition into
irreducible components--we {\em can} prove the following partial result.

\begin{thm}
If $m,n\ge 2$, then the space of elliptic bigrids is a component of the
space of perfect bigrids; that is, the perfect bigrids in a sufficiently
small neighborhood of a generic elliptic bigrid are all elliptic.
\end{thm}

\begin{proof}
It suffices to show that the tangent space to the space of perfect bigrids
at a generic elliptic bigrid has the same dimension as the space of
elliptic bigrids; that is, that it is seven-dimensional for $(m,n)=(2,2)$,
and otherwise eight-dimensional.  For the small cases $(2,2)$, $(2,3)$,
$(2,4)$, $(3,2)$, $(4,2)$, $(3,3)$, we can verify this as follows (again a
computation facilitated by using the $I_1$ degeneration).  First, choose a
random instance defined over $\Q$, and verify that the tangent space has
the required dimension.  Moreover, the ideal generated by appropriate
minors determines a finite set of primes; reducing modulo a prime outside
that set verifies the result in that characteristic.  Then for each
remaining characteristic, we can choose a random instance over a field of
that characteristic, and again verify the dimension.  This proves the
result in those cases.

We now proceed by induction; let $(m,n)$ be a pair not on the above list
such that the theorem holds for all smaller pairs.  If $m\ge 4$, the two
natural sub-bigrids of shape $(m-1)^n$ overlap in a sub-bigrid of shape
$(m-2)^n$, and thus by induction, for any bigrid in a small neighborhood of the
original bigrid, these sub-bigrids are elliptic.  By assumption, either
$m\ge 5$ or $n\ge 3$, and thus a generic elliptic bigrid of shape
$(m-2)^n$ uniquely determines its parameters. It
follows that the sub-bigrids of shape $(m-1)^n$ have compatible
parameters, and thus join to form an elliptic bigrid of shape $m^n$ as
required.  The proof for $n\ge 4$ is analogous, and thus the theorem
follows.
\end{proof}

In particular, it follows that, just as we needed to generalize Okounkov's
notion of interpolation polynomial to obtain elliptic special functions,
any attempt to further generalize the elliptic theory must necessarily
involve a further generalization of balanced interpolation polynomials;
there is no room in the space of perfect bigrids to add more parameters to
elliptic bigrids.

Now that we have ruled out generalization, it is natural to turn to
specialization.  The space of elliptic bigrids is clearly not closed; it is
therefore of interest to describe the various degenerations of elliptic
bigrids, i.e., to give explicit descriptions of the points in the closure
of the space.  Again, we will be content with a partial result, a
description of the four (over the algebraic closure) particular classes of
degenerations that can produce arbitrarily large regular bigrids.

The key to understanding these degenerations is the curve $C_0\subset
\P^1\times\P^1$ constructed in proving Theorem \ref{thm:birational_bigrid}.
If this curve is smooth, then the bigrid is elliptic, so assume it is
singular; we also assume it is reduced (since the nonreduced case can only
be regular if $m=2$, and corresponds to the case $q^2=1$).  Thus $C_0$ has
isolated singular points, and since it has arithmetic genus one, consists
of a collection of $\P^1$s.

Suppose first that $C_0$ is irreducible, and thus has only one singular point
over $\overline{k}$.  The singular point must lie on the diagonal of
$\P^1\times\P^1$, and by simultaneous ($k$-rational) linear fractional
transformation can be taken to be $(\infty,\infty)$.  We thus obtain a
curve on $\A^1\times \A^1=\A^2$ of the form
\[
a_{20}(x^2+y^2)+a_{11}xy + a_{10}(x+y)+a_{00}.
\]
If $a_{20}=0$, any sequence of points $x_0$, $x_1$, $x_2$, $x_3\in \P^1$
with $(x_0,x_1)$, $(x_1,x_2)$, $(x_2,x_3)\in C$ will contain the same point
twice; it again follows that this case cannot arise from large regular
bigrids.  We may thus assume $a_{20}=1$.  If $a_{11}\ne -2$, then we can
eliminate the linear term; then $a_{00}\ne 0$, lest $(0,0)$ be another
singular point.  There are then two cases, depending on whether the
quadratic term has roots over $k$.  If it does, let $q$ be such that
$q^2+a_{11}q+1=0$.  Then we find that the smooth $k$-rational points on
$C$ are precisely those points of the form
\[
(\frac{u+a_{00}/u}{q-1/q},\frac{qu+a_{00}/qu}{q-1/q})
\]
with $u\in k^*$.  In other words, the smooth points of $C_0$ can be
identified with the multiplicative group in such a way that translation by
$q_0$ is multiplication by $q$.  Similarly, the translation by $t_0$
induced by the bigrid is multiplication by some element $t\in k^*$, as is
the translation by $b/a$.  We thus obtain the following class of (perfect)
bigrids, upon restoring the linear fractional transformation freedom.

\begin{defn}
A degenerate elliptic bigrid of type $I_1$ is a bigrid given by the formula
\[
\gamma(0,i,j) = \phi(a q^jt^{n-i})\quad
\gamma(1,i,j) = \phi(b q^{-j}t^{i-1}),
\]
where $a$, $b$, $q$, $t\in k^*$, and $\phi$ is a $k$-rational degree 2 function
taking the same values at $0$ and $\infty$.
\end{defn}

\begin{drems}
In particular, it is quite straightforward to compute the tangent space to
a generic such bigrid of shape $2^2$, and thus show that that space is
7-dimensional; since the closure of the space of elliptic bigrids contains
this space, its dimension must be at least as large.
\end{drems}

\begin{drems}
The notation $I_1$ is motivated by the fact that the curve $C_0$ is a
degenerate elliptic curve of Kodaira symbol $I_1$, and similarly for the
labels we use for the other untwisted classes.
\end{drems}

We similarly obtain a twisted version of this, when $x^2+a_{11}xy+y^2$ is
irreducible over $k$.  Let $l$ be the splitting field of this polynomial.

\begin{defn}
A degenerate elliptic bigrid of type $I'_1$ is a bigrid given by the formula
\[
\gamma(0,i,j) = \phi(a q^jt^{n-i})\quad
\gamma(1,i,j) = \phi(b q^{-j}t^{i-1}),
\]
where $a$, $b$, $q$, $t$ are in the subgroup of $l^*$ of norm 1 and $\phi$
is a nonconstant function on this subgroup of the form
\[
\phi(x) = \frac{\alpha \Tr(x)+\beta}{\gamma \Tr(x)+\delta},
\]
$\alpha$, $\beta$, $\gamma$, $\delta\in k$.  Note that $\phi$ takes values in
$k\cup\{\infty\}$.
\end{defn}

The final case with a single singular point is when the quadratic term has
a multiple root; i.e., the singular point is a cusp.  In finite
characteristic $p$, we find that the translation map has period $p$ (period
4 in characteristic 2), and thus only leads to bigrids of bounded size; we
may thus assume $k$ of characteristic 0.  We obtain the following.

\begin{defn}
A degenerate elliptic bigrid of type $\II$ is a bigrid given by the formula
\[
\gamma(0,i,j) = \phi(a+jq+(n-i)t)\quad
\gamma(1,i,j) = \phi(b-jq+(i-1)t),
\]
where $a$, $b$, $q$, $t\in k$, and $\phi$ is a degree two function on $k$
of the form
\[
\phi(x) = \frac{\alpha x(\tau-x)+\beta}{\gamma x(\tau-x)+\delta},
\]
with $\alpha$, $\beta$, $\gamma$, $\delta$, $\tau\in k$.
\end{defn}

\begin{rem}
This case can also be thought of as a limit of the $I_1$ case, by taking
all of the points to 1 at comparable rates, including the zeros and poles
of $\phi$.  It corresponds (for sufficiently large $m$, $n$) to a
6-dimensional subvariety of the space of elliptic grids.
\end{rem}

There are two cases to consider when $C_0$ is reducible.  Again, to produce
regular bigrids, each factor must have bidegree at least $(1,1)$, so there
are two factors, both of bidegree $(1,1)$.  Either the individual factors
are symmetrical, and distinct, or the two factors form an orbit under the
exchange map.  Note that since the curve possesses a large number of
rational points, each component must in fact be rational.

In the second case, we find that either there are two points of
intersection (on the diagonal), possibly defined over a quadratic
extension, or there is one rational point of intersection where the two
components are tangent.  (I.e., the curve has Kodaira symbol $I_2$ or
$\III$ respectively.)  For $I_2$, we can transform to obtain the
curve $(x-qy)(y-qx)$, while for $\III$ we can transform to
$(x-y)^2-q^2$.

Now, suppose the curve associated to $j$ translation is reducible, and
consider the curve associated to $i$ translation.  We find that it must
also be reducible, and further that there are compatibility constraints.
Further analysis reveals the following possibilities for bigrids, covering
all of the reducible cases.

Choose an element $\kappa\in k^*$, and construct a group $k^*.2(\kappa)$ on
the set $k^*\times \Z_2$ by the rule
\[
(x,0)(y,0)=(xy,0)\quad
(x,0)(y,1)=(xy,1)\quad
(x,1)(y,0)=(xy,1)\quad
(x,1)(y,1)=(xy/\kappa,0).
\]

\begin{defn}
A degenerate elliptic bigrid of type $I_2$ is a bigrid given by the formula
\[
\gamma(0,i,j) = \phi(a q^j t^{n-i})\quad
\gamma(1,i,j) = \phi(b q^{-j} t^{i-1}),
\]
where $a$, $b$, $q$, $t\in k^*.2(\kappa)$, and $\phi$ is defined by
\[
\phi(x,0) = \phi(\tau/x,1) = \frac{\alpha x+\beta}{\gamma x+\delta},
\]
for some nonconstant, rational, linear fractional transformation and some
$\tau\in k^*$.
\end{defn}

\begin{drems}
Similarly, there is a twisted version with $k^*$ replaced by the norm 1
subgroup of a quadratic extension, and $\phi$ replaced by a nonconstant
$k\cup\{\infty\}$-valued map
\[
\phi(x,0) = \phi(\tau/x,1) = \frac{\alpha x+\overline{\alpha}}{\gamma x+\overline{\gamma}}.
\]
\end{drems}

\begin{drems}
In truth, this is comprised of a total of 8 different cases, depending on
in which coset of $k^*$ each of $a/b$, $q$, $t$ lies.  The result, for
sufficiently large $m$, $n$, is a 7-dimensional space of bigrids, unless
$q$, $t\in k^*$, in which case the space is 6-dimensional.  In particular,
only the latter two cases (in which the full group structure is not being
used) can be obtained by degeneration of $I_1$.  (The corresponding
biorthogonal functions are multivariate analogues of the biorthogonal
functions of \cite{AlSalamWA/IsmailMEH:1994}.)  Apparently the other $I_2$
cases have never been studied, even at the univariate level (where for
$q\notin k^*$ we obtain series such that the term ratios $t_{2m}/t_{2m-1}$
and $t_{2m+1}/t_{2m}$ are {\em different} rational functions of $q^{2m}$).
One can thus view such a series as a sum of two basic hypergeometric series
(by separating out the terms of even index from those of odd index).  For
instance, one such identity (a limiting case of the Frenkel-Turaev
summation identity \cite{FrenkelIB/TuraevVG:1997}) is that
\begin{align}
{}_6\phi_5\biggl(&
\genfrac{}{}{0pt}{}{a,b_1,b_2,b_3,b_4,a^2/b_1b_2b_3b_4}
                   {a/b_1,a/b_2,a/b_3,a/b_4,q^2b_1b_2b_3b_4/a}\biggm|
q^2,q^2\biggr)\notag\\
-
\frac{a^2(1-b_1)(1-b_2)(1-b_3)(1-b_4)}{(a-b_1)(a-b_2)(a-b_3)(a-b_4)}
{}_6\phi_5\biggl(&
\genfrac{}{}{0pt}{}{a,q^2b_1,q^2b_2,q^2b_2,q^2b_3,q^2b_4,a^2/b_1b_2b_3b_4}
           {q^2a/b_1,q^2a/b_2,q^2a/b_3,q^2a/b_4,q^2b_1b_2b_3b_4/a}\biggm|
q^2,q^2\biggr)\notag\\
&{}=
\frac{
(a,a/b_1b_2,a/b_1b_3,a/b_1b_4,a/b_2b_3,a/b_2b_4,a/b_3b_4,a/b_1b_2b_3b_4;q^2)
}{
(a/b_1,a/b_2,a/b_3,a/b_4,a/b_1b_2b_3,a/b_1b_2b_4,a/b_1b_3b_4,a/b_2b_3b_4;q^2)
}
\notag
\end{align}
as long as both sums terminate.  Note that
the above hypergeometric sums are balanced but not {\em quite} well-poised.
\end{drems}

\begin{drems}
This is superficially similar to the degeneration of ordinary interpolation
polynomials denoted IV in \cite{OkounkovA:1998c}, which exists only
for $n=2$; the latter is not, however, a degeneration of our $I_2$, but
rather corresponds to a degenerate curve of Kodaira symbol $I_4$ (which
never produces a regular bigrid).
\end{drems}

Finally, we have the additive analogue of type $I_2$; we define
$k.2(\kappa)$ analogously.

\begin{defn}
A degenerate elliptic bigrid of type $\III$ is a bigrid given by the
formula
\[
\gamma(0,i,j) = \phi(a+jq+(n-i)t)\quad
\gamma(1,i,j) = \phi(b-jq+(i-1)t),
\]
where $a$, $b$, $q$, $t\in k.2(\kappa)$, and $\phi$ is given by
\[
\phi(x,0)=\phi(\tau-x,1)=\frac{\alpha x+\beta}{\gamma x+\delta}
\]
for some nonconstant linear fractional transformation.
\end{defn}

\begin{rem}
Again, aside from degenerations of type $\II$, and aside from a brief
mention of the version of this degeneration for ordinary interpolation
polynomials in \cite{OkounkovA:1998c} (there denoted by III(abc)), this
case does not appear to have been studied.  Each $\III$ case is one
dimension smaller than the corresponding $I_2$ case.
\end{rem}

\section{The trigonometric degeneration}

Of course, the classification of degenerate elliptic bigrids is just a
first step towards an understanding of the various possible limit cases of
our results.  The ``trigonometric'' limit (corresponding to bigrids of type
$I_1$) is of particular interest, as this is where the Koornwinder
polynomials \cite{KoornwinderTH:1992} and Okounkov's interpolation
polynomials \cite{OkounkovA:1998a} (and thus the results
of \cite{bcpoly}) are to be found.  Although we will only be considering
those degenerations necessary to connect our results to those of
\cite{bcpoly}, it is striking how many different limits arise there.  We
thus second Rosengren's call in \cite{RosengrenH:2003b} for a more systematic
classification of limits of elliptic hypergeometric structures.

The only subtle point in extending the algebraic construction to the $I_1$
case is the construction of the difference operator; the point is that the
divisor $\la\tau^2\ra$ on a smooth $C$ degenerates to a divisor that hits
the node.  As a result, a simple condition on poles is not quite enough to
specify $\delta^{(n)}$; we must add a condition comparing the asymptotics
along the two branches of the node.  We thus obtain the following
expression for the degenerate $\delta^{(n)}$, where
$t^{n-1}u_0u_1u_2u_3=\tau^2$, and all parameters lie in $k^*$.
\begin{align}
\delta^{(n)}(\dots x_i\dots;&u_0,u_1,u_2,u_3;t;G_m,\tau)\\
&=
\prod_{1\le i\le n}
\frac{u_0 t^{n-i}(1-u_0x_i/\tau)(1-u_1x_i/\tau)(1-u_2x_i/\tau)(1-u_3x_i/\tau)}
     {x_i(1-x_i^2/\tau)(1-u_1u_0 t^{n-i}/\tau)(1-u_2u_0 t^{n-i}/\tau)(1-u_3u_0 t^{n-i}/\tau)}
\prod_{1\le i<j\le n}
\frac{1-t x_ix_j/\tau}{1-x_ix_j/\tau}.\notag
\end{align}
(Here $G_m$ refers to the multiplicative group scheme, or rather its
natural compactification to a nodal $\P^1$.) Then, if we define
$A^{(n)}_m(u_0;q;G_m,\tau)$ analogously to the elliptic case, we have the
following analogue of Theorem \ref{thm:modular_diff_op}.

\begin{thm}
Define a difference operator on $BC_n(\tau)$-symmetric rational functions
by
\begin{align}
({\tilde D}^{(n)}(u_0,u_1,u_2;q,t;&G_m,\tau)f)(\dots x_i\dots)\\
&=
\prod_{1\le i\le n} (1+\iota_{i,\tau}^*)
\delta^{(n)}(\dots x_i\dots;u_0,u_1,u_2,\frac{\tau^2}{t^{n-1}u_0u_1u_2};t;G_m,\tau)
f(\dots x_i\dots),\notag
\end{align}
where $\iota_\tau(x)=\tau/x$.  Then for any integer $m\ge 0$,
\[
{\tilde D}^{(n)}(u_0,u_1,u_2;q,t;G_m,\tau)
A^{(n)}_m(u_0;q;G_m,\tau/q)
\subset
A^{(n)}_m(u_0;q;G_m,\tau).
\]
Moreover,
\begin{align}
({\tilde D}^{(n+k)}(u_0,u_1,u_2;q,t;G_m,\tau)f)(z_1,\dots, z_n,&u_0,\dots, t^{k-1}u_0)\\
&=
({\tilde D}^{(n)}(u_0,u_1,u_2;q,t;G_m,\tau)f(\underline{\quad},u_0,\dots, t^{k-1}u_0))
(z_1,\dots, z_n).\notag
\end{align}
\end{thm}

\begin{proof}
The only way in which the proof differs is that we need to show that the
image function is smooth at the node.  This is true by assumption for $f$,
as is the condition (by $\iota_\tau$-invariance) that either branch of the
node gives the same limit.  But then multiplying by $\delta^{(n)}$ gives a
function with linear growth at the node, in opposite directions along the
two branches.  Thus symmetrizing kills off this linear growth, giving
smoothness as required.
\end{proof}

Similarly, the $\Delta^0$, $\Delta$ and diagonal binomial coefficient
symbols are all straightforward to define.  We find, if $u_0u_1\cdots
u_{2m-1}=\tau^m$,
\[
\omega(x{:}u_0,\dots, u_{2m-1};G_m,\tau)
=
\prod_{0\le r<2m} \frac{1-x/u_r}{1-u_rx/\tau},
\]
and thus if $v_0v_1\cdots v_{2m-1} = q^m a^m$,
\begin{align}
\Delta^0_\lambda(a|v_0,\dots, v_{2m-1};q,t;G_m)
&=
\prod_{(i,j)\in \lambda}
\prod_{0\le r<2m} \frac{1-q^{1-j}t^{i-1}/v_r}{1-v_r q^{-j}t^{i-1}/a}\\
&=
\frac{C^0_\lambda(v_0,\dots, v_{2m-1};q,t)}
     {C^0_\lambda(qa/v_0,\dots, qa/v_{2m-1};q,t)}
\end{align}
Similarly, if $v_0v_1\cdots v_{2m-1}=(t/q) (qa)^{m-1}$
\[
\Delta_\lambda(a|v_0,\dots, v_{2m-1};q,t;G_m)
=
\frac{
q^{|\lambda|} t^{2n(\lambda)}
C^0_{2\lambda^2}(qa;q,t)C^0_\lambda(v_0,\dots, v_{2m-1};q,t)}
{C^-_\lambda(q,t;q,t)C^+_\lambda(a,qa/t;q,t)C^0_\lambda(qa/v_0,\dots, qa/v_{2m-1};q,t)}
\]
Finally,
\[
\binomE{\lambda}{\lambda}_{[a,b];q,t;G_m}
=
\frac{b^{|\lambda|}C^+_\lambda(a;q,t)C^0_\lambda(1/b,qa/b;q,t)}
     {C^+_\lambda(a/b;q,t) C^0_\lambda(b,qa;q,t)}
\]
(Here we use the $C^*$ notations of \cite{bcpoly}, which are simply the
limits $p\to 0$ of the elliptic notations.)

It is then trivial to state the analogues of the elliptic results for this
most general trigonometric case; simply replace $C$ by $G_m$ in each
identity of Section \ref{sec:modular}.  This in particular gives
multivariate analogues of the ${}_{10}\phi_9$ biorthogonal rational
functions of Rahman \cite{RahmanM:1986} (or, more precisely, since we are
considering only discrete biorthogonality here, they are multivariate
analogues of the biorthogonal rational functions of Wilson
\cite{WilsonJA:1991}; in any event, continuous biorthogonality {\em is}
preserved in the trigonometric limit).  The new thing that appears is that
we can take limits as various parameters approach 0 or infinity, i.e., the
node of the curve.  There are, it turns out, three different limits of
interpolation functions and binomial coefficients that need to be
considered for a full understanding of \cite{bcpoly} from our current
perspective.

Of course, one of these limits is simply the ordinary interpolation polynomials
$\bar{P}^{*(n)}_\lambda(;q,t,s)$.  More precisely, we have the following
limit theorem.

\begin{thm}
The interpolation functions and ordinary interpolation polynomials are related by
the identity
\[
\bar{P}^{*(n)}_\lambda(\dots x_i\dots;q,t,a)
=
\frac{t^{2n(\lambda)} C^0_\lambda(t^n,t^{n-1} av,a/v;q,t)}
     {(-t^{n-1} a)^{|\lambda|} q^{n(\lambda')} C^-_\lambda(t;q,t)}
\lim_{b\to \infty} R^{*(n)}_\lambda(\dots x_i\dots;a,b(v);q,t;G_m,\la 1\ra^2).
\]
Moreover, replacing the limit $b\to\infty$ by $b\to 0$ gives the same result.
\end{thm}

\begin{proof}
Fix $a$, $q$, $t$, and let
\[
f^{(n)}_\lambda(\dots x_i\dots;v)
\]
denote the corresponding right-hand side.  We claim that for each $n$,
$f^{(n)}_\lambda(;v)$ is a well-defined $BC_n$-symmetric Laurent polynomial
of degree at most $|\lambda|$.  Indeed, from the branching rule, we find
\[
f^{(n)}_\lambda(\dots x_i\dots;v)
=
\sum_\kappa
\frac{C^0_{\lambda/\kappa}(t^{n-1} a x_n^{\pm 1};q,t)}
     {C^0_{\lambda/\kappa}(t^{n-1} a v^{\pm 1};q,t)}
\lim_{b\to\infty}
\obinomE{\lambda}{\kappa}_{[t^{n-1}a/b,t](qa/tb,qt^{n-1}v/b,t^{n-1}a/v);q,t;G_m}
f^{(n-1)}_\kappa(\dots x_i\dots;vt),
\]
so by induction it suffices to show that the binomial coefficient has a
well-defined limit.  This then follows readily from the algebraic analogue
of Corollary \ref{cor:binom_spec_recur}.

Thus $f^{(n)}_\lambda(;v)$ is a $BC_n$-symmetric Laurent polynomial of at
most the correct degree, vanishing at points of the form $\dots a
t^{n-i}q^{\mu_i}\dots$ for $\mu\not\subset\lambda$, and is thus a multiple
of the ordinary interpolation polynomial as required.  That this multiple is 1
follows by comparing values at the point $\dots v t^{n-i}\dots$.

The argument for $b\to 0$ is analogous.
\end{proof}

\begin{rems}
One can also use the point $\dots a t^{n-i} q^{\lambda_i}\dots$ to set the
normalization; the resulting computation is somewhat more complicated, but
there is the slight advantage that the proof in \cite{bcpoly} of the
latter normalization does not use any properties of Macdonald polynomials.
As a result, we see that we can deduce such properties from the elliptic
theory.
\end{rems}

\begin{rems}
It ought to be possible to prove directly, using only the vanishing
conditions, that the interpolation polynomials are limits of interpolation
functions.  More precisely, given any partial order on partitions refining
the inclusion partial order, the complementary vanishing conditions define
a corresponding filtration on the space of $BC_n$-symmetric rational
functions; the main question is then to show that one of these filtrations
has the correct limit as $b\to 0,\infty$.  (Two cases of particular
interest are the dominance order, for which the limiting filtration is
simply that given by dominance of monomials, so in particular is
independent of $q$ and $t$, and the inclusion order itself, for which the
limiting filtration is, by results of \cite{bcpoly}, that induced by the
Macdonald polynomials.)
\end{rems}

Plugging into the respective definitions of binomial coefficients gives the
following result (in which the right-hand side uses the notation of Section
4 of \cite{bcpoly}).

\begin{cor}
We have the following identity of binomial coefficients.
\[
\lim_{b\to 0,\infty}
\binomE{\lambda}{\mu}_{[a,b];q,t;G_m}
=
\frac{(-1)^{|\mu|} t^{n(\mu)} C^+_\mu(a;q,t)}
     {q^{n(\mu')} C^0_\mu(qa;q,t)}
\binomQ{\lambda}{\mu}_{q,t,\sqrt{a}}
\]
\end{cor}

Now, it then follows that we have a corresponding limit for the inverse
binomial coefficients of \cite{bcpoly}.

\begin{cor}
We have the following identity of binomial coefficients.
\[
\lim_{b\to 0,\infty}
\binomE{\lambda}{\mu}_{[a/b,1/b];q,t;G_m}
=
\frac{(-1)^{|\lambda|} q^{n(\lambda')} C^0_\lambda(qa;q,t)}
     {t^{n(\lambda)} C^+_\lambda(a;q,t)}
\binomI{\lambda}{\mu}_{q,t,\sqrt{a}}
\]
\end{cor}

Now, the limiting bigrid on the left is in fact a regular bigrid of type
$I_2$; we thus obtain a direct description of the inverse binomial
coefficients in terms of vanishing conditions.  In particular, this
explains what happened to the evaluation symmetry condition when extending
from shifted Macdonald polynomials \cite{SahiS:1996,KnopF:1997} to ordinary
interpolation polynomials \cite{OkounkovA:1998a}: rather than relate two
functions of the same type, the symmetry relates functions from two
different degenerations.  The resulting family of degenerate interpolation
functions is simply the limit as $\tau\to 0$ or $\infty$ of the general
trigonometric interpolation functions.  In particular, for these functions,
the generalized eigenvalue problem simplifies to an eigenvalue problem,
relative to the family of difference operators (parametrized by $c$ with
$a$, $b$, $q$, $t$ fixed):
\begin{align}
\sum_{I\subset \{1,2,\dots, n\}}&
\frac{t^{|I|(|I|-1)/2}
\prod_{i\in I}    (x_i(1-a/x_i)(1-b/x_i))
\prod_{i\notin I} (c(1-x_i/c)(1-t^{n-1}ab/cx_i))}
{\prod_{1\le i\le n} (c(1-a t^{n-i}/c)(1-b t^{n-i}/c))}\notag\\
&\prod_{i\in I,j\notin I} \frac{x_i-t x_j}{x_i-x_j}
f(\dots q^{-[i\in I]} x_i\dots),
\end{align}
with eigenvalues
\[
\frac{q^{-|\lambda|}C^0_\lambda(aqt^{n-1}/c,qc/b;q,t)}
     {C^0_\lambda(t^{n-1}a/c,c/b;q,t;p)}.
\]
(Replace $c$ by $\tau/c$ in the difference equation for trigonometric
interpolation functions, then take $\tau\to 0,\infty$.  Note that the
further limit $c\to 0,\infty$ makes this the identity operator.)

If we attempt to obtain the bulk branching rule for ordinary interpolation
polynomials (Theorem 3.9 of \cite{bcpoly}) as a limit of our theory, we
find a third limit of binomial coefficients appearing.

\begin{thm}
We have the following identity of binomial coefficients.
\[
\lim_{a\to 0,\infty}
\binomE{\lambda}{\mu}_{[a,b];q,t;G_m}
=
\frac{b^{|\mu|} t^{n(\mu)} C^-_\lambda(t;q,t) C^0_\mu(1/b;q,t)}
     {t^{n(\lambda)} C^-_\mu(t;q,t) C^0_\lambda(b;q,t)}
P_{\lambda/\mu}([\frac{1-b^k}{1-t^k}];q,t)
\]
\end{thm}

\begin{proof}
Take the appropriate limit of the bulk branching rule for ordinary
interpolation functions to obtain a bulk branching rule for ordinary
interpolation polynomials.  Taking a further limit produces the Macdonald
polynomials (i.e., polynomials satisfying Macdonald's difference equation);
the result then follows essentially by the definition of $P_{\lambda/\mu}$.
\end{proof}

\begin{rems}
If we take a further limit $b\to 0$, the binomial coefficient can be
expressed in shifted Macdonald polynomials; we thus obtain Theorem I.4 of
\cite{GarsiaAM/HaimanM/TeslerG:1999}, which expresses the corresponding
binomial coefficients via the plethystic expression
\[
P_{\lambda/\mu}([\frac{1}{1-t^k}];q,t).
\]
\end{rems}

\begin{rems}
Of our three types of degenerate trigonometric binomial coefficients, this
is the only one that includes the four important special cases
\[
\binomE{\lambda}{\mu}_{[a,1];q,t},
\binomE{\lambda}{\mu}_{[a,1/q];q,t},
\binomE{\lambda}{\mu}_{[a,t];q,t},
\binomE{\lambda}{\mu}_{[a,q/t];q,t}.
\]
That the first three specializations of skew Macdonald polynomials have
nice factorizations is well known (the first from the fact that
$P_{\lambda/\mu}$ has positive degree, so constant term 0, the second and
third from the coefficients of the Pieri identities and branching rule).
That
\[
P_{\lambda/\mu}([\frac{1-(q/t)^k}{1-t^k}];q,t)
\]
can be expressed as a product of binomials does not appear to have been
observed before.
\end{rems}

These also correspond to a degenerate elliptic bigrid of type $I_2$, this
time with $a/b\not\in k^*$.  The difference operators are the same as for
the previous $I_2$ case, except that now the family has $a$, $c$, $q$, $t$
fixed and $b$ variable; in particular, we return to a generalized
eigenvalue problem scenario.

We thus find that all of the ordinary interpolation polynomial identities
of \cite{bcpoly} are indeed limits of identities satisfied by our
interpolation functions.  Similarly, the Koornwinder polynomials are
limiting cases of our biorthogonal functions:
\[
\lim_{u_0\to 0,\infty}
\tilde{R}^{(n)}(\dots
x_i\dots;t_0{:}t_1,t_2,t_3;u_0,q/t^{2n-2}t_0t_1t_2t_3u_0;q,t;G_m,\la 1\ra^2)
=
\frac{K^{(n)}_\lambda(\dots x_i\dots;q,t;t_0,t_1,t_2,t_3)}
     {K^{(n)}_\lambda(\dots t_0 t^{n-i}\dots;q,t;t_0,t_1,t_2,t_3)},
\]
where the normalization factor on the right can be deduced by taking
coefficients of $P^{*(n)}_\lambda(;q,t,t_0)$ on both sides.  In particular,
we note that Theorems \ref{thm:conn_biorth1}, \ref{thm:quasi-branch}, and
\ref{thm:quasi-Pieri} give rise to new results at the Koornwinder level,
generalizing similar results from \cite{bcpoly} in which the relevant
trigonometric binomial coefficient has second parameter $\in \{1/q,t\}$
(i.e., corresponds to a difference or integral operator).

\medskip

As an indication of some of the subtle issues that can arise when
degenerating the elliptic theory, we consider the following ``rational''
limit.  Here we take $q,t\to 1$ along a one-parameter family in such a way
that the corresponding maps of tangent spaces have ratio $\alpha$ at the
limit point.  In other words, we blow up $\Pic^0(C)^2$ at the point $(1,1)$
and consider the interpolation function at a point on the
exceptional divisor.  Note that the result is essentially independent of
every parameter {\em except} $\alpha$, including the choice of $C$ itself.

\begin{thm}
Let $(C,\tau)$ be an arbitrary genus 1 hyperelliptic curve over a field of
characteristic 0.  Then we have the following identity of interpolation
functions, for $a$, $b$, $v\in C$, $\alpha\in k^*$ generic.
\[
\lim_{\substack{q,t\to 1\\dq/dt\to \alpha}}
R^{*(n)}_\lambda(\dots z_i\dots;a,b(v);q,t;C,\tau)
=
\frac{J_\lambda(\dots\chi(z_i,a,v,b;C,\tau)\dots;\alpha)}
     {J_\lambda(\dots 1\dots;\alpha)},
\]
where $J_\lambda(;\alpha)$ is a Jack polynomial.
\end{thm}

\begin{proof}
It suffices to prove the result for $(C,\tau)$ of the form $(\C^*/p,\la
1\ra^2)$, in which case taking a limit in the branching rule gives the
desired result.
\end{proof}

There are several points to observe here.  The first is that, unlike in the
case of bigrids, we obtained a nontrivial degeneration without degenerating
the curve itself.  In fact, while the corresponding elliptic bigrid has a
well-defined limit, that limit is very far from being regular; indeed, the
limiting bigrid is actually independent of $i$ and $j$, having two constant
values independent of $\alpha$.  In other words, the map from elliptic
bigrids to interpolation polynomials is merely birational; a full
classification of degenerate cases of the one is insufficient to give a
classification the other.  The second point is that, compared to the
Macdonald polynomial case, it is much harder to obtain a useful limit from
the difference equation.  Indeed, the difference operator actually becomes
the identity in the limit, and even the generalized eigenvalue problem,
while it has a nontrivial limit when properly scaled, simply gives (after a
nontrivial calculation) the equation
\[
\sum_i \frac{\chi(z_i,a,v,b;C,\tau)\partial
  f(z_i)}{f(z_i)\partial\chi(z_i,a,v,b;C,\tau)} = |\lambda|,
\]
corresponding to the fact that $J_\lambda$ is homogeneous of degree
$|\lambda|$.  The third point is that we can similarly obtain limits of
binomial coefficients analogous to the three trigonometric limits we
obtained above; we find again that those limits are essentially independent
of those parameters not approaching 1.  We finally note that using the
binomial formula, we can obtain an analogous limit for the biorthogonal
functions; if $t_0t_1$, $t_2t_3$, $u_0u_1$ all tend to 1 simultaneously
with $q$, $t$, then the biorthogonal function becomes a multivariate Jacobi
polynomial, up to an analogous change of variables.  In the case $t_0\sim
1$, $t_2\sim -1$, the continuous biorthogonality density of \cite{xforms}
easily simplifies to the usual Jacobi ensemble (just as with the analogous
limit for the Koornwinder polynomials); more generally, if $|t_0|=1=|t_2|$,
we can at least formally obtain the appropriate ensemble, using the remark
following Lemma 3.3 of \cite{xforms}.

\section{Open problems}

It seems fitting to conclude by discussing several open problems arising
from the above theory.  We have already discussed partial results for one
such problem in the previous section, namely the question of classifying
degenerations of interpolation and biorthogonal functions; along these
lines, we also note the question of finding a direct proof that the
interpolation polynomials are limits of the interpolation functions.  For
the remaining problems we consider, we have even less in the way of
explicit results or even conjectures.

The most popular approach to the theory of Koornwinder polynomials involves
the so-called ``double affine Hecke algebra'' \cite{SahiS:1999}, along with
the associated representation of the affine Hecke algebra.  A major open
question, therefore, is how to generalize this theory to the elliptic
level, both to understand how the Koornwinder theory relates to the
elliptic theory, and ideally to allow the elliptic theory to be extended to
other root systems.  Via evaluation symmetry, this should be related to the
question of finding a general Pieri identity for the biorthogonal
functions.  We note the following partial result.

\begin{prop}
For otherwise generic parameters satisfying $t^{2n-2}
t_0t_1t_2t_3u_0u_1=q\tau^3$, and generic $u,v\in C$, let $c_{\lambda\kappa}$
be the coefficient of
\[
\tilde{R}^{(n)}_\kappa(\dots z_i\dots;t_0{:}t_1,t_2,t_3;u_0q,u_1/q;q,t;C,\tau)
\]
in the product
\[
\prod_{1\le i\le n}
\chi(z_i,u,v t^{n-i},u_0;C,\tau)
\tilde{R}^{(n)}_\lambda(\dots z_i\dots;t_0{:}t_1,t_2,t_3;u_0,u_1;C,\tau).
\]
Then $c_{\lambda\kappa}=0$ unless $\lambda\subset 1^n+\kappa\subset
2^n+\lambda$.
\end{prop}

\begin{proof}
That $\kappa\subset 1^n+\lambda$ follows from the Pieri identity for
interpolation functions and the fact that the biorthogonal functions are
triangular in the appropriate basis of interpolation functions.  That
$\lambda\subset 1^n+\kappa$ follows by the fact that for
$t^{n-1}q^m t_0t_1=\tau$, $c_{\lambda\kappa}$ can be computed via
biorthogonality, and thus satisfies a symmetry switching $\lambda$ and
$\kappa$.
\end{proof}

In the cases $u\in \{t_0,t_1,t_2,t_3,t_0/q,t_1/q,t_2/q,t_3/q\}$, we can
combine the quasi-Pieri identity, Theorem \ref{thm:quasi-Pieri}, with the
connection coefficient identity, to express the above coefficient as a sum
over binomial coefficients with $b=1/q$, which is thus explicitly supported
on the stated set of partitions.  The question is thus whether a
corresponding formula exists for general $u$.  (An expression in a double
sum is straightforward, but it is then nonobvious that the coefficients
have the support they do.)  This would give a family difference operators
via evaluation symmetry, which (by consideration of the Koornwinder case)
would give the analogue of the degree 1 subspace of center of the affine
Hecke algebra.  (See also the remark following Theorem 8.9 of
\cite{xforms}.)

Similar considerations hold for a general branching rule.

\begin{prop}
For otherwise generic parameters satisfying $t^{2n-2}
t_0t_1t_2t_3u_0u_1=q\tau^3$, and generic $u\in C$, let $c_{\lambda\kappa}$
be the coefficient of
\[
\tilde{R}^{(n-1)}_\kappa(\dots
z_i\dots;t_0{:}t_1,t_2,t_3;u_0,u_1t^2;q,t;C,\tau)
\]
in the specialization
\[
\tilde{R}^{(n)}_\lambda(\dots z_i\dots,u;t_0{:}t_1,t_2,t_3;u_0,u_1;C,\tau).
\]
Then $c_{\lambda\kappa}=0$ unless $\kappa\prec'_2 \lambda$.
\end{prop}

\begin{proof}
Using the Cauchy identity, Theorem \ref{thm:Cauchy_biorth}, one can
express the coefficients of the branching rule in terms of the
coefficients of the Pieri identity (and vice versa); the stated vanishing
property follows.
\end{proof}

In particular, we see that it is hopeless to expect anything better than a
sum for these coefficients, as even in the case $n=1$ we obtain a
nontrivial hypergeometric sum.  We also note that such an expression for
the branching rule would be a new result even for the Koornwinder case.

\medskip
Another question, which is less important, but likely to have interesting
consequences, concerns the classification of regular perfect bigrids.  The
four non-elliptic components we have identified above all correspond to
important special cases of the interpolation functions, and thus of the
corresponding hypergeometric identities.  It thus seems likely that other
components would give similarly interesting special cases.  One such
component would appear to correspond to the following.

\begin{conj}
Any bigrid such that $\gamma(1,i,j)=\gamma(0,i+1,j)$ whenever both sides
are defined is perfect.
\end{conj}

Note that this would give a $((n+1)(m+1)-2)$-dimensional space of perfect
grids, and would thus be larger even than the $m(n+1)$-dimensional Cauchy
case.  When $m=n=2$ (in which case the conjecture is straightforward to
verify), we obtain a 7-dimensional space, distinct from the 7-dimensional
space of elliptic bigrids; in contrast, a bigrid of shape $2^2$ coming from
any of the four known cases is necessarily elliptic.  In the elliptic case,
this corresponds to binomial coefficients with $b=1/t$, and thus should
correspond in some sense to an inverse of the integral operator of
\cite{xforms}.

The other conjectured component we have found (by computing tangent spaces
to random elliptic bigrids over small finite fields) is the following.

\begin{conj}
Let $\tau$ be a linear fractional transformation of order 2.  Then any
bigrid satisfying
\[
\gamma(0,i,j)=\tau(\gamma(1,i,j+1))=\tau(\gamma(0,i+1,j+1))
\]
whenever both sides are defined is perfect.
\end{conj}

\medskip
Finally, we observe that we have only extended about half of the results of
\cite{bcpoly} to the elliptic case.  The remaining results fall into two
main classes.  The first of these is the construction of families of
symmetric functions algebraically continuing interpolation and Koornwinder
polynomials in the number of variables (via the quantity $t^n$).  Although
it seems unlikely that these results extend cleanly to the elliptic level
(given the absence of an appropriate theory of symmetric functions), there
should still be some weak version; after all, the connection coefficients
that arise above depend on the number of variables only via $t^n$, and
similarly for the bulk branching rule.  The other class of results concerns
the ``vanishing integrals'' conjectured in \cite{bcpoly}, which generalize
quadratic transformations of univariate hypergeometric series; elliptic
analogues of these would likely be interesting even at the univariate
level.  Note in particular that there are two potential kinds of quadratic
transformation at the elliptic level, depending on whether the relevant
isogeny is a multiplication map (squaring, in our notation), or simply an
isogeny of degree 2.

\bibliographystyle{plain}

\end{document}